\documentclass[11pt,a4paper,sunlimits]{article}

\usepackage[francais]{babel}
\usepackage{indentfirst}
\usepackage{amsmath}
\usepackage{amsthm}
\usepackage{amsfonts}
\usepackage{amssymb}
\usepackage{amscd}
\usepackage{textcomp}
\usepackage{graphicx}
\usepackage[usenames, dvipsnames]{color}
\usepackage[colorlinks]{hyperref}

\usepackage[all]{xy}
\newtheorem{theo}{Th\' eor\`eme}[subsection]
\newtheorem{lem}{Lemme}[subsection]
\newtheorem{prop}{Proposition}[subsection]
\newtheorem{Def}{D\' efinition}[subsection]
\newtheorem{cor}{Corollaire}[subsection]

\newtheorem{rem}{Remarque}[subsection]
\newenvironment{dem}{\noindent $ Preuve:$}%
{\unskip\hfill\null\nobreak\hfill\carre\vskip1em\par}
\newcommand{\carre}{\rule{1 ex}{1 ex}}

\newcommand{\al}{\alpha}
\newcommand{\be}{\beta}
\newcommand{\De}{\Delta}
\newcommand{\de}{\delta}
\newcommand{\ep}{\epsilon}

\newcommand{\ga}{\gamma}
\newcommand{\ka}{\kappa}
\newcommand{\la}{\lambda}
\newcommand{\om}{\omega}
\def\Om{\Omega}
\def\vi{\varpi}

\newcommand{\te}{\theta}

\def\a{\mathfrak{a}}
\def\b{\mathfrak{b}}
\def\ce{\mathfrak{c}}

\def\g{\mathfrak{g}}
\def\h{\mathfrak{h}}

\def\l{\mathfrak{l}}
\def\m{\mathfrak{m}}
\def\n{\mathfrak{n}}

\def\q{\mathfrak{q}}
\def\qr{\mathfrak{q}^{reg}}
\def\r{\mathfrak{r}}

\def\z{\mathfrak{z}}

\newcommand{\R}{\mathbb{R}}
\newcommand{\C}{\mathbb{C}}
\newcommand{\N}{\mathbb{N}}
\newcommand{\Z}{\mathbb{Z}}

\def\me{\medskip}
\def\no{\noindent}
\newcommand{\qq}{\quad}
\newcommand{\dis}{\displaystyle}

\def\mah{\left (\begin{array}{c|c} A&0\\ \hline
0&B\end{array}\right )}
\def\maq{\left (\begin{array}{c|c}
0&Y\\ \hline Z&0\end{array}\right )}
\def\com{\left(
\begin{array}{c|c}
0&\begin{array}{cc} \tau&-\te\\ \te&\tau\end{array}\\ \hline
\begin{array}{cc} \tau&-\te\\ \te&\tau\end{array}&0
\end{array}
\right)}
\def\app{\a_{+,+}}
\def\apm{\a_{+,-}}
\def\amp{\a_{-,+}}
\def\amm{\a_{-,-}}
\def\ad{\a_{2}}
\def\are{\a^{reg}}

\def\sec{sous-espace de Cartan $\space$}
\def\w{W_{H_\C}(\a_\C)}
\def\wg{W(\g_\C,\a_\C)}

\def\mul{m_\alpha}
\def\ra{\Delta(\mathfrak{g}_{\mathbb{C}},\mathfrak{a}_{\mathbb{C}})}
\def\rap{\Delta(\mathfrak{g}_{\mathbb{C}},\mathfrak{a}_{\mathbb{C}})^+}

\def\dq{\partial (Q)}
\def\ds{\partial (S)}
\def\rdp{\Delta_{\mathfrak{a}}(\partial P)}
\def\rdq{\Delta_{\mathfrak{a}}(\partial Q)}
\def\rds{\Delta_{\mathfrak{a}}(\partial S)}
\def\el{L_\alpha}
\def\elu{L_{\alpha_{1}}}
\def\eld{L_{\alpha_{2}}}
\def\L1loc{L_{loc}^1(\mathfrak{q})^H}
\def\Lu{L_{loc}^1(\mathfrak{q})^H}
\def\cqh{\mathbb{C}[\mathfrak{q}_\mathbb{C}]^{H_\mathbb{C}}}
\def\zt{\z_3\cap\q}
\def\mq{\m\cap\q}
\def\ztp{(\zt)_+}
\def\fr{\frac{\partial}}

\newcommand{\p}{^\backprime}
\newcommand{\D}{\mathcal{D}}
\newcommand{\im}{\imath}

\def\ad{\textrm{ad}}
\def\Ad{\textrm{Ad}}
\def\Mf{\mathcal{M}_H(f)}
\def\Mfm{\mathcal{M}f_\mathfrak{m}}
\def\Mf2{\mathcal{M}f_2}
\def\E{\mathcal{E}_{sing}}

\begin{document}
\author{Pascale Harinck\footnote{Ecole Polytechnique, CMLS- CNRS UMR 7640, Route de Saclay 91128 Palaiseau C\' edex, harinck@math.polytechnique.fr  }\and Nicolas Jacquet\footnote{Lyc\' ee Jean-Baptiste Corot, 9 Place Davout,  91600  Savigny-sur-Orge, nclsjacquet@gmail.com  }}
\title{Distributions propres invariantes sur la paire sym\' etrique  $\big(\mathfrak{gl}(4,\R), \mathfrak{gl}(2,\R)\times \mathfrak{gl}(2,\R)\big)$}
\date{}
\maketitle
\begin{abstract} We study orbital integrals and invariant eigendistributions for the symmetric pair $(\g,\h)=\big(\mathfrak{gl}(4,\R),\mathfrak{gl}(2,\R)\times \mathfrak{gl}(2,\R)\big)$. Let  $\q=\g/\h$ and  let $\mathcal{N}$ be the set of nilpotents of $\q$. We first obtain an asymptotic behavior of orbital integrals around nonzero semisimple  elements of $\q$. We study eigendistributions around such elements and give  an explicit  basis of  eigendistributions on $\q-\mathcal{N}$  given by a locally integrable function on $\q-\mathcal{N}$.

\no {\it Keywords~: symmetric pair, orbital integral, invariant eigendistribution}
\end{abstract}

\maketitle
\section*{Introduction}

Dans  cet article, nous \' etudions  les int\' egrales orbitales  des fonctions de classe ${\mathcal C}^\infty$ \`a support compact et  les distributions propres invariantes  pour la paire sym\' etrique  $(\mathfrak{gl}(4,\Bbb R),\mathfrak{gl}(2,\Bbb R)\times\mathfrak{gl}(2,\Bbb R))$. 

 De  tels objets sont parfaitement d\' ecrits  pour les alg\`ebres de Lie  r\' eductives $\g_1$ ( correspondant \`a la paire sym\' etrique $(\g_1\times\g_1,diagonale)$)  par les travaux d'Harish-Chandra puis de A. Bouaziz, et  pour  les paires sym\' etriques de rang $1$ \`a travers les travaux de nombreux auteurs dont   P.D.~M\' eth\' e, A.~Tengstrand, J.~Faraut et G.van Dijk. Par ailleurs, 
J.~Sekiguchi a d\' etermin\' e la liste des paires sym\' etriques pour lesquelles il n'existait pas de distributions propres invariantes \`a support dit singulier c'est-\`a-dire dans le compl\' ementaire  dans $\q$ de l'ouvert dense $\qr$ des \' el\' ements r\' eguliers de $\q$. 

Hormis ces situations, les propri\' et\' es des  int\'egrales orbitales et des distributions propres invariantes 
sont peu connues dans le cadre g\' en\' eral des paires sym\' etriques quelconques.
Le travail  pr\' esent\' e ici  donne un  premier exemple d'\'etude pr\' ecise et explicite  de tels objets pour une paire sym\'etrique de rang sup\' erieure \`a $2$.\me

On  consid\`ere le groupe $G=GL(4,\R)$ et son alg\`ebre de Lie  $\g =\mathfrak{gl}(4,\Bbb R)$ munis de
l'involution $\sigma$ d\' efinie par $\displaystyle \sigma(X)=\left
(\begin{array}{cc}
I_{2}&0\\0&-I_{2}
\end{array}\right )X
\left (
\begin{array}{cc}I_{2}&0\\0&-I_{2}
\end{array}\right ),$
 o\`u $I_{2}$ est la matrice identit\' e $2\times 2$. 
Soit
 $\g=\h\oplus\q$ la d\' ecomposition     de $\g$ relative en sous-espaces propres pour $\sigma$ associ\' es respectivement  \`a $+1$ et $-1$ et $H=GL(2,\R)\times GL(2,\R)$ le sous-groupe de $G$ des points fix\' es par $\sigma$.
 
L'int\' egrale orbitale    d'une fonction
$f\in\mathcal{C}^\infty_c(\mathfrak{q})$ est la moyenne de $f$ sur les $H$-orbites des \' el\' ements semi-simples r\' eguliers de $\q$. Nous pr\' ecisons tout d'abord  le comportement des int\' egrales  orbitales au voisinage de tout point semi-simple non nul.  On note $\mathcal{N}$ le c\^one des \' el\' ements nilpotents de $\q$ et on fixe un caract\`ere $\chi$ r\' egulier de l'alg\`ebre $S(\q)^H$ des op\' erateurs $H$-invariants \`a coefficients constants sur $\q$. Nous obtenons des conditions n\' ecessaires  et  suffisantes pour qu'une fonction localement int\' egrable sur  ${\mathcal U}=\q-\mathcal{N}$ d\' efinisse une distribution  $T_F$  invariante et propre pour le caract\`ere $\chi$ sur ${\mathcal U}$, c'est-\`a-dire sastisfaisant, pour tout $P\in S(\q)^H$, la relation $\partial(P)T_F=\chi(P)T_F$ sur  ${\mathcal U}$. Nous d\' eduisons une base
    de distributions   invariantes, propres  pour le caract\`ere $\chi$ sur ${\mathcal U}$ donn\' ees par une fonction
     localement int\'egrable sur ${\mathcal U}$ et nous
    ferons quelques commentaires sur l'extension de ces solutions sur $\q$.

 La paire sym\' etrique $(\g,\h)$ figure dans la liste de J.~Sekiguchi et correspond \`a   la situation tangente  de l'espace sym\' etrique $GL(4,\R)/GL(2,\R)\times GL(2,\R)$. Des r\' esultats partiels concernant int\' egrales orbitales et  distributions sph\' eriques sur $GL(4,\R)/GL(2,\R)\times GL(2,\R)$ ont \' et\' e annonc\' es par S.~Kato et S.~Aoki dans (\cite{ao1}), mais aucune preuve n'est explicitement donn\' ee  dans leur article ce qui ne nous a    pas permis de comparer  nos r\' esultats. \me
 
 Pr\'ecisons le contenu de cet article. 
 
 La premi\`ere partie  est consacr\' ee \`a l'\' etude  des orbites sous $H$ des \' el\' ements 
    semi-simples  et des sous-espaces de Cartan de $\q$ (c'est-\`a-dire form\' e d'\' el\' ements semi-simples de $\q$, ab\'
elien et maximal pour ces propri\' et\' es) . Les $H$-orbites semi-simples de $\q$ sont uniquement d\' etermin\'
ees par les polyn\^omes $Q(X)=tr(X^2)/2$ et $S(X)=det(X)$
qui forment une base de  l'alg\`ebre $\C[\q]^H$ des polyn\^omes
$H$-invariants sur $\q$, ou encore par les valeurs propres $u(X) $
et $v(X)$ de $YZ$ si $X=\left (
\begin{array}{cc}0&Y\\Z&0
\end{array}\right )$. Les fonctions $u$ et $v$ satisfont les relations $u+v=Q$ et $uv=S$.

La paire sym\' etrique $(\g ,\h)$ est de rang $2$ et poss\`ede $4$
classes de conjugaison sous $H$ de sous-espaces de Cartan
selon la nature des
valeurs propres $u(X)$ et $v(X)$: $\a_{++}$,  $\a_{+-}$ et $\a_{--}$
sur lesquels les fonctions $u$ et $v$ prennent des valeurs r\'
eelles de signe constant et $\a_2$ sur lequel $u$ et $v$ sont
complexes conjugu\' ees.   L'ensemble  $\q^{reg}$  est form\' e des
\' el\' ements semi-simples $X$ dont  les valeurs propres sont deux
\`a deux distinctes ce qui \' equivaut \`a la condition
$u(X)v(X)(u(X)-v(X))\neq 0$. Le syst\`eme de
    racines de chaque sous-espace de Cartan   comporte deux racines positives de multiplicit\' e $1$ dont les carr\' es sont les fonctions $4u$ et $4v$, et deux racines de multiplicit\' e $2$ dont le produit est, \`a une fonction localement constante sur $\q^{reg}$  pr\`es,  la fonction $u-v$.  Nous mettons
 en \'evidence
    des sym\'etries entre ces sous-espaces de Cartan et  d\' efinissons un isomorphisme $H$-invariant  $\varpi$ de $\q$  qui renverse
    l'ordre d'Hira\"\i .\me

Nous souhaitons \'etudier  le comportement des int\' egrales orbitales et des distributions propres invariantes  au voisinage des points semi-simples non nuls et non r\' eguliers. Ces \'el\'ements, appel\'es semi-r\'eguliers, sont les \' el\' ements semi-simples $X\in\q$ pour lesquels une et une seule des trois valeurs $u(X)$, $v(X)$ ou $u(X)-v(X)$ est nulle ou de mani\`ere \' equivalente qui
    sont \`a l'intersection d'exactement deux sous espaces de Cartan. 
    
    Le paragraphe
   1.4 est consacr\'e \`a l'\'etude des centralisateurs des  \'el\'ements semi-r\'eguliers.
Tout \' el\' ement semi-r\' egulier  annul\'e
    par une racine de multiplicit\'e $1$ (ce qui correspond \`a la situation $u(X)=0$ ou bien $v(X)=0$ avec $u(X)\neq v(X)$) est conjugu\' e sous $H$ \`a un \' el\' ement de  $\a_{++}\cap\a_{+-}$ ou de  $\a_{+-}\cap \a_{--}$.  Soient $\z_1$ et  $\z_2$ les centralisateurs dans $\g$ de chacun de ces deux espaces. Alors $\m=\z_1+\z_2$ est une sous-alg\`ebre de Lie $\sigma$-stable de $\g$ et la  sous-paire sym\'etrique
    $(\m,\m\cap\h)$ est isomorphe au produit de deux copies $H$-conjugu\' ees de $(\mathfrak{gl}(2,\R),\mathfrak{gl}(1,\R)\times \mathfrak{gl}(1,\R))$.  Les \' el\' ements semi-r\' eguliers annul\' es par une racine de multiplicit\' e $2$ ( ce qui correspond \`a la situation $u(X)-v(X)=0$ et $u(X)v(X)\neq 0$)  sont tous $H$-conjugu\' es \`a un \' el\' ement de $\a_{++}\cap \a_2$ ou de $\varpi(\a_{++}\cap \a_2)=\a_{--}\cap h_0\cdot\a_2$ o\`u $h_0\in H$. Soit $\z_3$ et $\z_4$ les centralisateurs repectifs de ces deux espaces. On a alors les isomorphismes $(\z_3,\z_3\cap\h)\simeq(\mathfrak{gl}(2,\R)\times\mathfrak{gl}(2,\R), diagonale)$ et
$(\z_4,\z_4\cap\h)\simeq (\mathfrak{gl}(2,\C),\mathfrak{gl}(2,\R))$. De plus, on a $\q=\z_3+\z_4$ .\me

L'int\' egrale orbitale  $\mathcal{M} (f)$  d'une fonction
$f\in\mathcal{C}^\infty_c(\mathfrak{q})$ est la fonction
$H$-invariante  de classe
$\mathcal{C}^\infty$ sur $\q^{reg}$ donn\' ee pour $X\in\q^{reg}$ par $\displaystyle
\mathcal{M}(f)(X)=\vert  \delta(X)\vert \int_{H/Z_H(X)}f(h\cdot
X)d\dot{h},$ o\`u $\delta(X)=u(X)-v(X)$ et $Z_H(X)$ est le
centralisateur de $X$ dans $H$.

Pour \'etudier
le comportement  de $\mathcal{M} (f)$  au voisinage des points semi-simples non nuls,  nous g\' en\' eralisons tout d'abord   les r\' esultats
de T. Tengstrand et J.~Faraut concernant les int\' egrales orbitales
sur les hyperbolo\" ides (espace sym\' etrique de rang $1$) . Ceci fait l'objet de la partie $2$.  Etant
donn\' e une forme quadratique $\mathcal{Q}$ sur $\R^n$, la
  fonction moyenne $M_{\mathcal{Q}}(f)$ de  $f\in\mathcal{C}_c^\infty(\R^n)$ introduite par J.~Faraut 
( \cite{fa3} appendice)  est l'application
    $M_\mathcal{Q}(f)\in\mathcal{C}^\infty(\R^*)$ v\' erifiant   pour tout $F\in\mathcal{C}_c^\infty(\R)$,
    la relation
$\displaystyle \int_{\R^n}F\circ \mathcal{Q}(x)f(x)dx=\int_\R
    F(t)M_\mathcal{Q}f(t)dt.$
    Pour $t\in\R^*$, on a
    $M_{\mathcal{Q}}(f)(t)=\phi_0(t)+\eta(t)\phi_1(t), $
    o\`u $\phi_0$ et $\phi_1$ sont des fonctions de classe
    $\mathcal{C}^\infty$ \`a support compact dans $\R$ et $\eta$ est une fonction
    singularit\'e d\'ependant  de la signature de la forme quadratique
    $\mathcal{Q}$.
 
Nous g\' en\' eralisons ce r\' esultat  pour
$f\in\mathcal{C}_c^\infty(\R^m\times\R^n)$ en montrant qu'il existe $\Phi_0$ et
$\Phi_1$ dans $\mathcal{C}_c^\infty(\R^m\times \R)$ v\' erifiant, pour tout
$(x,t)\in\R^m\times\R^*,$ la relation
$M_\mathcal{Q}(f(x,\cdot))(t)=\Phi_0(x,t)+\eta (t) \Phi_1(x,t),$
et toute fonction du type  $\Phi_0(x,t)+\eta (t) \Phi_1(x,t)$ avec
$\Phi_0$ et $\Phi_1$ dans $\mathcal{C}_c^\infty(\R^m\times \R)$ est l'image
$M_\mathcal{Q}(f(x,\cdot))(t)$ d'une fonction
$f\in\mathcal{C}_c^\infty(\R^m\times\R^n)$ (Th\'eor\`emes \ref{comp} et \ref{surj2}).
Nous consid\' erons \' egalement la moyenne
$M_{\mathcal{Q}'}\left(x\mapsto M_{\mathcal{Q}}(
    f(x,.))\right)(t',t)$ relative \`a une forme quadratique $Q'$ de $\R^m$ (Th\'eor\`emes \ref{comp2} et \ref{surj3}).
\me

En adaptant la m\' ethode de descente d'Harish-Chandra \`a la paire
sym\' etrique $(\m,\m\cap \h)$, nous en d\' eduisons le comportement
des int\' egrales orbitales simultan\' ement  au voisinage des
points non nuls de  $\a_{++}\cap \a_{+-}$ et $ \a_{+-}\cap
\a_{--}$  (th\' eor\`eme \ref{iom}): 

Soit $\mathcal{H}^2_{\log}$ l'espace des fonctions sur $(\R^*)^2-diag$ de la forme 
$\varphi_0(t_1,t_2)+\log|t_1|\varphi_1(t_1,t_2)+\log|t_2|\varphi_1(t_2,t_1)+\log|t_1|\log|t_2|\varphi_2(t_1,t_2)$ avec 
$\varphi_0,\varphi_1,\varphi_2 \in\D(\R^2-diag)$ et $
\varphi_0, \varphi_2  \textrm{ sym\' etriques}$.  

Soit $\mathcal {U}_\m=H\cdot \p\m\cap\q=\{ X\in\q; Q^2(X)-4S(X)>0\}$.
Alors, pour tout $f\in {\mathcal C}_c^\infty(H.\p\m\cap\q)$, il existe une
unique fonction $F_f\in \mathcal{H}^2_{\log}$ telle que pour $X\in
H\cdot \p\m\cap\q^{reg}$, l'on ait $\mathcal{M}(f)(X)=F_f\big(
u(X),v(X)\big)$ et l'application $f\mapsto F_f$ est surjective de
$\mathcal{C}_c^\infty(H.\p\m\cap\q)$  dans $ \mathcal{H}^2_{\log}$.

Les autres  points semi-r\' eguliers correspondent \`a la situation
$u(X)=v(X)$ et sont conjugu\' es \`a un \' el\' ement de
$\a_{++}\cap\a_2$ ou de $h_0.\a_2\cap \a_{--}$ avec  $h_0\in H$.
Contrairement \`a la
situation pr\' ec\' edente, ces deux cas ne peuvent \^etre trait\'
es simultan\' ement mais se d\' eduisent l'un de l'autre par
l'isomorphisme $\varpi$ de $\q$. La
m\' ethode de descente classique de Harish-Chandra donne le r\' esultat suivant.

Soit $Y$ la fonction d'Heaviside et soit $\mathcal{H}_Y^{pair}$ l'espace des fonctions sur $\R^*\times \R$ de la forme 
   $ a(\tau,w)+Y(-w)|w|^{1/2}b(\tau,w)$ avec $a,b\in  {\mathcal C}_c^\infty(\R^*\times\R\cap\{(\tau,w); \tau^2>w\}$   paires par rapport \`a $\tau$. On note  ${\mathcal
U}_3=H\cdot\p(\z_3\cap\q)=\{X\in\q; S(X)>0 \textrm{ et }
Q(X)>-2\sqrt{S(X)}\}$. Alors, pour tout $f\in  {\mathcal C}_c^\infty({\mathcal U}_3)$,  il existe une unique fonction
    $G_f\in\mathcal{H}_Y^{pair}$ telle que   $\displaystyle \mathcal{M}(f)(X)=G_f(\tau,w)$, ceci
     pour tout $X\in {\mathcal U}_3\cap\qr$ tels que, pour  $\dfrac{Q(X)}{2}+\sqrt{S(X)}=\tau^2$ et $\dfrac{Q(X)}{2}-\sqrt{S(X)}=w$ avec  $\tau^2>w$. De plus l'application $f\mapsto G_f$ est surjective de $\mathcal{C}_c^\infty({\mathcal U}_3)$ dans $\mathcal{H}_\eta^{pair}$.\\

Les deux derni\`eres parties sont consacr\' ees \`a l'\' etude des distributions propres invariantes sur $\q$. On fixe un caract\`ere $\chi$ de $\cqh=\C[Q,S]\simeq S(\q)^H$ donn\' e par deux complexes $\la_1$ et $\la_2$ v\' erifiant $\chi(Q)=\la_1+\la_2$ et $\chi(S)=\la_1\la_2$ avec $\la_1\la_2(\la_1-\la_2)\neq 0$. \me

Dans la partie 4, nous d\' ecrivons les  solutions analytiques
du syst\`eme   suivant~: $\big(E_\chi(\qr)\big)$    pour tout $P\in \cqh$ alors $ \partial(P)F=\chi(P)F$ sur $q^{reg}.$

On fixe un sous-espace de Cartan $\a$ de $\q$ dont on note
$W_{H_\C}(\a_\C)$ le quotient du normalisateur  de $\a_\C$ dans
$H_\C$ par son centralisateur. Les   op\' erateurs de Dunkl
permettent de d\' ecrire l'action des composantes radiales des op\'
erateurs $\partial(Q)$ et $\partial(S)$ sur l'espace des fonctions
$W_{H_\C}(\a_\C)$- invariantes . Plus pr\' ecis\' ement, notons   $\al_1$ et $\al_2$ les racines de multiplicit\' e $1$ de $\a$ dans $\g$, alors, la fonction
$H$-invariante $F$ sur $\q^{reg}$ satisfait le syst\`eme
$\big(E_\chi(\qr)\big)$ si et seulement si, pour tout sous-espace de Cartan
$\a$ de $\q$, on a\me

$\left\{\begin{array}{c}\big(\partial(Q) F\big)_{/\a^{reg}}=\de^{-1}\circ(\elu+\eld)\circ\de(F_{/\a^{reg}})=(\la_1+\la_2)F_{/\a^{reg}}\\
 \big(\partial(S) F\big)_{/\a^{reg}}=\de^{-1}\circ(\elu\eld)\circ\de(F_{/\a^{reg}})=(\la_1\la_2)F_{/\a^{reg}}\end{array}\right.$\me

\no o\`u $L_{\al_i}=\frac{1}{4\al_i}\partial
(\al_i)(\al_i\partial (\al_i)).$ Cet op\' erateur est la composante
radiale de l'op\' erateur de Laplace de la paire sym\' etrique de
rang un
    $(\mathfrak{so}(2,1),\mathfrak{so}(1,1))$ (appendice A de
    \cite{fa2}) lorsque la racine $\al_i$ est r\'eelle ou imaginaire. 

Les
op\'erateurs $L_{\al_1}$ et $L_{\al_2}$ commutent car  $\al_1$ et $\al_2$ sont orthogonales, ainsi  le syst\`eme
pr\' ec\' edent se ram\`ene \`a l'\' etude de $\ker(L_{\al_j}-\la id)$
pour $j=1,2$. Via le changement de variable $z=\al_j^2$, nous
montrons que ses solutions   s'expriment \`a l'aide des solutions
d'une part sur $\C$ et d'autre part sur $\R$  de l'\' equation de
type Bessel  $(E_\C)\quad L_c\Phi= 4z\frac{\partial^2\Phi }{\partial
z^2}+4\frac{\partial\Phi}{\partial z}=\lambda \Phi$ .

Un syst\`eme fondamental de l'espace ${\mathcal S}ol(L_c,\la)$ des
solutions de $(E_\C)$ sur $\C-\R_-$ est donn\' e par une fonction
$\Phi_\la$ analytique sur $\C$ et une fonction
$W_\la(z)=w_\la(z)+\log(z) \Phi_\la(z)$ o\`u $\log(z)$ d\' esigne la
d\' etermination principale du logarithme et  $w_\la$ est une
fonction analytique sur $\C$. L'espace des solutions ${\mathcal
S}ol(L,\la)$ sur $\R^*$ de l'\' equation $4t\frac{\partial^2\Phi
}{\partial t^2}+4\frac{\partial\Phi}{\partial t}=\lambda \Phi$ est
engendr\' e par les fonctions $\Phi_\la(t)$ et
$W_\la^r(t)=w_\la(t)+\log\vert t\vert \Phi_\la(t)$ avec $t\in\R^*$.
On constate en particulier que la solution $W_\la(z)$ ne se prolonge
pas sur $\R^*$ en $W_\la^r(t)$ ce qui joue un r\^ole important par
la suite.

La fonction $F$ solution de $\big(E_\chi(\qr)\big)$ s'\' ecrit alors sur chaque composante connexe de $\q^{reg} $
comme combinaison lin\' eaire de fonctions $\displaystyle \frac{A(u(X))B(v(X))}{u(X)-v(X)}$ et
$\displaystyle \frac{A(v(X))B(u(X))}{u(X)-v(X)}$ o\`u $(A,B)$ parcourt ${\mathcal S}ol(L,\la_1)\times{\mathcal S}ol(L,\la_2)$ lorsque $u(X)$ et $v(X)$
 sont r\' eelles et ${\mathcal S}ol(L_c,\la_1)\times{\mathcal S}ol(L_c,\la_2)$ sinon.\me

Dans la partie $5$, nous consid\`erons  une  solution  $F$ du syst\`eme$\big(E_\chi(\qr)\big)$.  Nous  cherchons   des conditions n\' ecessaires
et suffisantes pour que la distribution $T_{\partial(P)F}-\partial(P)T_F$ soit nulle pour tout $P\in S(\q)^H$.
Comme dans le cas du groupe (\cite{hi}), la formule d'int\' egration de Weyl permet une \' etude sur chaque sous-espace de Cartan.
 Une int\' egration par parties tenant compte du  comportement des int\' egrales orbitales obtenu pr\' ec\' edemment nous  permet de
 d\' egager des conditions n\'ecessaires  qui portent sur le comportement de $F$ au voisinage des points semi-r\' eguliers. 

Nous \' etudions tout d'abord les conditions n\' ecessaires  successivement sur les ouverts
 ${\mathcal U}_\m=H\cdot \p\m\cap\q,$
${\mathcal U}_3=H\cdot\p(\z_3\cap\q)$
et $ \varpi({\mathcal U}_3)=H\cdot\varpi(\p(\z_3\cap\q))$ qui forment un recouvrement ouvert de $\q-{\mathcal N}={\mathcal U}$.

Sur l'ouvert ${\mathcal U}_\m$, la partie $4$ assure qu'il existe une  fonction
$\Psi_\m\in {\mathcal C}^\infty((\R^*)^2-diagonale)$,  sym\' etrique telle que, pour $X\in{\mathcal U}_\m\cap
\q^{reg}$ , on ait $F(X)=\Psi_\m(t_1,t_2)/|t_1-t_2|$ avec $\{u(X),
v(X)\}=\{ t_1, t_2\}\subset \R^*)^2-diagonale$. 

 La fonction $\Psi_\m$ et les
int\' egrales orbitales d'une fonction de $\mathcal{C}_c^\infty({\mathcal U}_\m)$ ont
des singularit\' es de type $\log|t_j|$.  Nous
obtenons des conditions de recollement analogues \`a celles de J.Faraut  pour un
hyperbolo\"
    ide \` a une nappe (\cite{
     fa3}) qui  s'expriment de la mani\`ere suivante:

     Si $u$ est une fonction de la forme $u(t)=v(t)+\log|t|w(t),$ o\`u  $v, w\in \mathcal{C}^2(\R^*)$   admettent,
ainsi que leurs d\' eriv\' ees, des limites  \` a droite et \` a
gauche en $0$. On pose
$u^{[1]}(t)=tu'(t)\qq\textrm{et}\qq u^{[0]}(t)=u(t)-\log|t|u^{[1]}(t),$
de telle sorte que
$ \lim_{t\to0^{\pm}}u^{[1]}(t)=w(0^\pm)\qq\textrm{et}\space
\lim_{t\to0^{\pm}}u^{[0]}(t)=v(0^\pm).$
Alors pour tout $ j\in\{0,1\}$ et
pour tout $ t\in\R^*$, on a  (proposition \ref{recolm})
$$(\Psi_\m(t,.))^{[j]}(0^{+})=(\Psi_\m(t,.))^{[j]}(0^{-}) \;\textrm{ et }\; (\Psi_\m(., t))^{[j]}(0^{+})=(\Psi_\m(., t))^{[j]}(0^{-}).$$
   On note pour $f$ et $g$ deux fonctions d'une variable  complexe, $S^+(f,g)(z_1,z_2)=f(z_1)g(z_2)+f(z_2)g(z_1)$ 
et  $[f,g](z_1,z_2)=f(z_1)g(z_2)-f(z_2)g(z_1).$ On obtient alors  (corollaire \ref{expsolm}) que la fonction $\Psi_\m$ est combinaison lin\' eaire de fonctions 
$S^+(A,B)(t_1,t_2)$ et $signe(t_1-t_2)[A,B](t_1, t_2)$ lorsque  $(A,B)$ parcourt 
$\mathcal{S}ol(L,\la_1)\times\mathcal{S}ol(L,\la_2)\rangle.$\me

Pour l'\' etude sur l'ouvert $\mathcal{U}_3$,  on d\' efinit  la fonction  $(\Psi_\m)_r$  sur $\{(\tau,\te)\in (\R^*_+)^2;
\tau^2\neq \te^2\}$ par
 $(\Psi_\m)_r(\tau,\te)=\Psi_\m\big((\tau+\te)^2,(\tau-\te)^2\big)$    et  la fonction $\Psi_2$ paire en chaque variable sur $(\R^*)^2$ par   $\Psi_2(\tau,\te)=(u-v)(X_{\tau,\te})F(X_{\tau,\te})$ o\`u l'orbite $H\cdot X_{\tau,\te}$ est caract\' eris\' e  par $\{(\tau+\im\te)^2,(\tau-\im\te)^2\}$.
Les conditions n\'
ecessaires obtenues sont alors (proposition \ref{recolztp}) \me

$\forall \tau\in\R_+^*,\;\left\{\begin{array}{c}\frac{\partial}{\partial\te}(\Psi_\m)_r(\tau,0^+)-
    \frac{\partial}{\partial\te}\Psi_2(\tau,0^+)=0\\
    \Psi_2(\tau,0^+)=0\;\;\qq\qq\end{array}\right..$

\no ce qui donne les expressions suivantes des fonctions $\psi_\m$ et $\psi_2$ (corollaire
\ref{solztp}):\me

$\begin{array}{l}\mbox{ pour } t_1>t_2>0 \textrm{ alors }
    \Psi_\m(t_1,t_2)=\Psi^+(t_1,t_2)+\im\Psi_c(t_1,t_2) \\
    \mbox{ pour }\tau>0\mbox{ et }\te>0 \textrm{ alors }
     \Psi_2(\tau,\te)=\Psi_c(u(X_{\tau,\te}),v(X_{\tau,\te})), \end{array}$\me
     
\no o\`u  $\Psi^+\in Vect\langle S^+(A,B); (A,B)\in
    \mathcal{S}ol(L,\la_1)\times \mathcal{S}ol(L,\la_2)\rangle$
  et $ \Psi_c\in Vect \langle [A,B]; (A,B)\in
    \mathcal{S}ol(L_c,\la_1)\times \mathcal{S}ol(L_c,\la_2)\rangle $

On utilise ensuite l'application $\vi$,  pour obtenir les conditions de recollement et
les fonctions correspondantes sur $\vi(\mathcal{U}_3)$. \me

La synth\`ese des r\' esultats pr\' ec\' edents permet de d\' ecrire
une base de l'espace des distributions propres invariantes sur
${\mathcal U}=\q-{\mathcal N}$ d\' efinies par une fonction
localement int\' egrable sur  ${\mathcal U}$. Celle-ci  est  donn\'
ee par les fonctions suivantes (th\'eor\`eme \ref{solgen}):\me
    
    \no$(1)$
    La fonction analytique  $\dis F_{ana}=\frac{[\Phi_{\la_1},\Phi_{\la_2}](u,v)}{u-v}$
    qui s'exprime \' egalement sous  la forme $f(Q,S)$ o\`u $f$ est analytique sur $\R^2$. Cette distribution est propre invariante sur $\q$ tout entier.

   \no$(2)$ La fonction $\dis F_{sing}=\frac{[\Phi_{\la_1},w_{\la_2}](u,v)+[w_{\la_1},\Phi_{\la_2}](u,v)+\log|uv|[\Phi_{\la_1},\Phi_{\la_2}](u,v)}{u-v}$
    qui s'exprime \'egalement sous  la forme $f(Q,S)+\log|S|g(Q,S)$ o\`u $f$ et $g$ sont  analytiques sur $\R^2$. Cette fonction est localement  int\'egrable sur $\q$.
  
  \no$(3)$  Les fonctions du
    type $F^+_{A,B}=Y(Q^2-4S)\dfrac{A(u)B(v)+A(v)B(u)}{u-v}$, o\`u $Y$ d\' esigne la fonction d'Heaviside et $(A,B)$ parcourt $ {\mathcal S}ol(L,\la_1)\times{\mathcal S}ol(L,\la_2).$\me

   Ces r\' esultats sont partiels et  montrent la limite de la  m\' ethode de descente. D'une part, nous n'obtenons pas  le comportement  des int\' egrales orbitales au voisinage de $0$ et ceci est n\' ecessaire pour  savoir si les fonctions $F_{sing}$ et $F^+_{A,B}$ se prolongent en des distributions propres invariantes sur $\q$ tout entier. Par ailleurs, les restrictions \`a $\z_3\cap\q$ des fonctions obtenues n'ont aucun lien avec les distributions propres invariantes pour la paire sym\' etrique $(\z_3,\z_3\cap\h)\simeq\mathfrak{gl}(2,\R)$.

Ainsi,  notre \' etude montre que si $T$ est une distribution propre invariante sur $\q$ alors sa restriction $T_{/\qr}=F$ est une fonction localement int\' egrable sur ${\mathcal U}$. Mais les questions suivantes restent ouvertes: la fonction $F$ est-elle localement int\' grable sur $\q$?  La distribution $T_F$ est-elle propre invariante sur ${\mathcal U}$? sur $\q$?

\section{Pr\' eliminaires}

\subsection{Notations  }
 Si  $M$  est une vari\'et\'e
diff\'erentiable, on note $C^{\infty}(M)$ l'espace des
fonctions de classe $C^{\infty}$ sur $M$ ,
${\cal D}(M)$ le sous-espace de $C^{\infty}(M)$ des fonctions \`a support
compact et 
${\cal D}(M)'$ l'espace des distributions sur $M$, 
c'est-\`a-dire le dual de ${\cal D}(M)$.

Pour $N\subset M$ et f  une fonction d\' efinie sur M, on notera $f _{/N}$
 sa restriction \`a
$N$. 
\me

 Si $\Om$ est un ensemble fini, alors  $\mid \Om \mid$ d\' esigne son cardinal.

  Le spectre d'une matrice $X$ sera not\' e $Sp(X)$.
 
  On notera $\log$ la d\' etermination principale du logarithme d\'
efinie sur $\C-\R_-$ par $\log(z)=\log\vert z\vert +\im Arg(z)$ o\`
u l'argument $Arg(z)$ de $z$ appartient \` a $]-\pi,\pi[$.

 La fonction $Y(t)$ d\'esignera la fonction d'Heaviside sur $\R$. \me

Soit $\mathbf{G}$ un groupe alg\' ebrique lin\' eaire d\' efini sur $\R$. On
note $G$ l'ensemble de ses points r\' eels et $G_\C$ l'ensemble de
ses points complexes. Soit  $\g$ l'alg\`ebre de Lie de $G$. 
Pour $g\in G$ et $X,Y\in \g$, on notera $g\cdot X=\Ad(g)X$ (resp. $\ad(Y)X=[Y,X]$) l'action adjointe de $G$ (resp. $\g$) sur $\g$.\me

Soit $M$ un sous-groupe de $G$ et   $M_{\C}$ son
complexifi\' e dans
 $G_{\C}$.
Si  $U$ une partie de $\g$ alors, on note
$Z_{M}(U)=\{m\in M; m\cdot u=u \textrm{ pour tout } u\in U\} $ le centralisateur de $U$ dans $M$ et $N_{M}(U)=\{m\in M; m\cdot U\subset  U\} $ son
normalisateur. De m\^{e}me pour tout sous-espace vectoriel
$\mathfrak{v}$ de $\g ,$ on note $\mathfrak{z}_{\mathfrak{v}}(U)$ le
centralisateur de $U$ dans $\mathfrak{v}$ et $\n_{\mathfrak{v}}(U)$
son normalisateur.\me

%%%%%%%%%%%%%%%%%%%%%%%%%%%%%%%%
Soit $V$ un espace vectoriel r\' eel de dimension finie. On notera $V^{*}$ son dual 
et $V_{\Bbb C}$ son complexifi\' e. L'alg\`{e}bre sym\' etrique $S[V]$ de $V$
 s'identifie d'une part \`a
l'alg\` ebre $\R[V^*]$ des fonctions polynomiales \` a coefficients
r\' eels sur $V^*$  et d'autre part,  \` a
l'alg\`{e}bre des op\' erateurs diff\' erentiels \`{a} coefficients
constants r\' eels  sur $V$.  De m\^eme, on a  $S[V_\C]=\C[V^*]$ et
cette alg\` ebre s'identifie  \` a l'alg\`{e}bre des op\' erateurs
diff\' erentiels \`{a} coefficients constants complexes  sur  $V_\C$. Si $u\in S[V]$ (resp. 
$S[V_\C]$), on note $\partial(u)$ l'op\' erateur diff\' erentiel
correspondant.

On
note $S[V]^{M}$ (respectivement $S[V_\C]^M$) la sous alg\`{e}bre de
$S[V]$ (respectivement $S[V_\C]$) constitu\' ee des \' el\' ements
invariants sous l'action d'un sous-groupe  $M$ de $G$. \me

On consid\`ere $G=GL(4,\R)$
et son alg\`ebre de Lie $\g=\mathfrak{gl}(4,\R)$ munis de 
l'involution $$ \sigma(X)=\left (
\begin{array}{cc}
I_{2}&0\\0&-I_{2}
\end{array}\right )X
\left (
\begin{array}{cc}I_{2}&0\\0&-I_{2}
\end{array}\right ),
$$ o\`u $I_{n}$ est la matrice identit\' e de taille
$n\times n$. 

Le groupe $\displaystyle H=\left\{\mah; A,B\in
GL(2,\R)\right\}$ est alors le sous-groupe des points de $G$ fixes
sous l'action de $\sigma$ dont l'alg\`ebre de Lie est $\h=\{X\in\g; \sigma(X)=X\}$. On note $\q=\{X\in\g; \sigma(X)=-X\}$ de telle sorte que $\g=\h\oplus\q$.  On a
$$\h=\left\{\mah; A,B\in\mathfrak{gl}(2,\R)\right\},\quad \textrm{ et }\quad\q=\left\{\maq; Y,Z\in\mathfrak{gl}(2,\R)\right\}.$$
L'espace $\q$ est stable sous l'action de $H$  et  plus pr\' ecis\' ement,
pour $\mah\in H$ et $\maq\in\q$, on a $\dis \mah\cdot\maq=\left
(\begin{array}{c|c} 0&AYB^{-1}\\ \hline BZA^{-1}&0\end{array}\right
).$\bigskip

On note $\q^{ss}$ l'ensemble des \' el\' ements   semi-simples  c'est-\`a-dire des matrices semi-simples, de
$\q$. L'alg\`ebre de Lie $\g$ \' etant de rang $4$, il existe des
polyn\^omes $d_4\neq 0, \ldots,d_{16}$ sur $\g$ tels que pour tout
$X\in\g$, l'on ait $\displaystyle
\det(t-\mathrm{ad}(X))=\sum_{j=4}^{16}d_j(X)t^j. $ Un \' el\' ement  $X\in \g$ semi-simple est dit r\' egulier dans $\g$ si
$d_4(X)\ne0$ ce qui \' equivaut \`a dire que   les  valeurs propres de $X$ sont deux \`a deux distinctes. 

\medskip

\begin{Def} Soit $l$ le plus petit entier tel que $d_{l /\q}\ne0$.
Un \' el\' ement $X$ de $\q$ est $\q$-r\' egulier s'il est semi-simple et satisfait   $d_l(X)\ne0$.
On note $\q^{reg}$ l'ensemble des \' el\' ements $\q$-r\' eguliers
de $\q$.\end{Def}

Les valeurs propres de la matrice 
$Y=\left(\begin{array}{llll}0&0&\la_1&0\\0&0&0&\la_2\\\la_1&0&0&0\\0&\la_2&0&0\end{array}\right)\in\q $
sont $\pm \la_1$ et $\pm \la_2$. Ainsi, on a $l=4$ et   le r\' esultat suivant est imm\' ediat.

\begin{lem} Soit $X\in \q$. Alors $X$  est $\q$-r\' egulier si et seulement si $X$  est semi-simple et 
r\' egulier dans $\g$.

\end{lem}

On appelle sous-espace de Cartan $\a\subset\q$ un sous-espace ab\'
elien, form\' e d'\' el\' ements semi-simples et maximal pour cette
propri\' et\' e. On note $car(\q)$ l'ensemble des sous-espaces de
Cartan de $\q$.

Soit $\a$ un sous-espace de Cartan. On note $\ra$ le syst\`eme de
racines de $\a_\C$ dans $\g_\C$ et lorsque l'on choisit un ordre sur
celui-ci, on note $\rap$ le syst\`eme de racines positives
correspondant. Soit $\wg$ le groupe de Weyl de $\ra$ et 
$W_H(\a)=\{ Ad(h)_{/\a}; h\in N_H(\a)\}\simeq N_H(\a)/Z_H(\a)$ le groupe de Weyl de $\a$ dans $H$.  

Pour $\al\in\ra$, on note 
$\dis \g_\C^\al=\{Y\in\g_\C,\;\forall X\in\a, [X,Y]=\al(X)Y\},$ l'espace
radiciel correspondant et  $m_\al=dim_\C \g_\C^\al $ la multiplicit\' e de
$\al$. \\

On rappelle les r\'
esultats classiques suivants:

\begin{prop} (paragraphe
1 de \cite{or}) (a) Les \' enonc\' es suivants sont \' equivalents:
\begin{enumerate}
\item $X$ est $\q$-r\' egulier.
\item $\z_\q(X)$ est un sous espace de Cartan de $\q$.
\item Si $\a$ est un sous-espace de Cartan de $\q$ contenant $X$,
alors $X$ n'annule aucune racine de $\a_\C$.
\end{enumerate}
\noindent (b) Tout \' el\' ement r\' egulier appartient \`{a} un
unique sous-espace de Cartan de $\q$.
\end{prop}

On note $\om$ la forme bilin\' eaire sym\' etrique sur $\g$ d\' efinie par 
par
$\om(X,X')=\frac{1}{2}tr(XX')$ pour $X $ et $X' $  dans $\g.$
La restriction de $\om$ \`a $[\g,\g]$ co\" incide, \`a un facteur
multiplicatif pr\`es, avec la forme de Killing. La forme $\om$ est  $H-$ invariante et  non d\' eg\' en\' er\' ee. Il est en de m\^eme de ses restrictions  \`a $\h\times\h$,  \`a $\q\times\q$  et  \`a $\a\times \a$ pour tout $\a\in car(\q)$.
\medskip

Soit $\a\in car(\q)$. 
Ainsi, pour tout $\al\in\ra$, il existe un
unique \' el\' ement  $h_\al\in \a_\C$  tel que
$\om(h_\al,X)=\al(X)$  pour tout $X\in\a$. On d\' efinit la partie compacte
$\a_I=\left(\sum_{\al\in\ra}\im\R h_\al)\right)\cap\a$ et la partie
d\' eploy\' ee $\a_R=\left(\sum_{\al\in\ra}\R h_\al)\right)\cap\a$
de $\a$
de telle sorte que l'on ait $\a=\a_R+\a_I.$\\
Une racine $\al$ est dite r\' eelle (respectivement imaginaire) si
$\al(\a)\subset\R$ (respectivement $\al(\a)\subset\im\R$), ceci est
\' equivalent \`a $h_\al\in\a_R$ (respectivement $h_\al\in\im\a_I$).\me

On note $\langle \a\rangle$ un repr\' esentant
de la classe de conjugaison de $\a$ sous  l'action de $H$.
L'ordre d'Hira\"i sur les classes modulo
$H$ des
sous-espaces de Cartan est d\' efini de la mani\`ere suivante: 
soient $\a$ et $\b$ deux sous-espaces de Cartan. On dit que $\langle
\a\rangle\leq\langle \b\rangle$ si et seulement si il existe $h\in
H$ tel que $(h\cdot\a)_R\subset\b_R$. Lorsque cette inclusion est
stricte, on dit alors que $\langle \a\rangle<\langle \b\rangle$.

\subsection{Sous-espaces de Cartan}\label{cartans}

On d\' efinit les sous-espaces de Cartan $\a_{\ep_1,\ep_2}$ pour $\ep_1$ et $\ep_2=\pm$, et $\a_2$ de $\q$ de  la mani\`ere suivante:
$$\a_{\ep_1,\ep_2}=\left\{X_{u_1, u_2}^{\ep_1,\ep_2}=\left(
\begin{array}{c|c}
0&\begin{array}{cc} u_{1}&0\\0&u_{2}\end{array}\\ \hline
\begin{array}{cc} \ep_1 u_{1}&0\\0&\ep_2 u_{2}\end{array}&0
\end{array}
\right); (u_1,u_2)\in\R^2\right\}$$
$$\a_{2}=\left\{X_{\tau, \te} =\com;(\te,\tau)\in\R^2\right\}.$$

\no Les sous-espaces $\a_{+,-}$ et $\a_{-,+}$ sont conjugu\'es par $K= \left(
\begin{array}{c|c}
\begin{array}{cc} 0&1\\1&0\end{array}&0\\ \hline
0&\begin{array}{cc} 0&1\\1&0\end{array}
\end{array}
\right)$.

\begin{lem}\label{vp} \begin{enumerate}\item  Soit $X=\left(\begin{array}{cc} 0& Y\\ Z&0\end{array}\right)\in\q^{reg}$ .
Si $\la$ est une valeur propre  de $X$ alors $\la\neq 0$ et $-\la$
est aussi une valeur propre de $X$. On a $ Sp(YZ)=\{ \la^2;\la\in
Sp(X)\}$.
\item La famille $<car(\q)>= \{\a_{+,+}, \a_{+,-}, \a_{-,-},\a_2\}$ est une
famille repr\' esentative des classes de conjugaison sous $H$ des
sous-espaces de Cartan de $\q$.
\end{enumerate}
\end{lem}

\begin{dem}  $\it{1}.$  Soient $U$ et $V$ dans $\C^2$ tels que $\left(\begin{array}{c} U\\ V
\end{array}\right)$ soit un vecteur propre non nul de $X$ pour la valeur
propre $\la$. On a alors $YV=\la U$ et $ZU=\la V$. Si $\la\neq 0$
alors $U$ et $V$ sont non nuls. Par suite, le vecteur
$\left(\begin{array}{c} U\\ -V \end{array}\right)$ est un vecteur
propre non nul de $X$ pour la valeur propre $-\la$.\me

Comme $X$ est r\' egulier dans $\q$, ses valeurs propres sont deux
\`a deux distinctes. D'autre part, si $0$ \' etait une valeur propre
de $X$ alors, par ce qui pr\' ec\`ede, elle serait double ce qui est
impossible puisque $X$ est r\' egulier. Ainsi, on a $Sp(X)=\{
\pm\la_1,\pm\la_2\}$ o\`u $\la_1$ et $\la_2$ sont deux nombres
complexes non nuls tels que $\la_1\neq\pm \la_2$. En particulier,
les matrices $X$, $Y$ et $Z$ sont inversibles.

Pour $j=1\mbox{ ou }2 $, on fixe un vecteur propre non nul
$\left(\begin{array}{c} U_j\\ V_j \end{array}\right)$ associ\' e \`a
la valeur propre $\la_j$. Les vecteurs $U_1$ et $U_2$ de $\C^2$ sont
non nuls puisque la famille $\{\left(\begin{array}{c} U_j\\ V_j
\end{array}\right),\; \left(\begin{array}{c} U_j\\ -V_j
\end{array}\right) \}_j$ est une base de $\C^4$ et ils v\' erifient
$YZU_1=\la_1^2 U_1$ et $ YZU_2=\la_2^2 U_2$. On obtient donc la
premi\`ere assertion.\me

${\it 2}.$  Le spectre $\{\la_1^2,\la_2^2\}$ de $YZ$ est donn\' e par les racines du polyn\^ome $x^2-tr(YZ)\; x+
det(YZ)\in\R[x]$. Ainsi, soit $\la_1^2$ et $\la_2^2$ sont r\' eels,
soit $\la_1^2\in\C\setminus\R$ et $\la_2=\pm\overline{\la_1}$.

Si $\la_1^2$ et $\la_2^2$ sont r\' eels, il existe une matrice $P\in
GL(2,\R)$ telle que $PYZP^{-1}= \left(\begin{array}{cc} \la_1^2& 0\\
0&\la_2 ^2\end{array}\right)$. Notons $sgn(x)$ le signe de $x\in\R$.
On peut donc \' ecrire $PYZP^{-1}=D D'$ o\`u
$D=\left(\begin{array}{cc}\vert \la_1\vert & 0\\ 0&\vert \la_2\vert
\end{array}\right)$ et $
D'=\left(\begin{array}{cc} sgn(\la_1^2)\vert \la_1\vert & 0\\
0&sgn(\la_2^2)\vert \la_2\vert \end{array}\right).$ On obtient alors 
$\left(\begin{array}{cc} P& 0\\
0&D'PZ^{-1}\end{array}\right)\cdot X=\left(\begin{array}{cc} 0&D\\
D'&0
\end{array}\right)$
et dans ce cas $X$ est $H$-conjugu\' e \` a un \' el\' ement de $\a_{+,+}$, $\a_{+,-}$ ou $\a_{-,-}$.\\

Si $\la_2=\pm\overline{\la_1}$, on \' ecrit $\la_1=\tau+i\te$ avec
$(\te,\tau)\in\R^2$. Il existe alors une matrice $P\in GL(2,\R)$
telle que $PYZP^{-1}=M^2$ avec $M=\left(\begin{array}{cc}\tau &
-\te\\ \te&\tau
\end{array}\right).$ On obtient donc $\left(\begin{array}{cc} P& 0\\
0&MPZ^{-1}\end{array}\right)\cdot X=\left(\begin{array}{cc} 0&M\\
M& 0
\end{array}\right)$
et dans ce cas $X$ est conjugu\' e \` a un \' el\' ement de $\a_2$.
\me

Soit $\a\in car(\q)$ et
$X\in \a^{reg}$ tel que $\a=\z_\q(X)$. Par ce qui pr\' ec\`ede, il
existe $h\in H$ tel que $h.X$ appartienne \`a $\a_{+,+}$,
$\a_{+,-}$, $\a_{-,-}$ ou $\a_2$. Comme $\z_{\q}(h.X)=h.\a$, on
obtient la deuxi\`eme assertion.  \end{dem}

Pour   $\a\in car(\q)$, on note $\al_1$ et $\al_2$ les  deux racines  de multiplicit\' e un dont on fixe les valeurs sur  $\a\in <car(\q)>$ de la mani\`ere suivante:

$$\begin{array}{|c|c|c|c|c|} \hline  & X_{u_1,u_2}^{++}\in \a_{++}& X_{u_1,u_2}^{+-}\in  \a_{+-}& X_{u_1,u_2}^{--}\in \a_{--}
& X_{\tau,\te}\in  \a_{2}\\
 & & & & \\
\hline \al_1(X)  & 2u_1     & 2u_1 & 2iu_1 
&2(\tau+\im \te)\\
\hline \al_2(X)  & 2u_2     &2iu_2  &2iu_2  
&2(\tau-\im \te)\\
\hline\end{array}$$

Dans tout cet article, on fixe  le  syst\`eme de racines positives en posant  $\rap=\{\al_1,\al_2,\be_1,\be_2\}$ o\`u    les deux racines $\be_1=\dfrac{\al_1-\al_2}{2}$ et $\be_2=\dfrac{\al_1+\al_2}{2}$ sont de multiplicit\' e deux. 

Ainsi,
suivant l'ordre d'Hira\"i, nous avons
$$\langle\amm\rangle<\langle\apm\rangle<\langle\app\rangle\quad \mbox{ et }\quad
\langle\amm\rangle<\langle\a_2\rangle<\langle\app\rangle.$$ 

\no avec  les relations suivantes $\dis \app\cap\apm=\R h_{\al_1}=Ker\;\al_2$, $\apm\cap\amm=\R \im h_{\al_2}=Ker\;\al_1$, $\app\cap\a_2=\R h_{\be_2} =Ker\;\be_1$  et  $\amm\cap h_0\cdot \a_2=\R \im h_{\be_1}=Ker\;\be_2,$
o\`u $h_0=\left(
\begin{array}{c|c}
\begin{array}{cc} 0&1\\1&0\end{array}&0\\ \hline
0&\begin{array}{cc} 1&0\\0&-1\end{array}
\end{array}
\right)$.

On d\' efinit l'involution  $H$-\' equivariante  $\varpi $ de $\q$ par $\varpi (\maq)=\left (\begin{array}{c|c}
0&Y\\ \hline -Z&0\end{array}\right).$ Cet automorphisme v\' erifie les propri\' et\'es suivantes:

 \begin{enumerate}
\item pour $(X,X')\in\q^2,$ on a $ [\varpi(X),\varpi(X')]=-[X,X'].$
\item pour $(X,X')\in\h\times\q,$ on a $[X,\varpi(X')]=\varpi([X,X']),$
ceci provenant de la $H$-\' equivariance de $\varpi.$
\item $\varpi(\app)=\amm\quad \varpi(\apm)=\amp\mbox{ et } \varpi(\a_2)=h_0\cdot \a_2,$ et plus pr\' ecis\' ement, pour $X_{\tau,\te}\in\a_2$, on a
$\varpi(X_{\tau,\te})=h_0\cdot X_{\te,\tau}.$

\end{enumerate}

En particulier, l'involution $\varpi$ renverse l'ordre d'Hira\" i sur $car(\q)$. \bigskip

Pr\' ecisons maintenant les groupes $W_H(\a)$ pour $\a\in < car(\q)>$. On pose $$ \ka=\Ad K=\Ad \left(
\begin{array}{c|c}
\begin{array}{cc} 0&1\\1&0\end{array}&0\\ \hline
0&\begin{array}{cc} 0&1\\1&0\end{array}
\end{array}
\right)\quad\mbox{ et } \quad
\varrho=\Ad\left(
\begin{array}{c|c}
\begin{array}{cc} 1&0\\0&-1\end{array}&0\\ \hline
0&\begin{array}{cc} 1&0\\0&1\end{array}
\end{array}
\right).$$

Sur $\a_{+,+}$ et $\a_{-,-}$, l'\' el\' ement $\ka$ \' echange les
racines $\al_1$ et $\al_2$ et l'\' el\' ement $\varrho$ \' echange
les racines $\be_1$ et $\be_2$ et on obtient alors:

 \begin{lem}\label{WH}

\begin{enumerate}
\item $W_H(\a_{+,+})=W_H(\a_{-,-})=\{\pm I_4, \pm \ka, \pm\varrho,\pm\ka\varrho\}$
\item $W_H(\a_{+,-})=\{\pm I_4, \pm \varrho\}$
\item $W_H(\a_{2})=\{\pm I_4, \pm \ka\}$
\end{enumerate}
\end{lem}

\no \begin{dem} Soit $h=\left(\begin{array}{c|c} A& 0\\ \hline 0&
B\end{array}\right)\in N_H(\app)$. Pour tout $X^{++}_{u_1,u_2}\in
\app$ il existe  $X^{++}_{v_1,v_2}\in
\app$ tel que $h\cdot X^{++}_{u_1,u_2}=X^{++}_{v_1,v_2}$. Ainsi,  on a
$A\left( \begin{array}{cc}u_1&0\\0&u_2\end{array}\right)B^{-1}=B\left( \begin{array}{cc}u_1&0\\0&u_2\end{array}\right)A^{-1}=\left( \begin{array}{cc}v_1&0\\0&v_2\end{array}\right)$. Par suite, pour $M=A$ ou $B$,  on a
$M\left(
\begin{array}{cc}u^2_1&0\\0&u^2_2\end{array}\right)M^{-1}=\left(
\begin{array}{cc}v^2_1&0\\0&v^2_2\end{array}\right)$ et donc,
$M$ est diagonale ou $KM$ est diagonale. La premi\` ere relation
implique que, soit $A$ et $B$ sont diagonales et $A^2=B^2$, soit
$\left(\begin{array}{cc} 0& 1\\ 1& 0\end{array}\right)A$ et
$\left(\begin{array}{cc} 0& 1\\ 1& 0\end{array}\right)B$ sont
diagonales et ont  m\^eme carr\' e. Maintenant, on remarque que
$Z_H(\app)=\left\{\left(\begin{array}{c|c}
\begin{array}{cc}a&0\\0&b\end{array}&0\\ \hline
0&\begin{array}{cc}a&0\\0&b\end{array}\end{array}\right)\in
GL(4,\R)\right\}$, ce qui donne  ais\' ement la premi\` ere
assertion. Les autres assertions se prouvent de m\^eme.\end{dem}

%%%%%%%%% ORBITES
\subsection{Orbites semi-simples}\label{OS}

 Nous allons d\' ecrire maintenant les orbites
semi-simples de $\q$ sous l'action de $H$.\me

La paire sym\' etrique $(\g,\h)$ est de rang $2$. Par suite, une base de $\cqh$ est donn\'ee par les polyn\^omes $Q$ et $S$ d\'efinis , pour $X=\maq\in\q$, par
$$Q(X)=\dfrac{1}{2}tr(X^2)=tr(YZ)=\om(X,X)  \textrm{ et }  S(X)=det(X)=det(YZ).$$ 

 Pour $X=\maq$,  les racines du polyn\^ome $t^2-Q(X) t+S(X)$ forment le  spectre de $YZ$. On pose 
  $$S_0=Q^2-4S\textrm{ et }\de=\im^{Y(-S_0)}\sqrt{|S_0|}$$ 
(o\`u $Y$ est la fonction d'Heaviside) de telle sorte que $ \de^2=S_0$.
On d\'efinit   les fonctions $H$-invariantes $u$ et $v$ sur $\q$ par
$$u=\frac{Q+\de}{2}\qq\textrm{et}\qq
v=\frac{Q-\de}{2},$$
 de telle sorte que, pour  $X=\maq\in \q^{ss}$, alors
$Sp(YZ)=\{u(X),v(X)\}.$ 

\begin{lem}\label{orbites}
 \begin{enumerate}
\item L'application  $$\left\{\begin{array}{lll}\{H\textrm{-Orbites
semi-simples
de}\;\q\}&\longrightarrow&\qq\R^2\\\qq\qq\qq\mathcal{O}=H\cdot
X&\longmapsto&(Q(X),S(X))\end{array}\right.$$
est bijective. 
\item  Nous notons $\sim$ la relation d'\' equivalence sur $\C^2$ qui
identifie les couples $(x,y)$ et $(y,x)$. Alors l'application
$$\left\{\begin{array}{lll}\{H\textrm{-Orbites semi-simples
de}\;\q\}&\longrightarrow&\qq\C^2/\sim\\\qq\qq\qq\mathcal{O}=H\cdot
X&\longmapsto&(u(X),v(X))\end{array}\right.$$ est injective et
d'image $\big(\R^2\cup\{(\la,\bar{\la}); \la\in\C-\R\}\big)/\sim$.
\end{enumerate}
\end{lem}

\begin{dem}  La description de  $\a_{+,+},\a_{+,-},\a_{-,-}$ et
$\a_2$ permet d'obtenir l'image des deux applications consid\' er\' ees. Par ailleurs, toute  $H$- orbite semi-simple rencontre $\a\in<car(\q)>$. L'\' etude de $W_H(\a)$ pour $\a\in<car(\q)>$ donn\' ee dans le lemme
\ref{WH} permet d'\' etablir l'injectivit\' e de ces applications.\end{dem}

\begin{rem}\label{rempoly} \begin{enumerate}   \item Soit $X\in\q^{ss}$. Alors, les fonctions $u$ et $v$ prennent des valeurs r\' elles sur $\a_{\ep_1\ep_2}$ pour $\ep_1,\ep_2=\pm$ et on a 
\begin{enumerate} \item $H\cdot X\cap \a_{+,+}\neq\emptyset$ si  et seulement si $u(X)\geq0$ et $v(X)\geq0$,
\item $H\cdot X\cap\a_{+,-}\neq\emptyset$ si  et seulement si  $u(X)v(X)\leq0$,
\item $H\cdot X\cap\a_{-,-}\neq\emptyset$ si  et seulement si  $u(X)\leq0$ et $v(X)\leq0$,
\item $H\cdot X\cap\a_2\neq\emptyset$ si  et seulement si  $u(X)$ et $v(X)$ sont dans $\C$ et  $\overline{u(X)}=v(X)$.
\end{enumerate}
 \item Soit $\a\in car(\q)$ et $X\in \a$. 
On a les relations suivantes entre les racines de $\ra$, les  polyn\^omes $Q, S$ et  $S_0$   et les fonctions $u$ et $v$:
$$Q(X)=\dfrac{\al_1(X)^2+\al_2(X)^2}{4}=u(X)+v(X),$$
 $$ S(X)=\dfrac{\al_1(X)^2\al_2(X)^2}{16}=u(X)v(X),$$ 
 $$ S_0(X)=\be_1(X)^2\be_2(X)^2=(u(X)-v(X))^2.$$

$$\delta(X)=\left\{\begin{array}{ll} \vert \be_1(X)\be_2(X) \vert & \mbox{ si } \a  \mbox{ est $H$-conjugu\' e \` a  }\app, \apm  \mbox{ ou } \amm\\
\im \vert \be_1(X)\be_2(X) \vert & \mbox{ si } \a  \mbox{ est
$H$-conjugu\' e \` a  }\a_2\end{array}\right.$$

 \item $\qr=\{X\in\q, S(X)S_0(X)\ne0\}.$\end{enumerate}
\end{rem}

\subsection{Points semi-r\' eguliers}\label{semireg}

Les distributions propres invariantes et les int\' egrales orbitales
sont des fonctions $H$-invariantes de classe ${\cal C}^\infty$ sur
$\q^{reg}$. Pour \' etudier ces objets au voisinage des points de
$\q-\q^{reg}$, nous nous inspirons de la m\' ethode de descente
d'Harish-Chandra qui consiste \` a ramener cette \' etude \` a celle
d'objets de m\^eme type d\' efinis sur des espaces sym\' etriques
plus petits construits \` a partir des centralisateurs des   points
semi-r\' eguliers.\medskip

On rappelle  qu'un \' el\' ement semi-r\' egulier est un \' el\'
ement de $\q^{ss}$ qui se trouve dans l'intersection d'exactement
deux sous-espaces de Cartan. De fa\c con \' equivalente c'est un \'
el\' ement d'un sous-espace de Cartan qui annule exactement une
racine positive.\medskip

Ainsi, tout point semi-r\' egulier est conjugu\' e par $H$ \` a un
\' el\' ement de $\app\cap \apm$, $\apm\cap \amm$, $\app\cap \a_2$
ou $\amm\cap h_0\cdot \a_2$ avec $h_0=\left(
\begin{array}{c|c}
\begin{array}{cc} 0&1\\1&0\end{array}&0\\ \hline
0&\begin{array}{cc} 1&0\\0&-1\end{array}
\end{array}
\right)$.\medskip

Pour  $\a\in car(\q)$, on rappelle que  $\rap=\{\al_1, \al_2, \be_1=\dfrac{\al_1-\al_2}{2}, \be_2=\dfrac{\al_1+\al_2}{2}\}$ o\`u $\al_1$ et $\al_2$ sont fortement orthogonales et   de multiplicit\' e $1$  et $\be_1$, $\be_2$ sont  de multiplicit\' e $2$, orthogonales
mais non fortement orthogonales. Ces diff\' erences entre racines de
multiplicit\' e $1$ et racines de multiplicit\' e $2$ interviennent
dans l'\' etude des centralisateurs des points semi-r\' eguliers et
seront essentielles lors de l'\' etude des distributions propres
invariantes .

\subsubsection{Points semi-r\' eguliers annul\' es par une racine de
multiplicit\' e un}\label{semireg1}

Soit $H_1=\left(
\begin{array}{c|c}
0&\begin{array}{cc} 1&0\\0&0\end{array}\\ \hline
\begin{array}{cc} 1&0\\0&0\end{array}&0
\end{array}
\right)=\frac{1}{2} h_{\al_1}$  et $H_2=\ka(H_1)=\frac{1}{2}
h_{\al_2}$ les  coracines respectivement de $\al_1$ et $\al_2$ sur $\app$.
Ainsi, on a  $\R H_1=\app\cap\apm, \qq\R\vi(H_1)=\amm\cap\amp\qq\R\ka(H_1)=\app\cap\amp, \textrm{ et }\qq\R\vi\circ\ka(H_1)=\amm\cap\apm.$

Soit  $\z_1=\z_{\g}(H_1)=\left\{\left(\begin{array}{c|c}
\begin{array}{cc}a&0\\0&b\end{array}&\begin{array}{cc}x&0\\
0&y\end{array}\\ \hline
\begin{array}{cc}x&0\\
0&z\end{array}&\begin{array}{cc}a&0\\0&c\end{array}\end{array}\right);
a,b,c,x,y,z\in\R\right\}$ et $ \z_2=\z_{\g}(H_2)=\ka(\z_1).$

Les alg\` ebre $\z_1$ et $\z_2 $ sont  r\' eductives  et
$\sigma$-stables. On note $\ce_j$ le centre de $\z_j$.
On a alors  $\ce_j=\ce_j\cap\h\oplus\R H_j$ et 
$\ce_1\cap\h=\ce_2\cap\h= \left\{\left(
\begin{array}{c|c}
\begin{array}{cc}a&0\\0&b\end{array}&\begin{array}{cc} 0&0\\0&0\end{array}\\ \hline
\begin{array}{cc}
0&0\\0&0\end{array}&\begin{array}{cc}a&0\\0&b\end{array}
\end{array}\right); a,b\in\R\right\}.$ Pour conserver la sym\' etrie des espaces, on pose
$A_1=\left(
\begin{array}{c|c}
\begin{array}{cc}0&0\\0&1\end{array}&\begin{array}{cc} 0&0\\0&0\end{array}\\ \hline
\begin{array}{cc}
0&0\\0&0\end{array}&\begin{array}{cc}0&0\\0&1\end{array}
\end{array}\right)$ de telle sorte que $\ce_1\cap\h=\R A_1\oplus \R \kappa(A_1)$ et on note $\m_1=\R A_1\oplus[\z_1,\z_1]$ et $\m_2=\kappa(\m_1)$. 

On d\' efinit alors l'espace $\m$ par
 $$\m=\z_1+\z_2=\ce_1\cap\h\oplus[\z_1,\z_1]\oplus[\z_2,\z_2]=\m_1\oplus\m_2.$$

\begin{lem}\begin {enumerate}\item Les alg\` ebres $\m_1$ et $\m_2$ sont r\' eductives, $\sigma$-stables et $[\m_1,\m_2]=\{0\}$,

\item $\m$ est une sous-alg\` ebre de Lie $\sigma$-stable de $\g$ de dimension $8$. La paire sym\' etrique $(\m,\m\cap\h)$ est de rang $2$ et c'est le produit des deux paires sym\' etriques de rang un $(\m_1,\m_1\cap\h)$ et $(\m_2,\m_2\cap\h)$, chacune
isomorphe \` a $(\mathfrak{gl}(2,\R), \mathfrak{gl}(1,\R)\times
\mathfrak{gl}(1,\R))$, que l'on peut permuter \`a l'aide de $\ka$.

\end{enumerate}
\end{lem}

\no \begin{dem}  Les alg\` ebres $\z_1$ et $\z_2$ sont  r\' eductives et
$\sigma$-stables et les \' el\' ements $A_1$ et $A_2$ sont centraux
dans $\m$ donc $\m_1$ et $\m_2$ sont  r\' eductives et
$\sigma$-stables. Les racines de $\app$ \' etant r\' eelles, on a
$\z_1=\z_\g(\a_{++})\oplus(\g_{\al_2}^\C\cap\g\oplus\g_{-\al_2}^\C\cap\g)\mbox
{et } [\z_1,\z_1]\subset
(\g_{\al_2}^\C\cap\g\oplus\g_{-\al_2}^\C\cap\g)$.  Les racines
$\al_1$ et $\al_2$ sont fortement orthogonales et  les \' el\'
ements $A_1$ et $\ka(A_1)$ sont centraux dans $\m$, ainsi on obtient
$[\m_1,\m_2]=\{0\}$.

Comme  $[\z_1,\z_1]=\big\{\left(
\begin{array}{c|c}
\begin{array}{cc}0&0\\0&-c\end{array}&\begin{array}{cc} 0&0\\0&y\end{array}\\ \hline
\begin{array}{cc}
0&0\\0&z\end{array}&\begin{array}{cc}0&0\\0&c\end{array}
\end{array}
\right), y, z, c\in\R\big\}$, la deuxi\`eme  assertion est imm\' ediate. \end{dem}

Pour $(a, b, c, d)\in\R^4$, on note $diag(a, b, c, d)= \left (\begin{array}{c|c} \begin{array}{cc} a& 0\\0& b\end{array}&0\\
\hline 0&\begin{array}{cc}c & 0\\0 & d\end{array}\end{array}\right).
$

On rappelle que $\ka= Ad(K)$ avec $K= \left (\begin{array}{c|c} \begin{array}{cc} 0&1\\1&0\end{array}&0\\
\hline 0&\begin{array}{cc} 0&1\\1&0\end{array}\end{array}\right)$

\begin{lem}\label{ceno} Soit $N=N_H(\mq)$.
\begin{enumerate}
\item On a $N=N^0\cup KN^0$ avec $N^0=\{ diag(a, b, c, d); a, b, c, d\in \R^*\}$.
\item Soit $h\in H$ et $X\in\mq$ tels que $S_0(X)\neq 0$ et $h\cdot X\in\mq$, alors $h\in N$. ( Le polyn\^ome $S_0\in\C[\q]^H$ est d\' efini dans le paragraphe \ref{OS}).
\item Si $X\in\mq^{reg}$ alors $Z_H(X)=\{ diag(\al, \be, \al, \be); \al, \be\in\R^*\}.$
\end{enumerate}
\end{lem}

\no \begin{dem} Pour $D\in GL(4,\R)$ une matrice diagonale, on a imm\' ediatement $D\in N$ et $KD\in N$.\\

Soient $h=\left (\begin{array}{c|c} A&0\\
\hline 0&B\end{array}\right)\in H$  et  $X\in\mq$ tels que  $h\cdot X=X'\in\mq$.  On \' ecrit $X=\left(
\begin{array}{c|c}
0&\begin{array}{cc} x&0\\0&y\end{array}\\ \hline
\begin{array}{cc} z&0\\0&t\end{array}&0
\end{array}
\right)$ et  $X'=\left(
\begin{array}{c|c}
0&\begin{array}{cc} x'&0\\0&y'\end{array}\\ \hline
\begin{array}{cc} z'&0\\0&t'\end{array}&0
\end{array}
\right).$ Ainsi, on a
$A\left(\begin{array}{cc}
x&0\\0&y\end{array}\right)B^{-1}=\left(\begin{array}{cc}
x'&0\\0&y'\end{array}\right)$ et $
B\left(\begin{array}{cc}
z&0\\0&t\end{array}\right)A^{-1}=\left(\begin{array}{cc}
z'&0\\0&t'\end{array}\right), $ ce qui donne $A\left(\begin{array}{cc}
xz&0\\0&yt\end{array}\right)A^{-1}=B\left(\begin{array}{cc}
xz&0\\0&yt\end{array}\right)B^{-1}=\left(\begin{array}{cc}
x'z'&0\\0&y't'\end{array}\right).$\\

Supposons   $S_0(X)\neq 0$. On a alors  $xz\ne yt$, et par suite $x'z'\ne y't'$. Si  $xz=x'z'$ et $yt=y't'$ alors  $A$ et $B$ sont
diagonales, donc $h$ est diagonale. Si $xz=y't'$ et $yt=x'z'$, alors
$\left(\begin{array}{cc} 0&1\\1&0\end{array}\right)A$ et
$\left(\begin{array}{cc} 0&1\\1&0\end{array}\right)B$ sont
diagonales, donc $Kh$ est diagonale.
On obtient ainsi les assertions  {\it 1}. et {\it 2}.

Si $X\in\m\cap\q^{reg}$ et $h.X=X$ alors $S_0(X)\neq 0$ et  $h\in N^0$ par ce qui pr\' ec\`ede. On
en d\' eduit facilement l'assertion {\it 3}. \end{dem}

\begin{lem}
La famille $<car(\mq)>=\{\a_{+,+}, \a_{+,-}, \a_{-,-}\}$ est
une famille repr\' esentative des classes de conjugaison sous $N$
des sous-espaces de Cartan de $\m\cap\q$.\end{lem}

\no \begin{dem} Par d\' efinition de $\m\cap\q$, les sous -espaces $\app$, $\apm$ et $\amm$ sont des sous-espaces de Cartan de $\m\cap\q$.  Si $\a$ est un sous-espace de Cartan de $\mq$, c'est un sous-espace de Cartan de $\q$ et donc, il 
existe $h\in H$ tel que
$h\cdot\a\in<car(\q)>.$ Maintenant, par le  point $\it{2.}$ de la remarque \ref{rempoly},
le polyn\^ome $S_0\in\C[q]^H$ ne prend que des valeurs n\' egatives sur $\a_2$.  Un \' el\' ement de $\m\cap\q$ s'\' ecrit $X=\left(
\begin{array}{c|c}
0&\begin{array}{cc} x&0\\0&y\end{array}\\ \hline
\begin{array}{cc} z&0\\0&t\end{array}&0
\end{array}
\right)$ et par suite  $S_0(X)=(xz-yt)^2\geq 0$. On  obtient  donc $h\cdot \a\in \{\a_{+,+}, \a_{+,-}, \a_{-,-}\}$. L'assertion $\it{2.}$ 
du lemme \ref{ceno}  donne alors $h\in N$ ce qui prouve le lemme.\end{dem}

\subsubsection{El\' ements semi-r\' eguliers annul\' es par une racine de multiplicit\' e deux}\label{semireg2}
On note $H_3=\left(
\begin{array}{c|c}
0&\begin{array}{cc} 1&0\\0&1\end{array}\\ \hline
\begin{array}{cc} 1&0\\0&1\end{array}&0
\end{array}
\right)$ la coracine de $\be_2$ dans $\app$.  On a  alors $\R H_3=\app\cap\a_2$
et
$\R\vi(H_3)=\amm\cap\vi(\a_2)=\amm\cap h_0\cdot\a_{2}.$
Ainsi, tout \' el\' ement semi-r\' egulier annul\' e par une racine de multiplicit\' e $2$ est $H$-conjugu\' e \` a un \' el\' ement de $\R H_3$ ou de $\R\varpi(H_3)$.\me

 On pose
$\z_3=\z_\g(H_3)=\left\{\left(\begin{array}{c|c} A&B\\ \hline
B&A\end{array}\right); A,B\in\g\l(2,\R)\right\}\textrm{  et  }
\z_4=\z_\g(\vi(H_3))=\left\{\left(\begin{array}{c|c}
A&B\\ \hline -B&A\end{array}\right); A,B\in\g\l(2,\R)\right\}.$
La paire sym\' etrique
$(\z_3,\z_3\cap\h)$ est donc isomorphe \`a  $\Big(\mathfrak{gl}(2,\R)\times\mathfrak{gl}(2,\R),diag(\mathfrak{gl}(2,\R)\times\mathfrak{gl}(2,\R))\Big)\simeq
\mathfrak{gl}(2,\R)$ par l'application $\left (\begin{array}{c|c}
A&Y\\ \hline Y&A\end{array}\right )\rightarrow (A+Y,A-Y)$ et la paire sym\' etrique
$(\z_4, \z_4\cap\h)$  est isomorphe ˆ
$(\mathfrak{gl}(2,\C),\mathfrak{gl}(2,\R))$ par l'application
$\left (\begin{array}{c|c} A&Y\\ \hline -Y&A\end{array}\right
)\rightarrow A+\im Y$.\medskip

Comme $\q=\z_3\cap\q\oplus\z_4\cap\q$, on ne peut pas effectuer de r\' eduction \`a un espace
plus petit comme au paragraphe pr\' ec\' edent. Bien que  les  espaces
vectoriels $\z_3\cap\q$ et $\z_4\cap\q$ ne soient pas conjugu\' es
sous $H$,  on remarque que $\varpi(\z_3\cap\q)=\z_4\cap\q.$ Ainsi
l'\' etude sur  $\z_3\cap\q$ suffira pour comprendre les ph\'
enom\`enes au voisinage de tous les \' el\' ements semi-r\' eguliers
annul\' es par une racine de multiplicit\' e deux.\medskip

On note $\tilde{K}=\left(\begin{array}{c|c} I_2& 0\\ \hline 0 &
-I_2\end{array}\right)$ et $\tilde{\ka}=Ad(\tilde{K})$ de telle
sorte que $\tilde{\ka}|\q=-Id|\q$.

\begin{lem}\label{ceno3} Soit $N_3=N_H(\z_3\cap q)$.
\begin{enumerate}

\item On a $N_3=N_3^0\cup \tilde{K}N_3^0$ avec $N_3^0=\left\{ \left(\begin{array}{c|c} A&0\\ \hline
0&A\end{array}\right ); A\in GL(2,\R)\right\}.$
\item La famille $<car(\z_3\cap\q)>=\{\app,\a_2\}$ est une famille repr\' esentative des classes de conjugaison sous $N_3$ des sous-espaces de Cartan de $\z_3\cap\q$.
\end{enumerate}
\end{lem}

\no \begin{dem} Si $h= \left(\begin{array}{c|c} A&0\\ \hline
0&A\end{array}\right )$ alors on a imm\' ediatement $h\in N_3$ et
$\tilde{K} h\in N_3$.

Si $h= \left(\begin{array}{c|c} A&0\\ \hline
0&B\end{array}\right)\in N_3$ alors pour tout $Y\in GL(2,\R)$, on a $AYB^{-1}=BYA^{-1}$ et donc la matrice $P=B^{-1}A$ v\' erifie $PYP=Y$ pour tout $Y\in GL(2,\R)$. En appliquant ceci \` a $Y=\left(\begin{array}{cc} 1&0\\
0&0\end{array}\right )$, puis \` a $Y=\left(\begin{array}{cc} 0&0\\
0&1\end{array}\right )$ on obtient $P=\left(\begin{array}{cc} \ep_1&0\\
0&\ep_2\end{array}\right )$ avec $\ep_j=\pm 1$. Comme $A=BP$, la
relation  $AYB^{-1}=BYA^{-1}$ pour tout $Y\in GL(2,\R)$  donne l'assertion {\it 1}.

Soit $\a$ un sous-espace de Cartan de $\z_3\cap \q$.  Il existe
$h\in H$ tel que $h\cdot \a\in<car(\q)>$. Comme le polyn\^ome $S$
prend des valeurs positives sur $\z_3\cap\q$ et des valeurs
n\'egatives sur $\apm$, on a $h\cdot\a\neq\apm$. Maintenant, si
$h\cdot\a=\amm$, par un raisonnement analogue \`a ce qui
pr\'ec\`ede, il existerait $P\in GL(2,\R)$ telle que pour toute
matrice diagonale $D$ l'on ait $PDP=-D$ ce qui n'est pas possible.
Ainsi on a  $h.\a\in\{\app,\a_2\}$. Par suite  $\a$ et $h.\a$ sont
deux sous-espaces de Cartan  $H$- conjugu\' es de $\z_3\cap\q$, il
existe donc  $h_1\in N_3$ tel que $h_1\cdot\a=h\cdot\a$. On a alors
$h_1h^{-1}\in N_H(h\cdot \a)$ et donc par le lemme \ref{WH} on
obtient $h\in N_3$.\end{dem}

 \section{G\' en\' eralisation de certains r\' esultats
sur les int\' egrales orbitales en rang un}

Nous rappelons tout d'abord certains  r\' esultats d'Harish-Chandra
concernant l'int\' egration sur les fibres dont on peut trouver la preuve dans  le lemme 1 du
chapitre 2 premi\` ere partie de \cite{va} ou dans le chapitre 3
page 192 de \cite{fa}.

\begin{prop}\label{hc}
Soient $M$ et $N$ deux vari\' et\' es $\mathcal{C}^\infty$ de
dimension $m$ et $n$. Soit $\psi : M\to N$ une submersion de classe
$\mathcal{C}^\infty$ surjective de $M$ sur $N$.  Soient $\omega_M$
(resp. $\omega_N$) une $m$-forme (resp. $n$-forme) de classe
${\mathcal C}^\infty $ sur $M$ (resp. $N$) non nulle en tout point.
On note $\mu_M$ et  $\mu_N$ les mesures de Radon positives
respectivement sur $M$et  $N$ associ\' ees . Alors, on a:
\begin{enumerate}
\item 
Pour  tout $f\in \mathcal{D}(M)$, il existe une   unique  fonction $\psi_\ast(f)\in\mathcal{D}(N)$ telle que
\begin{equation}\label{rel}\int_Mf\cdot F\circ\psi d\mu_M=\int_N\psi_\ast(f)Fd\mu_N \qquad
,\forall F\in\mathcal{D}(N).\end{equation}

La fonction $\psi_*(f)$ est d\' efinie comme suit:  pour tout  $y\in N$, il existe une unique mesure de radon
positive $\mu_y$ telle que
$\psi_*(f)(y)=\int_{\psi^{-1}(y)}f d\mu_y.$
Ceci est appel\' e \textsf{l'int\' egration sur les fibres}.
\item L'application $\psi_\ast$ est une application lin\' eaire
continue surjective de $\D(M)$ dans
$\D(N)$.
\item La relation (\ref{rel}) est vraie si $F$ est localement
int\' egrable pour $\mu_N$.
\item L'application $\psi_\ast$ se prolonge aux fonctions
$\mu_M$-int\' egrables $f$ et dans ce cas $\psi_\ast(f)$ est d\'
efinie presque partout et $\mu_N$-int\' egrable. La
relation (\ref{rel}) est encore vraie.\\
De plus si $f$ est \` a support compact, alors $\psi_\ast(f)$ est \`a support compact  et de
$\mathcal{C}^\infty$ sur tout ouvert $\Omega$ tel que $f$ est
$\mathcal{C}^\infty$ sur $\psi^{-1}(\Omega)$.
\end{enumerate}
\end{prop}

\subsection{Cas d'une forme quadratique }\label{q1}

Soient $p$ et $q$ deux entiers de $\N^*$ et  $n=p+q$.  Soit $\mathcal{Q}$ la forme
quadratique d\' efinie sur $\R^{n}$ par
$$\mathcal{Q}(y_1,...,y_n)=\sum_{i=1}^p y_i^2-\sum_{i=1}^q y_{i+p}^2.$$

L'application $\mathcal{Q} $  est submersive et surjective de $\R^n-\{0\}$ dans $\R$ et la proposition pr\'ec\'edente 
 permet de d\' efinir l'application $\mathcal{Q}_*: \D( \R^n-\{0\})\rightarrow \D(\R)$.
 Cette application \' etudi\' ee par A. Tengstrand dans \cite{te} est appel\' ee fonction moyenne par J.~Faraut
 dans \cite{fa2} et not\' ee   $M_\mathcal{Q}$. Dans toute la suite, nous noterons $\mathcal{Q}_*=M_\mathcal{Q}.$\me

 D'apr\` es le  point ${\it 4.}$ de la proposition \ref{hc}, comme $0$ est une valeur critique de ${\mathcal Q}$, l'application $M_\mathcal{Q}$ se prolonge en une application continue, not\' ee encore $M_\mathcal{Q}$, de $\D(\R^n)$ dans $\mathcal{C}^{\infty}(\R^*)\cap L^1(\R)$ et pour $f\in\D(\R^n)$, la fonction $M_\mathcal{Q}(f)$ satisfait  toujours  la relation (\ref{rel}) qui s'\' ecrit,  
\begin{equation}\label{relbis}\int_{\R^n}F\circ \mathcal{Q}(y)f(y)dy=\int_\R
F(t)M_\mathcal{Q}f(t)dt \qq \textrm{ pour tout
 } F\in\D(\R).\end{equation}
 
Nous rappelons succinctement  les   propri\' et\' es essentielles de $M_\mathcal{Q}(f)$. On pourra se r\' ef\' erer \` a l'appendice de \cite{fa2} ou au chapitre 4 pages 204 et 205 de \cite{fa} pour  les d\' etails.\me

Soit $f\in\D(\R^n)$.   Le comportement de $M_\mathcal{Q}(f)(t)$ au voisinage de $t=0$ se
d\' ecrit \` a l'aide d'une fonction $\eta$ d\' efinie  sur $\R^*$,
qui d\' epend de la signature $(p,q)$ de $\mathcal{Q}$. Cette
fonction,  appel\' ee "fonction singularit\' e" relative \` a
$\mathcal{Q}$, est d\' efinie de la mani\` ere suivante (la fonction
d'Heaviside est not\' ee $Y$):\me

$\qq$ Si $p$ est impair et $q$ est pair,
$\eta(t)=Y(t)t^{\frac{n}{2}-1}$\me

$\qq$ Si $p$ est pair et $q$ est impair,
$\eta(t)=Y(-t)(-t)^{\frac{n}{2}-1}$\me

$\qq$ Si $p$ et $q$ sont pairs,
$\eta(t)=\frac{1}{2}sgn(t)t^{\frac{n}{2}-1}$\me

$\qq$ Si $p$ et $q$ sont impairs,
$\eta(t)=t^{\frac{n}{2}-1}\log|t|.$\me

Suivant (\cite{fa} chapitre 4 pages 204 et 205), on introduit
l'espace

$$\mathcal{H}_\eta=\{t\in\R^*\mapsto
\phi_0(t)+\eta(t)\phi_1(t);\;\phi_0,\phi_1\in\D(\R)\}.$$ 
Gr\^ace au  lemme de Borel, une fonction $\varphi$ appartient \`a $\mathcal{H}_\eta$ si et seulement si $\varphi$ est  de classe ${\mathcal C}^\infty$ sur $\R^*$, \`a support compact et s'il existe une suite de nombres  $\big( B_k(\varphi)\big)_{k\in\N}$ telle que, pour tout $N\in\N$ et $m<N+(n/2)$, la fonction $\varphi(t)-\eta(t)\sum_{k=0}^N t^kB_k(\varphi)$ soit de classe ${\mathcal C}^m$ sur $\R$.

Pour $a>0$, on note $\mathcal{H}_{\eta,a}=\{ \varphi\in\mathcal{H}_\eta;  support(\varphi)\subset [-a,a]$\}. Soit $\chi\in{\mathcal D}(\R)$ \' egale \`a $1$ au voisinage de $0$. On consid\`ere, pour $N\in\N$ et $m<N+(n/2)$ la semi-norme 
$$\Vert \varphi\Vert_{N,m}=\sup_{t}\arrowvert\big(\frac{d}{dt}\big)^m\big[\varphi(t)-\chi(t)\eta(t)\sum_{k=0}^N B_k(\varphi) t^k\big]\big)\arrowvert.$$ Les semi-normes $\Vert \varphi\Vert_{N,m}$, pour $m<N+(n/2)$, et $\vert B_k(\varphi)\vert $ pour $k\in\N$ munissent  $\mathcal{H}_{\eta,a}$ d'une topologie d'espace de Fr\' echet  et $\mathcal{H}_{\eta}$
est muni de la topologie limite inductive.

\begin{theo}\label{faraut}(\cite{fa}, paragraphe 4,  th\' eor\` eme page 205)  L'application $M_{\mathcal{Q}}$ est une surjection continue de
$\D(\R^n)$ dans $\mathcal{H}_\eta.$
\end{theo}

Soit $n_1, n_2\in\N$ et $f$ une fonction d\' efinie sur
$\R^{n_1}\times\R^{n_2}$. Pour $x\in\R^{n_1}$, on note $f_x$ la
 fonction d\' efinie sur $\R^{n_2}$ par $f_x(y)=f(x,y)$ et pour
$y\in\R^{n_2}$, on note $f^y$ la
 fonction d\' efinie sur $\R^{n_1}$ par $f^y(x)=f(x,y).$

Soit $m\in\N$. Ainsi pour $f\in\D(\R^m\times\R^n),$ il existe  deux
fonctions $\phi_0$ et $\phi_1$ d\' efinies sur $\R^m\times\R$,
telles que pour $x\in\R^m$, on ait $(\phi_0)_x$ et $(\phi_1)_x$ dans
$\D(\R)$ et
$$M_\mathcal{Q}f_x(t)=\phi_0(x,t)+\eta(t)\phi_1(x,t),\qq\forall t\in\R^*.$$
Le but de ce paragraphe  est de montrer que $\phi_0$ et $\phi_1$
sont dans $\D(\R^m\times\R).$\ Pour cela, nous reprenons les arguments d\' evelopp\' es dans le paragraphe  4 de  \cite{fa}  pour \' etudier la r\' egularit\' e en $x$ des fonctions  $(\phi_0)_x$ et $(\phi_1)_x$.

 On introduit la submersion surjective 
$G$ de $\R^m\times(\R^n-\{0\})$ dans $\R^m\times\R$ d\' efinie par 
$G(x,y)=(x,\mathcal{Q}(y)).$

\begin{lem}\label{regu}

 \begin{enumerate}
 \item L'application $G_*$ se prolonge aux fonctions de $\D(\R^m\times\R^n)$. Pour $f\in \D(\R^m\times\R^n)$ et  $(x,t)\in \R^m\times\R^*$, on a $G_*(f)(x,t)=M_\mathcal{Q}f_x(t)$ et la fonction $G_*(f)$ est \`a support compact.
 \item Pour $i_1,i_2,...,i_m$
 dans $\N$, $x\in\R^m$  et $t\in\R^*$, on a
$$\frac{\partial^{i_{1}}}{\partial
x_{1}^{i_{1}}}...\frac{\partial^{i_{m}}}{\partial
x_{m}^{i_m}}M_\mathcal{Q}f_{x}(t)=M_\mathcal{Q}\left(\left(\frac{\partial^{i_{1}}}{\partial
x_{1}^{i_{1}}}...\frac{\partial^{i_{m}}}{\partial
x_{m}^{i_m}}f\right)_{x}\;\right)(t).$$
En particulier, la fonction $(x,t)\mapsto M_\mathcal{Q}f_{x}(t)$ est de classe ${\mathcal C}^\infty$ sur $\R^m\times\R^*$.
\end{enumerate}
\end{lem}
\begin{dem} La premi\` ere assertion d\' ecoule de la d\' efinition et des propri\' et\' es de $M_{\mathcal{Q}}$ rappel\' ees au d\' ebut de ce paragraphe.

Pour prouver  le deuxi\` eme point, il suffit de montrer que, pour
tout $i\in\{1,...,m\}$, on a
$\dis \fr{\partial x_i}M_\mathcal{Q}f_{x}=M_\mathcal{Q}\left(\big(\fr{\partial
x_i}f\big)_{x}\;\right).$  
Soit $e_i$ le i\` eme vecteur de la
base canonique de $\R^m$ et $x\in\R^m$. On introduit la suite de
fonctions
$f_k=k(f_{x+\frac{1}{k}e_i}-f_x)
\in \D(\R^n).$ Les supports des fonctions $f_k$ et $\left(\fr{\partial
x_i}f\right)_{x}$ sont inclus dans un m\^eme compact \`a savoir 
$supp(f)+[0,1]e_i$. D'autre part, pour $j_1,j_2,...,j_n$ dans $\N$,  l'in\' egalit\' e de Taylor donne

$$\left\|\frac{\partial^{j_{1}}}{\partial
y_{1}^{j_{1}}}...\frac{\partial^{j_{n}}}{\partial y_{n}^{j_n}}f_k-
\frac{\partial^{j_{1}}}{\partial
y_{1}^{j_{1}}}...\frac{\partial^{j_{n}}}{\partial
y_{n}^{j_n}}\left(\left(\fr{\partial
x_i}f\right)_{x}\right)\right\|_\infty\leq\frac{1}{2k}\left\|\frac{\partial^{j_{1}}}{\partial
y_{1}^{j_{1}}}...\frac{\partial^{j_{n}}}{\partial
y_{n}^{j_n}}\frac{\partial^2}{\partial x_i^2}f\right\|_\infty.$$
Ainsi ,  la suite de
fonctions $\{f_k\}_{k\in\N^*}$ converge vers $\left(\fr{\partial
x_i}f\right)_{x}$ pour la topologie usuelle de  $\D(\R^n)$.
 La continuit\' e de l'application $M_\mathcal{Q}$ (th\' eor\` eme \ref{faraut}) donne alors 
$$  \lim_{k\to+\infty}M_\mathcal{Q}f_k=M_\mathcal{Q}\left(\left(\fr{\partial
x_i}f\right)_{x}\right).$$

 Par ailleurs pour $t\in\R^*$, on a
$\dis M_\mathcal{Q}f_k(t)=k(M_\mathcal{Q}f_{x+\frac{1}{k}e_i}(t)-M_\mathcal{Q}f_x(t)), $ et donc 
$\dis \lim_{k\to+\infty}M_\mathcal{Q}f_k(t)=\fr{\partial
x_i}M_\mathcal{Q}f_{x}(t)$ ce qui donne  le r\' esultat voulu.\end{dem}

Le comportement de la fonction $(x,t)\mapsto M_\mathcal{Q}f_{x}(t)$
au voisinage de $t=0$ se d\' eduit du comportement asymptotique en
$+\infty$ de sa transform\' ee de Fourier  partielle $\mathcal{F}(M_\mathcal{Q}f_x) $ en $t$ dont la relation (\ref{relbis}) donne l'expression suivante:
\begin{displaymath}\label{fourier}\mathcal{F}(M_\mathcal{Q}f_x)(\la)=\int_{\R}e^{-2\imath \pi\la
t}M_\mathcal{Q}f_x(t)dt=\int_{\R^n}e^{-2\imath \pi\la
\mathcal{Q}(y)}f_x(y)dy.\end{displaymath}

%%%%%%%%%%%TAYLOR
\begin{lem}\label{taylor}
Soit $f$ dans $\mathcal{D}(\R^m\times\R^n)$. Pour tout  $N$ dans
$\N$, il existe une fonction $\rho_N(x,\lambda)$ d\' efinie sur
$\R^m \times \R^\ast$ et une constante strictement positive  $C_N$
ind\' ependante de $x$ et $\la$ telles que pour tout $x$ dans $\R^m$
et $\la$ dans $\R^\ast$, on ait $|\rho_N(x,\la)|\leq
\frac{C_{N}}{|\la|^{N+1}}$ et
$$\mathcal{F}(M_\mathcal{Q}f_x)(\la)=
\frac{e^{-\im\frac{\pi}{4}(p-q)sgn(\la)}}{|2\la|^\frac{n}{2}}\big(\sum_{k=0}^N\frac{1}{k!}
(\frac{1}{8\imath \pi\la})^k(\partial \mathcal{Q})^kf_x(0)+\rho_N(x,\la)\big).$$
\end{lem}
%%%%%%%%%%%%%%%%%
\begin{dem} Gr\^ace au lemme page 203 de \cite{fa} on a :
$$\mathcal{F}(M_\mathcal{Q}f_x)(\la)=\int_{\R^n}e^{-2\imath \pi\la
\mathcal{Q}(y)}f_x(y)dy=\frac{e^{-\im\frac{\pi}{4}(p-q)sgn(\la)}}{|2\la|^\frac{n}{2}}
\int_{\R^n}\mathcal{F}(f_x)(\xi)e^{\frac{\imath\pi
\mathcal{Q}(\xi)}{2\la}}d\xi.$$ Le d\' eveloppement  de Taylor
de la fonction $t\mapsto e^{\imath t}$ donne,  pour tout
$(\xi,\la)\in\R^n\times \R^*$, la relation 
$$\mathcal{F}(f_x)(\xi)e^{\frac{\imath\pi
\mathcal{Q}(\xi)}{2\la}}=\sum_{k=0}^N\frac{1}{k!}(\frac{1}{8\imath
\pi\la})^k(-4\pi^2 \mathcal{Q}(\xi))^k\mathcal{F}(f_x)(\xi)+R_N(\xi,\la)\mathcal{F}(f_x)(\xi),$$
avec $|R_N(\xi,\la)|\leq \frac{|\pi \mathcal{Q}(\xi)|^{N+1}}{|2\la|^{N+1}(N+1)!}$. \me

Comme $\dis \int_{\R^n}(-4\pi^2 \mathcal{Q}(\xi))^k\mathcal{F}(f_x)(\xi)d\xi=(\partial
\mathcal{Q})^{k}f_x(0),$ on obtient l'expression voulue en prenant
$\rho_N(x,\la)=\int_{\R^n}R_N(\xi,\la)\mathcal{F}(f_x)(\xi)d\xi$,
pour $x\in\R^m$ et $\la\in\R^*$.\me

Pour conclure, il suffit de  majorer
$\dis{\int_{\R^n}|\mathcal{Q}(\xi)|^{N+1}|\mathcal{F}(f_x)(\xi)|d\xi}$ ind\' ependamment de $x$. Soient $R, T\in\R^*_+$ tels que
$supp(f)\subset
\overline{B(0,R)}\times[-T;T]^n.$ 
Soit $\chi$ une fonction positive \` a support compact et \' egale \`a $1$ sur $[-T;T]^n$.
Soit $P_N$ le polyn\^ome d\' efini par
$P_N(\xi_1,\cdots,\xi_n)=(\xi_1^2+\cdots+\xi_n^2)^{N+1}(1+\xi_1^2)\cdots(1+\xi_n^2)$ de telle sorte que 
$P_N(\xi)\mathcal{F}(f_x)(\xi)=\mathcal{F}\left(\frac{1}{(2\im\pi)^n}\partial(P_N)f_x\right)(\xi)$ o\`u $\frac{1}{(2\im\pi)^n}\partial(P_N)=P_N\left(\frac{1}{2\im\pi}\frac{\partial}{\partial
y_1},\cdots,\frac{1}{2\im\pi}\frac{\partial}{\partial
y_n}\right)$. Comme,  pour tout $\xi$ dans $\R^n$ et pour tout $x$ dans
$\R^m$, on a
$$\left|\mathcal{F}\left(\frac{1}{(2\im\pi)^n}\partial(P_N)f_x\right)(\xi)\right|\leq\left\|\frac{1}{(2\im\pi)^n}\partial(P_N)f\right\|_\infty\int_{\R^n} \chi=A_N,$$

on obtient 
$$\int_{\R^n}|\mathcal{Q}(\xi)|^{N+1}|\mathcal{F}(f_x)(\xi)|d\xi\leq\int_{\R^n}\dfrac{P_N(\xi)|\mathcal{F}(f_x)(\xi)|
}{(1+\xi_1^2)\cdots(1+\xi_n^2)}d\xi\leq A_N\pi^n.$$

Cette majoration ach\`eve  la preuve .\end{dem}
Nous rappelons la
g\' en\' eralisation suivante du th\' eor\` eme de Borel.
\begin{lem}\label{borel} (\cite{ho0} Theorem 1.2.6) Soit $k\in\N^*$.
Soient $\{\phi_\al\}_{\al\in\mathbb{N}^k}$ des fonctions de
$\D(\mathbb{R}^m)$ dont les supports sont inclus dans un m\^eme
compact. Alors il existe $\phi$ dans $
\D(\mathbb{R}^m\times\mathbb{R}^k)$ telle que:
$$\frac{\partial^\al \phi}{\partial y^\al}(x,0)=\phi_{\al}(x), \textrm{ pour tout }
\al\in\mathbb{N}^k, x\in\mathbb{R}^m.$$
\end{lem}
\begin{dem} La preuve est donn\'ee pour $k=1$ dans le th\'eor\`eme 1.2.6 de \cite{ho0}. Le r\'esultat s'obtient ensuite par r\'ecurrence sur $k$.    \end{dem}

\begin{theo} \label{comp}
Pour toute fonction $f$ dans $\D(\R^m\times\R^n)$, il existe
$\phi_0$ et $\phi_1$ dans $\D(\R^m\times \R)$ telles que, pour $t$
dans $\R^\ast$ et $x$ dans $\R^m$, on ait
$$M_\mathcal{Q}f_x(t)=\phi_0(x,t)+\eta(t)\phi_1(x,t),$$
et \begin{equation}\label{tay}\frac{\partial^k\phi_1}{\partial
t^k}(x,0)=
\frac{c\pi^{\frac{n}{2}}}{4^k\Gamma(\frac{n}{2}+k)}(\partial
\mathcal{Q})^kf_x(0),\end{equation}

avec \\$\qq c=(-1)^{\frac{q}{2}}$ pour $q$
pair,\\
$\qq c=(-1)^{\frac{p}{2}}$ pour $p$ pair et $q$ impair\\
$\qq c=\frac{(-1)^{\frac{q+1}{2}}}{\pi}$ pour $p$ impair et $q$
impair.
\end{theo}
\begin{dem} Soit $f\in\D(\R^m\times\R^n)$. On pose $\phi_k(x)=
\frac{c\pi^{\frac{n}{2}}}{4^kk!\Gamma(\frac{n}{2}+k)}(\partial
\mathcal{Q})^kf_x(0),$ o\`u la constante  $c$ est donn\' ee  dans
l'\' enonc\' e du th\' eor\` eme. 

 Nous allons tout d'abord montrer que, pour tout $N\in\N$,  la fonction
$\dis(x,t)\mapsto M_\mathcal{Q}f_x(t)-\eta(t)\sum_{k=0}^N\phi_k(x)t^k$
peut se prolonger en une fonction de
$\mathcal{C}^N(\R^m\times\R).$ C'est une fonction de $\mathcal{C}^\infty(\R^m\times\R^*)$ par le lemme \ref{regu}. Ainsi, il suffit de prouver que, pour 
$\al\in \D(\R)$, \' egale \` a $1$ au voisinage de $0$, la fonction $\psi_N$ d\' efinie par 
$$\psi_N(x,t)=M_\mathcal{Q}f_x(t)-\eta(t)\al(t)\sum_{k=0}^N\phi_k(x)t^k$$
se prolonge en une fonction de
$\mathcal{C}^\infty(\R^m\times\R).$ Comme $(\psi_N)_x\in L^1(\R)$, nous   \' etudions le   comportement en l'infini de
sa transform\' ee de Fourier $\mathcal{F}((\psi_N)_x)$ pour  obtenir sa r\' egularit\' e en $t=0$.

Le lemme \ref{taylor} et le calcul de $\mathcal{F}(t^k\eta(t))$ au sens distribution (proposition page 206 de \cite{fa}), donnent  pour $\la\neq 0$, la relation
$$\mathcal{F}((\psi_N)_x)(\la)=\sum_{k=1}^N\mathcal{F}(t^k\eta(t)-t^k\alpha(t)\eta(t))(\la) \phi_k(x)+\frac{e^{-\im\frac{\pi}{4}(p-q)sgn(\la)}}{|2\la|^{\frac{n}{2}}}\rho_N(x,\la).$$

Maintenant, on remarque que, pour $l$ et $k$ deux
entiers tels que $l>k+\frac{n}{2}$, la fonction
$t\mapsto\frac{d^l}{dt^l}\big((1-\al(t))t^k\eta(t)\big)$ est une
fonction int\' egrable sur $\R$. Ainsi,  il existe un r\' eel
$A_N$ tel que, pour tout $k\in\{0,...,N\}$ et pour tout
$\la\in\R^*$, on ait
$$|\mathcal{F}\big((1-\al(t))t^k\eta(t)\big)(\la)|\leq\frac{A_N}{|\la|^{N+\frac{n}{2}+1}}.$$

Par le lemme \ref{taylor}, en posant
$E_N=C_N+A_N(\sum_{k=0}^N\|\phi_k\|_\infty)$, on obtient,   pour tout $x\in\R^m$ et $\la\in\R^*$:
\begin{equation}\label{ineg}|\mathcal{F}((\psi_N)_x)(\la)|\leq\frac{E_N}{|\la|^{N+\frac{n}{2}+1}}.\end{equation}

Soient $j_1,...j_m$ des entiers positifs. On applique ce qui pr\'
ec\` ede en rempla\c cant $f_x$ par
$\frac{\partial^{j_{1}}}{\partial
x_{1}^{j_{1}}}...\frac{\partial^{j_{m}}}{\partial
x_{m}^{j_m}}f_{x}$. On obtient alors qu'il existe un r\' eel
$E_N^{j_1,...j_m}$ ind\' ependant de $x$ et $t$ tel que pour $x$
dans $\R^m$ et $t$ dans $\R^\ast$, on ait
\begin{equation}\label{maj}
|\mathcal{F}(\frac{\partial^{j_{1}}}{\partial
x_{1}^{j_{1}}}...\frac{\partial^{j_{m}}}{\partial
x_{m}^{j_m}}(\psi_N)_x)(\la)|\leq\frac{E_N^{j_1,...j_m}}{|\la|^{N+\frac{n}{2}+1}}.\end{equation}
Ainsi, par les propri\' et\' es de la transformation de Fourier, la
fonction $\psi_N$ se prolonge en une fonction de
$\mathcal{C}^N(\R^m\times\R)$. Donc la fonction $(x,t)\mapsto
M_\mathcal{Q}f_x(t)-\eta(t)\sum_{k=0}^N\phi_k(x)t^k$ est aussi dans
$\mathcal{C}^N(\R^m\times\R).$\me

Maintenant concluons. \me

Par d\' efinition  des  $\varphi_k(x)\in \D(\R^m)$,  ces fonctions ont leur support
inclus dans un m\^eme compact. Par le lemme \ref{borel}, il existe une fonction $\phi$ de
$\mathcal{D}(\R^m\times\R)$ telle que, pour tout $x\in\R^m$
$$\frac{\partial^k\phi}{\partial t^k}(x,0)=k!\phi_{k}(x).$$

Soit $N$ un entier positif. Pour tout $x$ dans $\R^m$ et $t$ dans
$\R^\ast$, on a :
$$M_\mathcal{Q}f_x(t)-\eta(t)\phi(x,t)$$$$=M_\mathcal{Q}f_x(t)-
\eta(t)\sum_{k=0}^N\frac{\partial^k\phi}{\partial
t^k}(x,0)\frac{t^k}{k!}-\eta(t)
t^{N+1}\int_0^1\frac{\partial^{N+1}\phi}{\partial
t^{N+1}}(x,tu)\frac{(1-u)^N}{N!}du$$ $$=\big(M_\mathcal{Q}f_x(t)-
\eta(t)\sum_{k=0}^N\phi_k(x)t^k\big)- \eta(t)
t^{N+1}\int_0^1\frac{\partial^{N+1}\phi}{\partial
t^{N+1}}(x,tu)\frac{(1-u)^N}{N!}du.$$
 Par ce qui pr\' ec\` ede, les deux termes de cette somme
sont dans $\mathcal{C}^N(\R^m\times\R)$.

 Ainsi la fonction $(x,t)\mapsto M_\mathcal{Q}f_x(t)-\eta(t)\phi(x,t)$ est dans
 $\mathcal{C}^N(\R^m\times\R)$ pour tout entier $N$, donc dans
$\mathcal{C}^\infty(\R^m\times\R)$, avec
$$\frac{\partial^k\phi}{\partial t^k}(x,0)=
\frac{c\pi^{\frac{n}{2}}}{4^k\Gamma(\frac{n}{2}+k)}(\partial
\mathcal{Q})^kf_x(0).$$ De plus gr\^ace au lemme \ref{regu}, la
fonction $(x,t)\mapsto M_\mathcal{Q}f_x(t)$ est \` a support born\'
e tout comme la fonction $(x,t)\mapsto\eta(t)\phi(x,t)$, donc la
fonction $(x,t)\mapsto M_\mathcal{Q}f_x(t)-\eta(t)\phi(x,t)$ est
dans
 $\D(\R^m\times\R).$
\end{dem}

Soit $\mathcal{H}^m_\eta$ l'espace des fonctions de $\R^m\times\R^*$
dans $\R$ d\' efini par
$$\mathcal{H}^m_\eta:=\{(x,t)\mapsto\phi_0(x,t)+\eta(t)\phi_1(x,t);\qq
\phi_0,\phi_1\in\D(\R^m\times\R)\}.$$

Comme pour  $\mathcal{H}_\eta$,  par le lemme de Borel g\' en\' eralis\' e, on a   $\varphi\in \mathcal{H}^m_\eta$ si et seulement si $\varphi$ est \`a support  compact, de classe $\mathcal{C}^\infty$ sur $\R^m\times\R^*$ et s'il existe une suite $(\varphi_k)_{k\in\N}$ de $\mathcal{D}(\R^m)$ dont les termes sont de support inclus dans un m\^eme compact telle que, pour tout $N\in\N$, la fonction $(x,t)\mapsto \varphi(x,t)-\eta(t)\sum_{k=0}^N \varphi_k(x) t^k$ soit de classe $\mathcal{C}^l$ sur $\R^m\times \R$ pour tout $l\in\N$ v\' erifiant $l<N+(n/2)$.

Pour $K\subset \R^{m+1}$ un compact, on note $\mathcal{H}_{\eta,K}=\{ \varphi\in\mathcal{H}_\eta; support(\varphi)\subset K\}$. Soit $\chi\in{\mathcal D}(\R)$ \' egale \`a $1$ au voisinage de $0$ et de support contenu dans $K$. On pose $ T_N(\varphi)(x,t)=\varphi(x,t)-\chi(t) \eta(t)\sum_{k=0}^N \varphi_k(x) t^k$. 

Pour $(N,l)\in\N^2$ avec  $l<N+(n/2)$, $k\in \N$  et $\mathbf{i}=(i_1,\ldots,i_m)\in\N^m$, on d\' efinit les semi-normes suivantes 
$$\Vert \varphi\Vert_{N,l,\mathbf{i}}=\sup_{(x,t)}\arrowvert\big(\frac{\partial}{\partial x_1}\big)^{i_1}\ldots \big(\frac{\partial}{\partial x_m}\big)^{i_m}\big(\frac{\partial}{\partial t}\big)^lT_N(\varphi)(x,t)\arrowvert$$
$$\quad\textrm{ et }\quad B_{k,\mathbf{i}}(\varphi)=\sup_x \arrowvert \big(\frac{\partial}{\partial x_1}\big)^{i_1}\ldots \big(\frac{\partial}{\partial x_m}\big)^{i_m}\varphi_k(x)\arrowvert.$$ 
Les semi-normes $\Vert \varphi\Vert_{N,l,\mathbf{i}}$ et $ B_{k,\mathbf{i}}(\varphi)$  munissent  $\mathcal{H}^m_{\eta,K}$ d'une topologie d'espace de Fr\' echet  et $\mathcal{H}^m_{\eta}$
est muni de la topologie limite inductive.

\begin{theo}\label{surj2}
L'application $$\left\{\begin{array}{lll}
\D(\R^m\times\R^n)&\to&\qq\qq\mathcal{H}^m_\eta\\
\qq\qq f&\mapsto&(x,t)\mapsto
M_\mathcal{Q}f_x(t)\end{array}\right.$$ est continue et surjective.
\end{theo}

\begin{dem} Montrons tout d'abord la surjectivit\' e de cette application. Soient $\phi_0,\phi_1$ dans $\D(\R^m\times\R)$. Par le lemme \ref{borel}, il existe  $g$ dans $\D(\R^m\times\R^n)$
telle que, pour $x\in\R^m$, l'on ait:
$$(\partial \mathcal{Q})^k g(x,0)=\frac{4^k\Gamma(\frac{n}{2}+k)}{c\pi^{\frac{n}{2}}}\frac{\partial^k\phi_1}{\partial
t^k}(x,0).$$ o\`u $c$ est la constante d\' ependant de la signature
de $\mathcal{Q}$ d\' efinie dans le th\' eor\` eme \ref{comp}.

Par la relation (\ref{tay}) du th\' eor\` eme  \ref{comp}
il existe  $\psi_0$ et $\psi_1$ dans $\D(\R^m\times\R)$ telles que pour $x$ dans $\R^m$ et $t\ne0$,
$$M_\mathcal{Q}g_x(t)=\psi_0(x,t)+\eta(t)\psi_1(x,t),\qq
\textrm{avec}\qq \frac{\partial^k\psi_1}{\partial
t^k}(x,0)=\frac{\partial^k\phi_1}{\partial t^k}(x,0).$$ 
Pour $(x,t)\in\R^m\times \R^*$, on
pose
$h(x,t)=\psi_1(x,t)-\phi_1(x,t)$   de telle sorte que

$$\frac{\partial^kh}{\partial t^k}(x,0)=0,\mbox{  pour tout entier positif } k $$
Ainsi la fonction $(x,t)\mapsto\eta(t)h(x,t)$ appartient \` a
 $\D(\R^m\times\R)$ et
donc la fonction $(x,t)\mapsto
M_\mathcal{Q}g_x(t)-\phi_0(x,t)-\eta(t)\phi_1(x,t)=
\psi_0(x,t)-\phi_0(x,t)+\eta(t)h_1(x,t).$$ $
se prolonge en une fonction de $\D(\R^m\times\R).$\me

Maintenant, l'application $G_* : \D(\R^m\times
(\R^n-\{0\})\rightarrow \D(\R^m\times\R)$ d\' efinie par
$G_*(f)(x,t)= M_{\mathcal{Q}}(f_x)(t)$ est surjective (lemme \ref{regu}). Il existe
donc $u$ dans $\D(\R^m\times(\R^n-\{0\}))$ telle que pour tout $x$
dans $\R^m$ et $t$ dans $\R$,
$$M_\mathcal{Q}u_x(t)=M_\mathcal{Q}g_x(t)-\phi_0(x,t)-\eta(t)\phi_1(x,t),$$
ce qui donne la surjectivit\' e.\me

 La continuit\' e s'obtient comme dans la situation $m=0$ (\cite{fa} page page 207) par le th\' eor\`eme du graphe ferm\' e . En effet, soit  $(f_k)_k$ une suite de $\D(\R^{n+m})$ qui converge vers $f$ dans $\D(\R^{n+m})$ et telle que $M_{\mathcal Q}( f_k)$ converge vers $\phi$ dans $\mathcal{H}^m_\eta$. Soit $\be\in\D(\R)$ nulle au voisinage de $0$. Alors la suite de fonctions $g_k(x,y)=\be\circ Q(y) f_k(x,y)$ converge vers $g(x,y) =\be\circ Q(y) f(x,y)$ dans $\D(\R^{n+m})$. Comme $g_k\in \D(\R^n\times\R^m-\{0\})$ , la continuit\' e de $G_*$ donne $\be(y)M_{\mathcal Q}( f_k)(x,y)$ converge vers $\be(y)M_{\mathcal Q}( f)(x,y)$ dans $\D(\R^n\times\R)$ et donc dans $\mathcal{H}^m_\eta$. Ceci \' etant valable pour toute fonction $\be$ comme ci-dessus, on obtient le r\' esultat voulu.\end{dem}

\subsection{Cas de deux formes quadratiques}\label{qq}

Nous allons g\' en\' eraliser les r\' esultats obtenus dans le paragraphe pr\' ec\' edent \` a la situation suivante. Pour $j\in\{1,2\}$, on fixe 
 $p_j$ et $q_j$ dans $\N^*$ et on pose
$n_j=p_j+q_j$. Soit $\mathcal{Q}_j$ la forme quadratique d\' efinie sur
$\R^{n_j}$ par
$$\mathcal{Q}_j(x_1,...,x_{n_j})=\sum_{i=1}^{p_j} x_i^2-\sum_{i=1}^{q_j} x_{i+p_j}^2.$$
 On consid\` ere l'application
$$H:\left\{\begin{array}{lll}\R^{n_1}\times\R^{n_2}&\to&\qq\qq\qq\R^2\\
\qq(x,y)&\to&(\mathcal{Q}_1(x),\mathcal{Q}_2(y))\end{array}\right..$$

La restriction $H|_{\R^{n_1}-\{0\}\times\R^{n_2}-\{0\}}$ est une
surjection submersive de $\R^{n_1}-\{0\}\times\R^{n_2}-\{0\}$ dans
$\R^2$. Par la proposition \ref{hc}, l'application $H_*$ se prolonge
en une application $M_{\mathcal{Q}_1,\mathcal{Q}_2}$ de
$\D(\R^{n_1}\times\R^{n_2})$ dans
$\mathcal{C}^\infty(\R^*\times\R^*)\cap L^1(\R^2)$. Pour $f\in
\D(\R^{n_1}\times\R^{n_2})$ la fonction
$M_{\mathcal{Q}_1,\mathcal{Q}_2}f$  v\' erifie,
 pour tout $F\in\D(\R^2)$, la relation suivante
\begin{equation}\label{rel2}\int_{\R^n_1\times\R^{n_2}}F\circ
H(x,y)f(x,y)dxdy=\int_{\R^2}
F(t_1,t_2)M_{\mathcal{Q}_1,\mathcal{Q}_2}f(t_1,t_2)dt_1dt_2.\end{equation}

Pour $i\in\{1,2\}$ on pose $\eta_i$ la "fonction singularit\' e"
relative \` a $\mathcal{Q}_i$.

\begin{lem}\label{double} Soit $f\in\D(\R^{n_1}\times\R^{n_2})$. Alors pour
$(t_1,t_2)\in\R^*\times\R^*$, on a
$$M_{\mathcal{Q}_1}((M_{\mathcal{Q}_2}f_x)^{t_2})(t_1)=M_{\mathcal{Q}_2}((M_{\mathcal{Q}_1}f^y)_{t_1})(t_2)
=M_{\mathcal{Q}_1,\mathcal{Q}_2}f(t_1,t_2).$$
\end{lem}

\begin{dem}  Par la relation (\ref{rel2}), pour tout
$F\in\D(\R^2)$, on a :
  $$
\int_{\R^2}F(t_1,t_2)M_{\mathcal{Q}_1,\mathcal{Q}_2}f(t_1,t_2)dt_1dt_2=\int_{\R^n_1\times\R^{n_2}}
F(\mathcal{Q}_1(x),\mathcal{Q}_2(y))f(x,y)dxdy=$$
$$\int_{\R^{n_1}}\int_\R F(\mathcal{Q}_1(x),t_2)M_{\mathcal{Q}_2}f_x(t_2)dt_2dx=
\int_{\R^2}F(t_1,t_2)M_{\mathcal{Q}_1}((M_{\mathcal{Q}_2}f_x)^{t_2})(t_1)dt_1dt_2,$$
Ces \' egalit\' es \' etant vraies pour toute fonction $F$ de
$\D(\R^2)$, on obtient
$$M_{\mathcal{Q}_1}((M_{\mathcal{Q}_2}f_x)^{t_2})(t_1)=M_{\mathcal{Q}_1,\mathcal{Q}_2}f(t_1,t_2).$$
De m\^eme, on montre que
$$M_{\mathcal{Q}_2}((M_{\mathcal{Q}_1}f^y)_{t_1})(t_2)
=M_{\mathcal{Q}_1,\mathcal{Q}_2}f(t_1,t_2).$$\end{dem}

\begin{theo}\label{comp2}
Pour $f\in\D(\R^{n_1}\times\R^{n_2})$, il existe des fonctions
$\psi_0, \psi_1, \psi_2$ et $\psi_3$ dans $\D(\R^2)$ telles que pour
$(t_1,t_2)\in\R^*\times\R^*$, on ait
$$M_{\mathcal{Q}_1,\mathcal{Q}_2}f(t_1,t_2)=\psi_0(t_1,t_2)+\eta_1(t_1)\psi_1(t_1,t_2)+
\eta_2(t_2)\psi_2(t_1,t_2)+\eta_1(t_1)\eta_2(t_2)\psi_3(t_1,t_2).$$

\end{theo}

\begin{dem}  Le th\' eor\` eme  \ref{comp} nous donne l'existence de
fonctions $\phi_0$ et $\phi_1$ dans $\D(\R^{n_1}\times\R)$ telles
que pour $(x,t_2)\in\R^{n_1}\times\R^*$, on ait
$\dis M_{\mathcal{Q}_2}f_x(t_2)=\phi_0(x,t_2)+\eta_2(t_2)\phi_1(x,t_2)$
et donc
$$ M_{\mathcal{Q}_1}((M_{\mathcal{Q}_2}f_x)^{t_2})(t_1)=M_{\mathcal{Q}_1}\big((\phi_0)^{t_2}\big)(t_1)+
\eta_2(t_2)M_{\mathcal{Q}_1}\big((\phi_1)^{t_2}\big)(t_1).$$
 Gr\^ace
au th\' eor\` eme  \ref{comp}, il existe des fonctions
$\phi_2,\phi_3,\phi_4$ et $\phi_5$ dans $\D(\R^2)$ telles que pour
$(t_1,t_2)\in\R^*\times\R$, on ait
$$M_{\mathcal{Q}_1}\big((\phi_0)^{t_2}\big)(t_1)=\phi_2(t_1,t_2)+\eta_1(t_1)\phi_3(t_1,t_2)$$
et
$$M_{\mathcal{Q}_1}\big((\phi_1)^{t_2}\big)(t_1)=\phi_4(t_1,t_2)+\eta_1(t_1)\phi_5(t_1,t_2).$$
Le lemme \ref{double} permet de conclure. \end{dem}

Soit $\mathcal{H}_{\eta_{1},\eta_2}$ le sous ensemble des fonctions
de $\R^*\times\R^*$ dans $\R$ d\' efini par
$$\mathcal{H}_{\eta_{1},\eta_2}=\{(t_1,t_2)\mapsto
\psi_0(t_1,t_2)+\eta_1(t_1)\psi_1(t_1,t_2)+
\eta_2(t_2)\psi_2(t_1,t_2)$$
$$+\eta_1(t_1)\eta_2(t_2)\psi_3(t_1,t_2);\qq
\psi_0,\psi_1,\psi_2,\psi_3\in\D(\R^2)\}.$$
Comme pr\' ec\' edemment,  on a $\varphi\in \mathcal{H}_{\eta_{1},\eta_2}$ si et seulement si $\varphi\in{\mathcal C}^\infty(\R^*\times\R^*)$  et \`a support compact et s'il existe des suites $\big( A_{k,l}(\varphi)\big)_{k,l\in\N}$, $\big( B_{k,l}(\varphi)\big)_{k,l\in\N}$ et $\big( C_{k,l}(\varphi)\big)_{k,l\in\N}$ telles que la fonction
 $\dis R_N(\varphi)(t_1,t_2) =\varphi(t_1,t_2)-\eta_1(t_1)\sum_{k,l\leq N} t_1^k t_2^l A_{k,l}(\varphi)$ $-\eta_2(t_2)\sum_{k,l\leq N} t_1^k t_2^l  B_{k,l}(\varphi)
  -\eta_1(t_1)\eta_2(t_2)\sum_{k,l\leq N} t_1^k t_2^l C_{k,l}(\varphi) $
  soit de classe ${\mathcal C}^m$ pour tout $m<N+\dfrac{n}{2}.$ 
 
 Pour $K$ une partie compacte de $\R^2$, on note $\mathcal{H}_{\eta_{1},\eta_2,K}$ le sous-espace de $\mathcal{H}_{\eta_{1},\eta_2}$ form\' e des fonctions \`a support dans $K$. On fixe $\chi_1$ et $\chi_2$ dans $\D(\R)$, \`a support contenu dans $K$ et \' egale \`a $1$ au voisinage de $0$. On pose
 $$\begin{array}{c}T_N(\varphi)(t_1,t_2) =\varphi(t_1,t_2)-\chi_1(t_1) \eta_1(t_1)\dis\sum_{k,l=0}^N t_1^k t_2^l A_{k,l}(\varphi)\\ -\chi_2(t_2) \eta_2(t_2)\dis\sum_{k,l=0}^N t_1^k t_2^l  B_{k,l}(\varphi)  -\chi_1(t_1)\chi_2(t_2) \eta_1(t_1)\eta_2(t_2)\dis\sum_{k,l\leq N} t_1^k t_2^l C_{k,l}(\varphi).\end{array}$$
 On introduit les  semi-normes suivantes $$\Vert \varphi\Vert_{l,k,N}=\sup_{(t_1,t_2)} \vert \big( \dfrac{\partial}{\partial t_1}\big)^k \big( \dfrac{\partial}{\partial t_2}\big)^l T_N(\varphi)(t_1,t_2)\vert \quad \textrm{ pour } l+k\leq N+(n/2).$$
Les semi-normes $\Vert \varphi\Vert_{l,k,N}$, $\vert A_{k,l}(\varphi)\vert$, $\vert B_{k,l}(\varphi)\vert$  et $\vert C_{k,l}(\varphi)\vert$   munissent $\mathcal{H}_{\eta_{1},\eta_2,K}$  d'une structure d'espace de Fr\' echet. L'espace  $\mathcal{H}_{\eta_{1},\eta_2}$  est munie de la topologie limite inductive.

\begin{prop}\label{surj3}
L'application $$\left\{\begin{array}{lll}
\D(\R^{n_1}\times\R^{n_2})&\to&\qq\qq\qq\mathcal{H}_{\eta_{1},\eta_2}\\
\qq\qq f&\mapsto&(t_1,t_2)\mapsto
M_{\mathcal{Q}_1,\mathcal{Q}_2}f(t_1,t_2)\end{array}\right.$$ est continue et 
surjective.
\end{prop}

\begin{dem} Montrons la surjectivit\' e. Soient $\psi_0,\psi_1,\psi_2,$ et $\psi_3\in\D(\R^2)$.
Gr\^ace au th\' eor\` eme \ref{surj2}, il existe $f$ et $g$ dans
$\mathcal{D}(\R^{n_1}\times\R)$ telles que pour
$(t_1,t_2)\in\R^*\times\R^*$, on ait
$$\begin{array}{lc}  & M_{\mathcal{Q}_1}(f^{t_2})(t_1)=\psi_0(t_1,t_2)+\eta_1(t_1)\psi_1(t_1,t_2)\\ 
\textrm{ et }\quad &  M_{\mathcal{Q}_1}(g^{t_2})(t_1)=\psi_2(t_1,t_2)+\eta_1(t_1)\psi_3(t_1,t_2).\end{array}$$
Toujours gr\^ace au th\' eor\` eme \ref{surj2}, il existe
$k\in\mathcal{D}(\R^{n_1}\times\R^{n_2})$ telle que pour
$(x,t_2)\in\R^{n_1}\times\R^*$, on ait
$M_{\mathcal{Q}_2}k_x(t_2)=f(x,t_2)+\eta_2(t_2)g(x,t_2).$
Le lemme \ref{double} donne alors la surjectivit\' e. La continuit\' e est imm\' ediate par le th\' eor\`eme  \ref{surj2} . \end{dem}

\section{Comportement des int\' egrales orbitales}

\subsection{Normalisation des int\' egrales orbitales}

L'int\' egrale orbitale d'une fonction $f\in\D(\q)$ est la fonction d\' efinie sur $\q^{reg} $  comme la moyenne de $f$ sur chaque  orbite $H\cdot X$  pour  $X\in\q^{reg}$, ceci pour une mesure $H$-invariante sur $H\cdot X$ que nous allons maintenant pr\' eciser.

Soit $p$  la dimension de $\q$. On note  $dZ$ la mesure $H$-invariante sur $\q$ d\' efinie par la densit\' e  $\mu$ donn\' ee par
$\mu(\xi_{1},...,\xi_{p})=|det(\om(\xi_{i},\xi_{j})_{i,j})|^{\frac{1}{2}}$.\me

 Pour $X\in\q^{reg}$ et $\a=Z_\q(X)$, on note $\dis\Pi=\big|\prod_{\al\in\rap}\al^{m_\al}\big|.$

 Par le lemme 1.20 et le th\' eor\` eme 1.27 de \cite{or}, il existe une unique mesure $H$-invariante $d\dot{h}$ sur $H/Z_H(X)$ normalis\' ee de telle sorte  que, pour tout $f\in D(\q)$, l'on ait
 $$\int_{H.\a}f(Z)\; dZ=\frac{1}{\vert W_H(\a)\vert}\int_{\a^{reg}}\big( \int_{H/Z_H(X)}f(h\cdot X)d\dot{h}\big)\Pi(X)\; dX.$$

 On rappelle que la fonction $\delta$, d\' efinie dans le paragraphe \ref{OS}, v\' erifie

 $$\vert \delta(X)\vert =\big|\prod_{\al\in\rap}\al^{m_\al-1}(X)\big| \quad\mbox{ et }\quad \delta^2=S_0. $$
 \begin{Def} Pour $f\in\D(\q)$, on d\' efinit son int\' egrale orbitale sur $\q^{reg}$ par
 $$\mathcal{M}_H(f)(X)=  |\delta(X)|  \int_{H/Z_H(X)}f(h\cdot X)d\dot{h}.$$
 \end{Def}

La normalisation des int\' egrales orbitales choisie ci-dessus co\"
incide, \`a une fonction localement constante sur $\q^{reg}$ pr\`
es,  avec celle de J. Faraut dans \cite{fa2} pour les espaces
hyperboliques et celles de Harish-Chandra pour les alg\` ebres de
Lie r\' eductives r\' eelles (\cite{va} partie I, paragraphe 3).

Par ailleurs, nous verrons que la fonction $\delta$ joue un r\^ole essentiel  dans le calcul des composantes radiales des op\' erateurs diff\' erentiels $H$-invariants \`a coefficients constants sur $\q$ et la normalisation choisie est particuli\` erement bien adapt\' ee pour l'\' etude ult\' erieure des distributions propres invariantes sur $\q$. 
\begin{lem}(\cite{or} Th\' eoreme 1.23). Soit $f\in\D(\q)$. Alors, pour tout $\a\in car(\q)$, la fonction $\mathcal{M}_Hf_{/\a^{reg}}$ est de classe $\mathcal{C}^\infty$ et \`a support born\' e dans $\a^{reg}$.
\end{lem}

\begin{Def}\label{notHinv}
Par le lemme \ref{orbites}, pour $(\la_1,\la_2)\in
\R^2\cup\{(\la,\overline{\la}); \la\in\C\}$ il existe une unique
$H$-orbite semi-simple $\mathcal{O}(\la_1,\la_2)$ de $\q$ telle que
pour tout $X\in\mathcal{O}(\la_1,\la_2)$, l'on ait $\{ u(X),
v(X)\}=\{\la_1,\la_2\}$. De plus, on a $\mathcal{O}(\la_1,\la_2)\cap\mq\neq\emptyset$ si et seulement si $(\la_1,\la_2)\in\R^2$ et $\mathcal{O}(\la_1,\la_2)\cap\a_2\neq\emptyset$ si et seulement si $\la_1=\overline{\la_2}\in\C$. 

Pour  $F$ une fonction $H$-invariante sur $\q^{reg}$, on d\' efinit
les  fonctions $F_\m$ et $F_2$ par:

$\begin{array}{lc}  \textrm{si } (t_1,t_2)\in (\R^*)^2-diag \textrm{ alors}&F_\m(t_1, t_2)= F(X) \mbox{ o\`u } X\in \mathcal{O}(t_1,t_2)\\
\textrm{et si } (\tau,\te)\in(\R^*)^2 \textrm{ alors} & F_2(\tau,\te)=F(X_{\tau,\te}) \mbox{ o\`u }
X_{\tau,\te}\in\a_2.\end{array}$

\no Pour $f\in\D(\q)$, on notera pour simplifier
$\Mfm=\big(\mathcal{M}_H (f)\big)_\m$ et $\Mf2=\big(\mathcal{M}_H
(f)\big)_2$.
\end{Def}

Avec les notations pr\' ec\' edentes, la formule d'int\' egration de
Weyl s'\' ecrit sous la forme suivante:

 \begin{lem}\label{intweyl} Soit $\Phi\in L^1_{loc}(\q)^H$ et $f\in\D(\q)$. On a la formule d'int\' egration  suivante:
 $$\int_\q \Phi(X) f(X)dX
 =\int_{t_1>t_2}(|\de| \Phi)_\m(t_1,t_2) \Mfm(t_1,t_2) dt_1 dt_2$$
 $$+8\int_{\begin{array}{c}Ê\tau>0\\ \te>0\end{array}} (|\delta| \Phi)_2(\tau,\te) \Mf2(\tau,\te) (\tau^2+\te^2)d\te d\tau.$$

 \end{lem}

\begin{dem} Avec la normalisation des mesures choisies pr\'ec\'edemment, la formule d'int\' egration de Weyl s'\' ecrit $\dis  \int_\q \Phi(X) f(X)dX=\sum_{\a\in <car(\q)>}\int_{H\cdot\a} \Phi(Z) f(Z) dZ$ avec 
$\dis  \int_{H\cdot\a} f(Z) dZ= \frac{1}{|W_H(\a)|}\int_\a \Phi (X) \mathcal{M}_H(f)(X) |\prod_{\al\in\De^+}\al(X)| dX.$

Pla\c{c}ons-nous sur  $\a=\a_{++}$.  Un \' el\' ement de $\app$ s'\'
ecrit  $X_{u_1,u_2}=u_1H_1+u_2\ka(H_1)$. Comme la base $\{
H_1,\ka(H_1)\}$ est orthonormale, on obtient donc

$$ \int_{H\cdot\app} \Phi(Z) f(Z) dZ=\frac{1}{8}\int_{\R^2} \Phi (X_{u_1, u_2}) \mathcal{M}_H(f)(X_{u_1, u_2}) |4u_1 u_2(u_1^2-u_2^2)| du_1du_2.$$

 Comme les fonctions consid\' er\' ees sont $H$-invariantes,   un   simple changement de variables donne
$$\int_{H\cdot\a} \Phi(Z) f(Z) dZ=\int_{t_1>t_2>0}|t_1-t_2| \Phi_\m(t_1,t_2) \Mfm(t_1,t_2) dt_1 dt_2.$$

Un raisonnement analogue sur $\apm$ et $\amm$ conduit \`a
$$\sum_{\a\in<car(\mq)>}\int_{H\cdot\a} \Phi(Z) f(Z) dZ=\int_{t_1>t_2}(|\de| \Phi)_\m(t_1,t_2) \Mfm(t_1,t_2) dt_1 dt_2.$$

Pour $\a_2$, on consid\` ere la base $H_3$ et
$H_4=\left(\begin{array}{c|c} 0&\begin{array}{cc} 0& -1\\1&
0\end{array} \\ \hline \begin{array}{cc} 0& -1\\1& 0\end{array}
&0\end{array}\right)$ qui est orthogonale et de volume $2$ pour
$\omega$.    Ainsi, on a
$$\int_{H\cdot\a_2} \Phi(Z) f(Z) dZ=2\int_{\R^2}  \Phi (X_{\tau, \te}) \mathcal{M}_H(f)(X_{\tau,\te}) |4\tau \te| (\tau^2+\te^2) d\tau d\te$$
$$=8 \int_{\begin{array}{c}Ê\tau>0\\ \te>0\end{array}} (|\delta|  \Phi)_2 (\tau, \te) \Mf2(\tau,\te) (\tau^2+\te^2)d\tau d\te.$$\end{dem}
La suite de cette partie est consacr\' ee \`a l'\' etude du
comportement des int\' egrales orbitales $\mathcal{M}_H(f)$ au
voisinage des points semi-r\' eguliers de $\q$.

\subsection{M\' ethode de descente}\label{descentegenerale}

Nous rappelons tout d'abord des r\'esultats de S. Sano (\cite{sa} paragraphe 2) concernant les d\'ecompositions radicielles.

On note $\tau$ la conjugaison de $\g_\C$ relative \`a la forme r\'
eelle $\g$. Soit $\a\in Car(\q)$. On a alors la d\'ecomposition
suivante de $\g_{\C}$:
$$\g_{\C}=\z_{\h_\C}(\a)\oplus\a_\C\oplus\n_\C\qq\mbox{o\`u}\quad  \n_\C=\sum_{\al\in  \De(\g_\C,\a_\C)}\g_\al^\C.$$

Il existe une base $(T_j)_{1\leq j\leq 12}$  de $\n_\C\cap\h$
v\'erifiant les propri\'et\'es suivantes: pour tout $j\in\{1,\ldots
12\}$, il existe  $\al\in  \De(\g_\C,\a_\C)^+$ telle que
$T_j=X_j+\sigma(X_j)$ o\`u  $X_j\in \g_\al^\C$ si $\al$ est r\'eelle
ou imaginaire et $X_j\in \g_\al^\C\oplus \g_{\al^\tau}^\C$ si $\al$
est complexe. On pose $\ga(T_j)=X_j-\sigma(X_j)$ si $\al$ est
r\'eelle ou complexe et   $\ga(T_j)=\im(X_j-\sigma(X_j))$ si $\al $
est imaginaire. La famille $\ga(T_j)$ est alors une base de
$\n_\C\cap\q$.

Si $X\in\a$, alors $\mathrm{ad}(X)$ induit un morphisme de
$\n_\C\cap\h$ dans $\n_\C\cap\q$. On note $det\big(\mathrm{ad}(
X)|^{\q/\a}_{\h/\z_{\h}(\a)}\big)$ le d\'eterminant de
$\mathrm{ad}(X)$ dans la base choisie pr\'ec\'edemment. Par un
calcul analogue \`a celui de S. Sano (\cite{sa}proposition 2.1), on
obtient  alors
$$ \Big|det\big(\mathrm{ad}(
X)|^{\q/\a}_{\h/\z_{\h}(\a)}\big)\Big|=
\Big|\prod_{\al\in\rap}\al(X)^{m_\al }\Big|.$$

%%%%%%%%%%%%%%%%%%%%
On consid\`ere maintenant $\l =\m$ ou $\l=\z_3$. On a alors
$car(\l\cap\q)=\{\a\in car(\q), \a\subset\l\cap\q\}$ et pour  $\a\in
car(\q\cap\l)$, on remarque que  $\z_{\l\cap \h}(\a)=\z_\h(\a)$. On
peut d\'ecomposer $\l_\C$ de mani\`ere analogue \`a $\g_\C$, ce qui
permet d'obtenir
$$\h=\l\cap\h\oplus (\r_\C\cap\h)\qq\mbox{et}\qq\q=\l\cap\q\oplus (\r_\C\cap\q),$$
o\`u $\displaystyle \r_\C=\sum_{\al\in  \De(\g_\C,\a_\C)-
\De(\l_\C,\a_\C)}\g_\al^\C.$
 La famille form\'ee des $T_j$ (respectivement $\ga(T_j)$) associ\'es \`a une racine $\al\in \De(\g_\C,\a_\C)^+- \De(\l_\C,\a_\C)^+$   d\'efinit  une base de $\r_\C\cap\h$ ( respectivement de $\r_\C\cap\q$).
  Si $X\in\a$, alors $\mathrm{ad}(X)$ induit un morphisme de $\r_\C\cap\h$ dans $\r_\C\cap\q$. On note $det\big(\mathrm{ad}(
X)|^{\q/\l\cap\q}_{\h/\l\cap\h}\big)$ son d\'eterminant dans les
bases choisies et on a
$$ \Big|det\big(\mathrm{ad}(
X)|^{\q/\l\cap\q}_{\h/\l\cap\h}\big)\Big|= \Big|\prod_{\al\in\rap-
\De(\l_\C,\a_\C)^+}\al(X)^{m_\al }\Big|.$$

\begin{lem}\label{descente}
 On pose  $\p\l\cap \q=\{X\in\l\cap\q,  \Big|det\big(\mathrm{ad}(
X)|^{\q/\l\cap\q}_{\h/\l\cap\h}\big)\Big|\neq 0\}$. L'application 
$\ga :H\times\l\cap\q \to \q$ d\' efinie par $\ga(h,X)=h\cdot X$ est submersive  en tout point de
$H\times\p\l\cap\q$. En particulier, l'int\' egration sur les fibres (lemme \ref{hc}) d\' efinit la surjection $\ga_\ast$ de $\D (H\times\p\l\cap\q)$ dans $\D(H\cdot\p\l\cap\q)$.

\end{lem}

\begin{dem}  Pour $(h,X)\in H\times \l\cap\q$, $(A,X')\in\h\times
\l\cap\q$ et $t\in\R$, on a $$\ga(h\; exp(tA), X+tX')=Ad(h)( X+t
X'+t[A,X]+t^2Y)$$
 avec $Y\in\q $. Par suite la diff\' erentielle  de $\ga$ en un point $(h,X)\in H\times \l\cap\q$ est donn\' ee   par $d_{(h,X)}\ga(A,X')=Ad(h)(X'+[A,X])$. Comme $[X,\l\cap\h]\subset \l\cap\q$ et $[X,\r_\C\cap\h]\subset \r_\C\cap\q$, la discussion pr\'ec\'edente donne le lemme.\end{dem}

\subsection{Etude de $\mathcal{M}_H(f)$ pour $f\in\D(H\cdot\p\m\cap\q)$}

On rappelle  que $N=N_H(\mq)$ et on consid\`ere  la submersion
surjective $\ga: H\times\p\mq\to H\cdot\p\mq$
  obtenue dans le paragraphe pr\' ec\' edent pour $\l=\m$. 

  Pr\' ecisons tout d'abord les fibres de $\ga$. Par le lemme \ref{descente} et la remarque \ref{rempoly}, pour tout $X\in\mq$, on a $\delta^2(X)=S_0(X)= \Big|det\big(\mathrm{ad}(
X)|^{\q/\m\cap\q}_{\h/\m\cap\h}\big)\Big|$ et donc $\p\mq=\{X\in\mq;
S_0(X)\neq 0\}$. Par le lemme  \ref{ceno},  pour
  $X\in \p\mq$, on obtient  donc $\ga^{-1}(X)=\{(h,h^{-1}\cdot X); h\in N\}.$\me

Pour
$\phi\in\D(H\times \p\mq)$ et $h\cdot X\in H\cdot \p\mq$, on pose:

$$\pi_\ast(\phi)(h.X)=\de(X)\ga_\ast(\phi)(h\cdot X)=\de(X)^{-1}\int_{N} \phi(hu^{-1},u\cdot X) du$$

L'application $\pi_\ast$ est surjective de $\D(H\times\p\mq)$ sur $\D(H\cdot \p\mq)$ puisque $\ga_\ast$ l'est.\me

Soit $p$ la projection de $H\times\p\m\cap\q$ sur
$\p\m\cap\q$. La surjection  $p_\ast: \D(H\times\p\m\cap\q)\to \D(\p\m\cap\q)$
est alors donn\' ee par $$p_\ast(\phi)(X)=\int_H\phi(h,X)dh.$$

L'int\' egrale orbitale sur $\mq$ d'une fonction
$f\in\D(\p\m\cap\q)$ est d\' efinie , pour $X\in\m\cap\qr$,  par:
$$\mathcal{M}_{N}(f)(X):=\int_{N/Z_H(X)}f(h\cdot X)d\dot{h}$$
Cette normalisation est compatible avec celle d\' efinie en d\' ebut
de paragraphe puisque, d'une part, $Z_{N}(X)=Z_H(X)\subset N$ par le
lemme \ref{ceno}  et d'autre part, pour tout $\a\in car(\mq)$, les racines de $\De(\m_\C,\a_\C)$ sont toutes  de multiplicit\' e un.
\begin{prop}\label{desc} Pour $\phi\in\D(H\times\p\m\cap\q)$ et $X\in\m\cap\qr$, on a
 $$\mathcal{M}_H(\pi_\ast(\phi))(X)=\mathcal{M}_{N}(p_\ast(\phi))(X).$$
\end{prop}

\begin{dem} Comme $H$ et $Z_H(X)$ sont unimodulaires, on
a:
$$\mathcal{M}_H(\pi_\ast(\phi))(X)=\de(X)\int_{H/Z_H(X)}\pi_\ast(\phi)(h\cdot
 X)d\dot{h}=\int_{H/Z_H(X)}\int_{N}\phi(hu^{-1},u\cdot
X)du\;d\dot{h}$$
$$=\int_{H/Z_H(X)}\int_{N/Z_H(X)}\int_{Z_H(X)}\phi(hv^{-1}u^{-1},u\cdot X)dv\;d\dot{u}\;d\dot{h}.$$
 
 $$=\int_{H}\int_{N/Z_H(X)}\phi(hu^{-1},u\cdot
 X)d\dot{u}\;dh=\mathcal{M}_{N}(p_\ast(\phi))(X).$$\end{dem}
Nous allons maintenant exprimer l'int\' egrale orbitale $\mathcal{M}_{N}(g)$  pour $g\in\D(\mq)$ en terme d'une  application moyenne $M_{Q_1,Q_2}(\overline{g})$ d\' efinie dans le paragraphe \ref{qq} o\`u $\overline{g}\in\D(\R^4)$ d\' epend de $f$.

On reprend les notations du paragraphe  \ref{semireg1},  en
particulier, on \' ecrit $\m=\z_1+\z_2=\m_1\oplus \m_2$. On pose
$Z_1=Z_H(H_1)=\{diag(a, b, a, d)\in GL(4,\R)\}$ et $Z_2=\ka(Z_1)$.

\begin{lem}\label{prod}
Si $X\in \mq^{reg}$ alors, on a les isomorphismes de groupes:
$$\begin{array}{rl}N^0/Z_H(X)&\simeq Z_1/Z_H(X)\times Z_2/Z_H(X)\\
 & \simeq \{diag(1, 1, 1, d); d\in\R^*\}\times \{diag(1, 1, c, 1); c\in\R^*\}\end{array}$$
\end{lem}

\begin{dem} Par le lemme \ref{ceno}, on a $Z_H(X)=\{diag(\al, \be, \al, \be)\in GL(4,\R)\}$. Les isomorphismes donn\' es sont alors clairs.\end{dem}

 On note $q$ la forme quadratique d\' efinie sur $\R^2$ par $q(x,y)=x^2-y^2$.  On d\' efinit l'isomorphisme   $\psi$ de $\R^2$ dans $\m_1\cap \q$  par $\psi(x,y)=\left(\begin{array}{c|c}
0&\begin{array}{cc} 0&0\\0&x+y\end{array}\\ \hline
\begin{array}{cc} 0&0\\0&x-y\end{array}&0
\end{array}
\right)$ de telle sorte que $Q\circ\psi =q$.

\begin{lem}\label{rg1}Pour  $g\in\D(\m_1\cap\q)$ et
$X_1\in(\m_1\cap\q)^{reg}$, on a
$$\int_{Z_1/Z_{ Z_1}(X_1)}g(u\cdot X_1)d\dot{u}=  M_{Q\circ \psi }(g\circ \psi) \big( Q(X_1)\big).$$
\end{lem}

\begin{dem} Pour $X_1\in(\m_1\cap\q)^{reg}$, le groupe   $Z_{ Z_1}(X_1)=\{diag(\al, \be, \al, \be)\in GL(4,\R)\}$ est ind\'ependant de $X_1$. On le note $C_1$. La formule d'int\'egration de Weyl sur $\m_1\cap\q$ assure que, pour tout $F\in\D(\R)$, on a $$\begin{array}{rl}\int_{\m_1\cap\q}  F\circ Q(X_1) g(X_1)d X_1& =\int_{\R_+} 2x F(x^2) \int_{Z_1/C_1}g(u\cdot\psi(x,0)) d\dot{u} \;dx\\
&+\int_{\R_+} 2yF(-y^2) \int_{Z_1/C_1}g(u\cdot\psi(0,y)) d\dot{u}
\;dy.\end{array}$$
Par ailleurs, on a $\int_{\m_1\cap\q}  g(X_1) F\circ Q(X_1) d X_1=\int_{\R^2} g\circ \psi(x,y) F(x^2-y^2) dx\; dy=\int_\R  M_{Q\circ \psi }(g\circ \psi )(t) F(t) dt$ ce qui donne le r\'esultat voulu. \end{dem}

De m\^eme, l'application $\ka\circ\psi$ d\' efinit un isomorphisme de $\R^2$ dans $\m_2\cap\q$ tel que $Q\circ\ka\circ\psi=q$ et on obtient un r\' esultat similaire sur $\m_2\cap\q$.

Pour $g\in\D(\mq)$ on d\' efinit l'application $\tilde{g}$ sur
$\R^2\times \R^2$  par
$\tilde{g}(u,v)=g\big( \psi(u)+\ka\circ \psi(v)\big).$

\begin{cor}\label{liend} Pour  $g\in\D(\m\cap\q)$ et
 $X=X_1+X_2\in(\m_1\oplus \m_2)\cap\q^{reg}$,  on a $$\mathcal{M}_{N^0}(g)(X)=\int_{Z_1/Z_H(X)\times Z_2/Z_H(X)}g(u_1\cdot X_1+u_2\cdot
X_2)d\dot{u_1}d\dot{u_2}=
 M_{q,q}\tilde{g}\big( Q(X_1),  Q(X_2)\big)$$ \end{cor}

\begin{dem} Ceci se d\' eduit imm\' ediatement du lemme pr\' ec\'
edent et du
lemme \ref{double}.\end{dem}

\begin{cor}\label{express} Soit $f\in\D(H\cdot\p\mq).$  Alors, il existe $f_0\in\D(\R^2\times \R^2)$ \`a support dans $\{(u,v)\in\R^2\times \R^2; q(u)\neq q(v)\}$, telle que, pour $X=X_1+X_2\in(\m_1\oplus\m_2)\cap q^{reg}$, on ait
$$\mathcal{M}_H(f)(X)=M_{q,q}f_0\big( Q(X_1),  Q(X_2)\big)+
M_{q,q}f_0\big( Q(X_2),  Q(X_1)\big)$$
\end{cor}

\begin{dem} Comme l'application $\pi_\ast$ est surjective,  il existe $\phi\in\D(H\times \p\mq)$  telle que $f=\pi_\ast(\phi)$. Par la proposition \ref{desc}, pour tout $X\in\mq^{reg}$,  on a $\mathcal{M}_H(f)(X)=\mathcal{M}_{N}(p_\ast(\phi))(X)=\mathcal{M}_{N^0}(p_\ast(\phi))(X)+\mathcal{M}_{N^0}(p_\ast(\phi))(\ka(X))$. Le corollaire pr\' ec\' edent permet de conclure en prenant $f_0= \widetilde{p_\ast(\phi)}. $ \end{dem}

Selon la d\' efinition \ref{notHinv}, pour tout $(t_1,t_2)\in(\R^*)^2-diag$, on note $\Mfm(t_1,t_2) =\mathcal{M}_H(f)(X_{t_1,t_2})$ o\`u $X_{t_1,t_2}\in\p\mq$ v\'erifie  $\{u(X_{t_1,t_2}),v(X_{t_1,t_2})\}=\{ t_1, t_2\}$. \me

On d\' efinit le sous-espace $\mathcal{H}^2_{\log}$ de $\mathcal{H}_{\eta,\eta}$  avec $\eta(t)=\log(t)$, par

 $$\mathcal{H}^2_{\log}=\{(t_1,t_2)\in(\R^*)^2-diag \mapsto\varphi_0(t_1,t_2)+\log|t_1|\varphi_1(t_1,t_2)$$
$$+\log|t_2|\varphi_1(t_2,t_1)+\log|t_1|\log|t_2|\varphi_2(t_1,t_2);\qq\varphi_0,\varphi_1,\varphi_2
\in\D(\R^2-diag) ;$$ $$ \varphi_0 \textrm{ et }  \varphi_2  \textrm{ sont sym\'
etriques par rapport aux deux variables}\}.$$

\begin{theo}\label{iom}
\begin{enumerate}
\item Pour tout $f\in\D(H\cdot\p\mq)$, on a $\Mfm \in \mathcal{H}^2_{\log}$
\item L'application
$$\left\{\begin{array}{lll}\D(H\cdot\p\mq)&\to&\mathcal{H}^2_{\log}\\
\qq\qq\qq f&\mapsto&\Mfm \end{array}\right.$$ est continue et surjective.
\end{enumerate}
\end{theo}

\begin{dem} Soit $f\in\D(H\cdot\p\mq)$. Par le corollaire \ref{express},  il existe
$f_0\in\D(\R^2\times\R^2)$, \`a support dans $\{(u,v)\in\R^2\times
\R^2; q(u)\neq q(v)\}$,  telle que pour
$X=X_1+X_2\in(\m_1\oplus\m_2)\cap\qr$ avec $t_1=Q(X_1)$ et
$t_2=Q(X_2)$, on ait
$\Mfm(t_1,t_2)=\mathcal{M}_H(f)(X)=M_{q,q}f_0(t_1,t_2)+
M_{q,q}f_0(t_2,t_1).$ Gr\^ace au th\' eor\`eme \ref{comp2} et
compte-tenu du support  de $f_0$, on obtient alors $\Mfm
\in\mathcal{H}_{\log}^2$. \me

Montrons maintenant la surjectivit\' e de l'application. Soit $F\in\mathcal{H}^2_{\log}$.  Il existe $\varepsilon>0$ tel que pour 
$|t_1-t_2|<\sqrt{\varepsilon}$ on ait $F(t_1,t_2)=0.$
D'apr\` es la proposition \ref{surj3}, il existe
$f\in\D(\R^2\times\R^2)$ tel que pour $(t_1,t_2)\in\R^*\times\R^*$,
on ait
$$F(t_1,t_2)=M_{q,q}f(t_1,t_2)=M_{q,q}(\frac{f}{2})(t_2,t_1)+M_{q,q}(\frac{f}{2})(t_1,t_2)$$
puisque les fonctions consid\' er\' ees sont sym\' etriques par
rapport aux deux variables.

Gr\^ace au corollaire \ref{liend}, la fonction  $g\in\D(\mq)$ d\'
efinie par $g(X_1+X_2)=\frac{1}{2 }f\big(\psi^{-1}(X_1), \ka\circ
\psi^{-1}(X_2)\big)$ pour $X_j\in\m_j\cap\q$, v\' erifie, pour tout
$X=X_1+X_2\in\mq^{reg}$, les relations

$\dis \mathcal{M}_{N^0} g(X)=\frac{1}{2}M_{q,q}f\big(Q(X_1), Q(X_2)\big)$ et $\dis
\mathcal{M}_{N^0} g\circ\ka(X)=\frac{1}{2}M_{q,q}f\big(Q(X_2),
Q(X_1)\big).$ Ainsi, on obtient
$\mathcal{M}_N g(X)=F(Q(X_1), Q(X_2)).$\me

Maintenant, on fixe une fonction $\chi$ dans $\mathcal{D}(\R)$ telle
que $\chi|_{[-\frac{\varepsilon}{2};\frac{\varepsilon}{2}]}=1$ et
$\chi|_{\R-[-\varepsilon;\varepsilon]}=0$. La fonction  $(1-\chi\circ S_0)g$ appartient cette fois-ci \`a
$\D(\p\mq)$ et v\' erifie
$\mathcal{M}_N\big((1-\chi\circ S_0)g\big)(X)=(1-\chi((t_1-t_2)^2))F(t_1,t_2)=F(t_1,t_2),$
pour $X=X_1+X_2$ tel que $Q(X_1)=t_1$ et $Q(X_2)=t_2$. 

Par la
surjectivit\' e de  $p_\ast$, il existe $\phi\in \D(H\times\p\mq)$
telle que $(1-\chi\circ S_0)g=p_\ast(\phi)$ et par la proposition
\ref{desc} on a
$\mathcal{M}_H(\pi_\ast(\phi))(X)=\mathcal{M}_{N}(p_\ast(\phi))(X)=F(t_1,t_2)$
ce qui donne la surjectivit\' e. La continuit\' e est imm\' ediate par le th\' eor\`eme du graphe ferm\' e.\end{dem}

\subsection{Etude de $\mathcal{M}_H(f)$ pour $f\in\D(H\cdot\p \z_3\cap\q)$}

 On rappelle (lemme \ref{ceno3}) que
$N_3=N_H(\z_3\cap\q)=N_3^0\cap\tilde{K}N_3^0$ avec
$\tilde{K}=\left(\begin{array}{c|c} I_2& 0\\ \hline 0 &
-I_2\end{array}\right)$ et $N_3^0=\left\{ \left(\begin{array}{c|c}
A&0\\ \hline
0&A\end{array}\right ); A\in GL(2,\R)\right\}$.
On note $\ga_3:  H\times \p\z_3\cap \q\rightarrow H.\p\z_3\cap \q$ la
submersion surjective obtenue dans le lemme \ref{descente} pour
$\l=\z_3$.

\begin{lem}\label{fibt}
\begin{enumerate}
\item On garde les notations du paragraphe \ref{descentegenerale}. Pour $X\in\z_3\cap\q$, on pose $\nu_3(X)=| det\big( \mathrm{ad} (X)|^{\q/\z_3\cap\q}_{\h/\z_3\cap\h}\big)|$.

Si    $X= \left (\begin{array}{c|c} 0&Y\\
\hline Y&0\end{array}\right) \in  \z_3\cap\q$ alors on a
$\nu_3(X)=4 (tr Y)^2\vert det(Y)\vert.$

En particulier, on a $\p\zt=\{X= \left (\begin{array}{c|c} 0&Y\\
\hline Y&0\end{array}\right); tr Y det Y\ne 0\}$.
\item Soit $\ztp=\left\{\left (\begin{array}{c|c} 0&Y\\ \hline
Y&0\end{array}\right);\qq (tr Y)^2\det(Y)>0\right\}.$ Alors
$$H\cdot\p\zt=H\cdot\ztp.$$

\item Si $h\cdot\ztp\cap\ztp\ne\emptyset$ alors $  h\in N_3.$

\end{enumerate}
\end{lem}

\begin{dem} Soit $\a\in <car(\zt)>$ et $X\in \a$. Comme $ \De((\z_3)_\C,\a_\C)^+=\{\be_1\}$,  le lemme \ref{descente} donne  $\nu_3(X)= \Big|(\al_1\al_2\be_2^2)(X)\Big|.$ Or, par la remarque \ref{rempoly}, on a $(\al_1\al_2)(X)=4\sqrt{S(X)}= 4\vert det(Y)\vert$ et $\be_2^2(X)=(tr Y)^2$,  et donc $\nu_3(X)=4( tr Y)^2 \vert det(Y)\vert$. L'assertion {\it 1.} s'en d\' eduit par densit\' e des \' el\' ements semi-simples.\me

Soit $ X= \left (\begin{array}{c|c} 0&Y\\
\hline Y&0\end{array}\right) \in\p\zt$. Si $det (Y)>0$ alors $X\in\ztp$ sinon la matrice $Y$ poss\` ede deux valeurs propres non nulles de signe contraire. En particulier, il existe $P\in GL(2,\R)$ telle que $PYP^{-1}= \left (\begin{array}{cc} u_1& 0\\
 0 & u_2\end{array}\right) $ avec $u_1u_2<0$. On a alors 
$$\varrho\left(\left (\begin{array}{c|c} P &0\\\hline 0&P\end{array}\right
)\cdot\left (\begin{array}{c|c} 0&Y\\ \hline
Y&0\end{array}\right)\right)= \left(
\begin{array}{c|c}
0&\begin{array}{cc} u_{1}&0\\0&-u_{2}\end{array}\\ \hline
\begin{array}{cc} u_{1}&0\\0&-u_{2}\end{array}&0
\end{array}
\right)$$ avec $\dis \varrho=\Ad \left(
\begin{array}{c|c}
\begin{array}{cc} 1&0\\0&-1\end{array}&0\\ \hline
0&\begin{array}{cc} 1&0\\0&1\end{array}
\end{array}
\right)$   ce qui donne l'assertion {\it 2.}\me

Soit $X\in\ztp$ et $h\in H$ tels que $h\cdot X\in\ztp$.

Si  $X$ est r\' egulier dans $\q$ alors, par le lemme \ref{ceno3},
il existe $h_0\in N_3$ tel que $\z_\h( h_0\cdot X)=\a$ avec
$\a=\app$ ou $\a= \a_2$. Comme $h\cdot X$ est \' egalement r\'
egulier et  l'action adjointe d'un \' el\' ement de $H$ ne modifie
pas  la nature r\' eelle ou complexe des racines, il existe $h_1\in
N_3$ tel que $\z_\h(h_1h\cdot X)=\a$. Ainsi, on a $h_1 h h_0^{-1}\in
N_H(\a)$. Comme $X\in\ztp$, le lemme \ref{WH} assure que $h_1 h
h_0^{-1}\in N_3$ et donc $h\in N_3$.

Si $X$ n'est pas r\' egulier, on \' ecrit  $ X= \left (\begin{array}{c|c} 0&Y\\
\hline Y&0\end{array}\right) $ et les deux valeurs propres de $Y$
sont \' egales.  Ainsi, la d\' ecomposition de Jordan de $Y$ s'\'
ecrit sous le forme
$Y=\left(\begin{array}{cc}
u&0\\0&u\end{array}\right)+Y_n,$ o\`u $Y_n$ est une matrice
nilpotente et $u$ est un r\' eel non nul.

En notant $h=\left (\begin{array}{c|c} A & 0\\
\hline 0 & B\end{array}\right)$, on a $h\cdot X=u\left (\begin{array}{c|c} 0&AB^{-1}\\
\hline BA^{-1}&0\end{array}\right)+h\cdot\left (\begin{array}{c|c} 0&Y_n\\
\hline Y_n&0\end{array}\right)$ et ceci est la d\' ecomposition de Jordan de
$h\cdot X$  dans $\z_3$ (\cite{va} Partie I,  lemme 3 du chapitre 1). On obtient donc  $u\left (\begin{array}{c|c} 0&AB^{-1}\\
\hline BA^{-1}&0\end{array}\right)\in\z_3$.
Ainsi, on a $AB^{-1}=BA^{-1}$
et par suite,  il existe $P$ dans $GL(2,\R)$ et
$\varepsilon_1,\varepsilon_2$ dans $\{\pm1\}$ tels que
$AB^{-1}=P\left(\begin{array}{cc}
\varepsilon_{1}&0\\0&\varepsilon_{2}\end{array}\right)P^{-1}.$
Ainsi
$$h\cdot X=\left (\begin{array}{c|c} 0&AYA^{-1}P\left(\begin{array}{cc}
\varepsilon_{1}&0\\0&\varepsilon_{2}\end{array}\right)P^{-1}\\
\hline P\left(\begin{array}{cc}
\varepsilon_{1}&0\\0&\varepsilon_{2}\end{array}\right)P^{-1}AYA^{-1}&0\end{array}\right).$$
Comme $h\cdot X$ est dans $\ztp$ alors
$
0<\det(AYA^{-1}P\left(\begin{array}{cc}
\varepsilon_{1}&0\\0&\varepsilon_{2}\end{array}\right)P^{-1})=\varepsilon_{1}\varepsilon_{2}\det(Y).$
Comme $X\in \ztp$, on a $\det(Y)>0$ et donc
$\varepsilon_{1}=\varepsilon_{2}$. On obtient alors  $h\in N_3$,
d'o\`u l'assertion {\it 3.}\end{dem}

Nous souhaitons pr\' eciser le comportement des int\' egrales
orbitales au voisinage des points $H$-conjugu\' es \`a
$\app\cap\a_2\subset\ztp$. Pour cette \' etude, nous consid\' erons
la submersion surjective
$$\pi_3 : H\times\ztp\to H\cdot\ztp.$$

Par le lemme \ref{fibt}, pour $h\cdot X\in  H\cdot\ztp$, on a
$\pi_3^{-1}(h\cdot X)=\{(hu^{-1},u\cdot X); u\in N_3\}$.
Ainsi, l'application surjective  $(\pi_3)_\ast: \D (H\times \ztp)\rightarrow \D(H\cdot\ztp)$  est donn\' ee par $$(\pi_3)_\ast (\varphi)(h\cdot  X)=\nu_3(X)^{-1}\int_{N_3}\varphi(h u^{-1}, u\cdot X) du$$ pour $\varphi\in \D (H\times \ztp)$ et $(h,X)\in H\times \ztp.$
\begin{lem}\label{centre}
Pour $\tau\in\R$ et $X'\in[\z_3,\z_3]\cap\q,$ tel que $\tau H_3+X'$
soit dans $\zt^{reg}$, on a:
$$Z_H(\tau H_3+X')=Z_{N_3^0}(X').$$
\end{lem}

\begin{dem} On note $X=\tau H_3+X'=\left(\begin{array}{c|c} 0& Y\\
\hline Y & 0\end{array}\right)$. Montrons dans un premier temps que
$Z_H(X)$ est inclus dans $N^0_3$.

Soit $h\in Z_H(X)$. Comme $X$ est r\' egulier, par le lemme
\ref{ceno3},  il existe  $h_1\in N_3$ tel que $X_1=h_1\cdot X$ soit
un \' el\' ement de  $\app$ ou de $\a_2$. Si $X_1\in \a_2$ alors
$X_1$ et donc $X$ appartiennent \`a $\ztp$ et le lemme \ref{fibt}
donne $h\in N_3$.  Si $X_1\in\app$ alors $h_1hh_1^{-1}\in
Z_H(X_1)\subset N_3$ par le lemme \ref{ceno} et donc, on a \'
egalement $h\in N_3$.

Maintenant, par r\' egularit\' e de $X$, on a  \' egalement $tr(Y)\ne 0$ et donc $Y$ et $-Y$ ne sont pas conjugu\' es ce qui implique que $h\in N_3^0$. 

Or, pour  $h$ dans $N^0_3$, on a
$h\cdot(\tau H_3+X')=\tau H_3+h\cdot X'$ ce qui donne le lemme.\end{dem}

On introduit la fonction $\de_3$ d\' efinie sur $\zt$ par:
$$\de_3(\tau H_3+X')=\im^{Y(-Q(X'))}\sqrt{2\vert Q(X')\vert}.$$
o\`u $\tau\in\R$ et $X'\in[\z_3,\z_3]\cap\q$.  Pour $\a\in car(\zt)$ et $X\in \a$, on a

$$|\de_3(X)|=|\prod_{\al\in\De^+((\z_3)_\C,\a_\C)}\al^{\mul-1}(X)|.$$

L'int\' egrale orbitale  d'une fonction
$f\in\D(\ztp),$ est  alors donn\' ee pour $X\in \ztp^{reg}$,  par:
$$\mathcal{M}_{N_3}(f)(X)=| \de_3(X)|\int_{N_3/Z_{N_3}(X)}f(h\cdot X)d\dot{h}.$$

 On note $p_3$ la projection de
$H\times\ztp$ sur $\ztp$. Compte-tenu de la normalisation des int\' egrales orbitales choisie,
on introduit l'application $(\widetilde{p_3})_\ast$ d\' efinie sur
$\D(H\times\ztp)$ par
$$(\widetilde{p_3})_\ast(
\phi)(X)=\frac{|\de(X)|}{|\de_3(X)|\nu_3(X)}(p_3)_\ast(\phi)(X),$$
 o\`u $(p_3)_\ast(\phi)(X)=\int_H\phi(h,X)dh$
pour $\phi\in\D(H\times\ztp)$ et
$X\in\ztp$.  L'application $(\widetilde{p_3})_\ast$  est surjective de $\D(H\times\ztp)$ dans $\D(\ztp)$.

\begin{prop}\label{desct} Pour $\phi\in\D(H\times\ztp)$ et $X\in\ztp^{reg}$, on a

$$\mathcal{M}_H((\pi_3)_\ast(\phi))(X)=\mathcal{M}_{N_3}((\widetilde{p_3})_\ast(\phi))(X).$$
\end{prop}

\begin{dem} Soit $X\in \ztp^{reg}$.  Par le lemme \ref{centre}, on a
$Z_H(X)=Z_{N_3}(X)$. La preuve est ensuite identique \`a celle de la proposition  \ref{desc}  \end{dem}

\begin{rem}\label{sansc} Soit $F\in\D(\zt)$. Comme $N_3=N_3^0\cup\tilde{K}N_3^0$ et $Ad(\tilde{K})|\zt=-Id|\zt$, pour $X\in\ztp^{reg}$, on a

$$\int_{N_3/Z_H(X)}F(l\cdot X)d\dot{l}=\int_{N^0_3/Z_H(X)}F(l\cdot X)d\dot{l}+\int_{N^0_3/Z_H(X)}F(-l\cdot X)d\dot{l}$$
Par ailleurs, gr\^ace au lemme \ref{centre},  en posant $X=\tau
H_3+X'$ avec $\tau\in\R$ et $X'\in[\z_3,\z_3]\cap\q,$  on a
$$\int_{N^0_3/Z_H(X)}F(l\cdot X)d\dot{l}=\int_{N_3^0/Z_{N_3^0}(X')}F(\tau H_3+l\cdot X')d\dot{l}.$$
\end{rem}\me

Nous exprimons maintenant  l'int\' egrale orbitale $\mathcal{M}_{N_3}(g)(\tau H_3+X')$  de  $g\in\D(\zt)$,   en terme d'une  application moyenne $M_{Q}(\tilde{g}_\tau)$ d\' efinie dans le paragraphe \ref{q1} o\`u $\tilde{g}\in\D(\R^4)$ d\' epend de $g$.

L'espace $[\z_3,\z_3]\cap\q$ est isomorphe \`a $\mathfrak{sl}(2,\R)$
par $i_1\left(\left(\begin{array}{c|c}0& Y\\ \hline Y& 0\end{array}
\right)\right)=Y$. Le groupe $ N^0_3$ est isomorphe \`a $GL(2,\R)$
et pour $A\in GL(2,\R)$, $Y\in \mathfrak{sl}(2,\R)$, on a $A \cdot
Y=i_1(\left(\begin{array}{c|c}A& 0\\ \hline 0& A\end{array}
\right)\cdot \left(\begin{array}{c|c}0& Y\\ \hline Y& 0\end{array}
\right))$. De plus,  la restriction de $Q$ \`a
 $[\z_3,\z_3]\cap\q$ est  la forme quadratique $N^0_3$-invariante donn\' ee par  $Q\left(\left(\begin{array}{c|c}0& Y\\ \hline Y& 0\end{array} \right)\right)=-2 det(Y)$ pour $Y\in\mathfrak{sl}(2,\R)$. \me

Notons $q_3$  la forme quadratique d\' efinie sur $\R^3$ par
$q_3(x,y,z)=x^2+y^2-z^2$ et $\psi_3$ l'isomorphisme de $\R^3$ dans
$[\z_3,\z_3]\cap \q$ d\' efini par
$$\psi_3(x,y,z)=\left(\begin{array}{c|c}
0&\begin{array}{cc} x&y+z\\y-z & -x \end{array}\\ \hline
\begin{array}{cc} x&y+z\\y-z & -x \end{array}&0
\end{array}
\right)$$
 de telle sorte que $Q\circ\psi_3=2 q_3$.

\begin{lem} Pour tout  $g\in\mathcal{D}([\z_3,\z_3]\cap\q)$ et pour tout
$X\in([\z_3,\z_3]\cap\q)^{reg}$, on a
$$|\de_3(X)| \int_{N^0_3/Z_{N^0_3}(X)}g(l\cdot
X)d\dot{l}= M_{q_3 }\big(g\circ \psi_3\big)(\frac{Q}{2}(X))$$
\end{lem}

\begin{dem} Ce r\' esultat s'obtient comme le lemme \ref{rg1} en \' ecrivant la formule d'int\' egration de Weyl pour la paire sym\' etrique $\big(\mathfrak{sl}(2,\R)\times \mathfrak{sl}(2,\R),\mathfrak{sl}(2,\R)\big)$.  \end{dem}

Pour $g\in\D(\zt)$ on d\' efinit la fonction $\tilde{g}$ sur
$\R\times \R^3$ par

$$\tilde{g}(\tau, u)=g\big(\tau H_3+\psi_3(u)\big)$$

\begin{cor}\label{liendt} Soit  $g\in\D(\zt)$ et $X=\tau H_3+X'\in(\z_3\cap\q)^{reg}$. On note $\displaystyle w=\frac{Q(X')}{2}$. On a alors $$\mathcal{M}_{N_3}(g)(X)= M_{q_3}\big(\tilde{g}_\tau+\tilde{g}_{-\tau}\big)(w).$$
\end{cor}

\begin{dem} Par la remarque \ref{sansc}, pour $X=\tau H_3+X'\in(\R
H_3\oplus [\z_3,\z_3])\cap \q^{reg}$, on a

$$\begin{array}{rl} \mathcal{M}_{N_3}g(X)&=|\de_3(X)|\int_{N^0_3/Z_H(X)}g(\tau H_3+l\cdot X')d\dot{l}\\
 &+|\de_3(-X)|\int_{N^0_3/Z_H(X)}g(-\tau H_3-l\cdot X')d\dot{l}.\end{array}$$
 
Comme   $Q(X')=Q(-X')$, le r\' esultat d\' ecoule du lemme pr\' ec\' edent.  \end{dem}

\begin{cor}\label{express3} Soit $f\in\D(H\cdot\ztp)$. Alors il existe une fonction $f_1\in\D(\R^4)$ paire par rapport \`a la premi\` ere variable et \`a support contenu dans $\{(\tau,u)\in\R^*\times\R^3; 2\tau^2-q_3(u)>0 \}$, telle que pour tout $(\tau, X')\in\R \times [\z_3,\z_3]\cap \q$ v\'erifiant $\tau H_3+X'\in\ztp^{reg}$, on ait
$$\mathcal{M}_H(f)(\tau H_3+X')=M_{q_3}(f_1)_\tau\big( Q(X')/2\big).$$
\end{cor}

\begin{dem} Comme l'application $(\pi_3)_\ast$ est surjective, il
existe $\phi\in\D(H\times\ztp)$ telle que $f=(\pi_3)_\ast(\phi)$. On
note $g=(\widetilde{p_3})_\ast(\phi)$ de sorte que
$\mathcal{M}_H(f)(X)=\mathcal{M}_{N_3}(g)(X)$.

Soit $f_1$ la fonction d\' efinie sur $\R\times \R^3$ par $f_1(\tau, u)=\tilde{g}(\tau, u)+\tilde{g}(-\tau, u)$. Comme $g\in\D(\ztp)$, la fonction $\tilde{g}$ est \`a support contenu dans $\{(\tau,u)\in\R^*\times\R^3; 2\tau^2-q_3(u)>0 \}$. Par le corollaire pr\' ec\' edent, on a $\mathcal{M}_H(f)(\tau H_3+X')= M_{q_3}(f_1)_\tau\big( Q(X')/2\big)$ et $f_1$ v\' erifie les propri\' et\' es demand\' ees. \end{dem}

Soit $(\tau,w)\in\R^*\times \R^*$ tel que $\tau^2>w$. Par le lemme
\ref{orbites} et la remarque \ref{rempoly}, il existe une unique
$H$-orbite r\' eguli\` ere $\mathcal{O}=H\cdot X(\tau,w)$ avec
$X(\tau,w)=\tau H_3+X'\in \q^{reg}\cap(\R H_3\oplus [\z_3,\z_3])_+$
tel que $Q(X(\tau,w))=2(\tau^2+w)$ et
$S(X(\tau,w))=\big(\tau^2-w\big)^2$. Ainsi, pour
$f\in\D(H\cdot\ztp)$, on d\' efinit la fonction $G_f$ sur $(\R^*)^2
$ par
$$G_f(\tau,w)=\left\{\begin{array}{ll} \mathcal{M}_H(f)(X(\tau,w)) & \mbox{ si } \tau^2>w\\
0 & \mbox{ sinon}\end{array} \right.$$
Par le corollaire \ref{express3}, la fonction $G_f$ est paire par rapport \`a $\tau$.

\begin{rem} Pour $ w>\tau^2>0$, on a imm\' ediatement $H\cdot X(\tau,w)=H\cdot X(\sqrt{w},\tau^2)$ et donc, dans ce cas, on obtient
$$ \mathcal{M}_H(f)(X(\tau,w)) =G_f(\sqrt{w},\tau^2).$$
\end{rem}

\no On d\' efinit le sous-espace $ \mathcal{H}_Y^{pair}$ de $\mathcal{H}^1_\eta$ o\`u $\eta(t)=Y(t)$ est la fonction d'Heaviside, par:
 $$\mathcal{H}_Y^{pair}=\{(\tau,w)\in\R^*\times \R\mapsto
a(\tau,w)+Y(-w)|w|^{1/2}b(\tau,w);$$
$$a,b\in\D(\R^*\times\R\cap\{(\tau,w); \tau^2>w\})\quad \textrm{paires par rapport \`a }\tau\}.$$
%%%%%%%%%%%
\begin{theo} \label{ioa2}

\begin{enumerate}

\item Pour tout $f\in\D(H\cdot\ztp)$, on a $G_f\in\mathcal{H}_Y^{pair}$
\item L'application
$$\left\{\begin{array}{lll}\D(H\cdot\ztp)&\to&\mathcal{H}_Y^{pair}\\
\qq\qq\qq f&\mapsto&G_f\end{array}\right.$$ est surjective et continue.

\end{enumerate}
\end{theo}

%%%%%%%%%%%
\begin{dem} Soit $f\in\D(H\cdot\ztp).$ Montrons tout d'abord que $G_f$ est
dans $\mathcal{H}_Y^{pair}$.

Gr\^ace au corollaire \ref{express3}, il existe $f_1\in\D(\R^4)$
paire par rapport \`a la premi\` ere variable et \`a support contenu
dans $\{(\tau,u)\in\R^*\times\R^3; 2\tau^2-q_3(u)>0 \}$, telle que
 pour $X=\tau H_3+X'\in\ztp^{reg}$ avec $\tau\in\R$ et $X'\in[\z_3,\z_3)\cap\q$ v\' erifiant
$Q(X')/2=w$, on ait
$$G_f(\tau,w)= \mathcal{M}_H(f)(X)=M_{q_3}(f_1)_\tau(w).$$
Gr\^ace au th\' eor\`eme \ref{surj2} et compte-tenu du  support de
$f_1$, on obtient alors $G_f\in\mathcal{H}_Y^{pair}.$
\medskip

Montrons maintenant la surjectivit\' e de l'application. 
Soit $G\in\mathcal{H}_Y^{pair}$. Il existe $\varepsilon>0$ tel
que
$$\tau^2<\varepsilon\;\textrm{ou}\;\tau^2-w<\varepsilon\Longrightarrow
G(\tau,w)=0.$$
Gr\^ace au th\' eor\`eme \ref{surj2}, il existe
$g\in\D(\R\times\R^3)$ tel que pour $\tau,w\in\R^*$, on ait
$$G(\tau,w)=M_{q_3}(g_\tau)(w)=M_{q_3}(\dfrac{1}{2}g_\tau)(w)+M_{q_3}(\dfrac{1}{2}g_{-\tau})(w).$$
Gr\^ace au corollaire \ref{liendt}, la fonction
$g_0\in\D(\z_3\cap\q)$ d\' efinie par $g_0(\tau
H_3+X')=g(\tau,\psi^{-1}_3(X'))+g(-\tau,\psi^{-1}_3(X'))$ v\'
erifie
$$\mathcal{M}_{N_3}(g_0)(\tau H_3+X')=G(\tau,Q(X')/2).$$
Soit $\chi$ une fonction de classe $\mathcal{C}^\infty$ nulle  sur $]-\infty;\frac{\varepsilon}{2}]$ et \' egale \`a $1$ sur
$[\varepsilon;+\infty[.$ Pour $X=\tau H_3+X'\in(\R
H_3+[\z_3,\z_3])\cap\q$, on pose
$\phi_1(X)=\tau^2$ et $\phi_2(X)=\tau^2-Q(X')/2.$
Soit $\te$ la fonction d\' efinie pour $X=\tau H_3+X'\in(\R
H_3+[\z_3,\z_3])\cap\q$ par
$$\te(X)=\phi_1(X)\phi_2(X)=\chi(\tau^2)\chi(\tau^2-Q(X')/2).$$
Cette fonction est $N_3$-invariante et pour
$\tau^2\leq\frac{\varepsilon}{2}$ ou
$\tau^2-Q(X')/2\leq\frac{\varepsilon}{2}$, on a $\te(X)=0$ et pour
$\tau^2\geq\varepsilon$ et $\tau^2-Q(X')/2\geq\varepsilon$, on a
$\te(X)=1$. Ainsi $\te g_0$ est \`a support dans $\ztp$ et v\'
erifie pour $Q(X')/2=w$ la relation
$$\mathcal{M}_{N_3}(\te
g_0)(X)=\chi(\tau^2)\chi(\tau^2-w)G(\tau,w)=G(\tau,w).$$ Par
surjectivit\' e de l'application $(\widetilde{p_3})_*$, il existe
$\phi\in\D(H\times\ztp)$ telle que $\te
g_0=(\widetilde{p_3})_*(\phi)$ et par la proposition \ref{desct}, on
a
$\mathcal{M}_H((\pi_3)_*(\phi))(X)=\mathcal{M}_{N_3}((\widetilde{p_3})_*(\phi))(X)=G(\tau,w)$,
ce qui ach\`eve la preuve.\end{dem}

\section{Op\' erateurs diff\' erentiels et solutions propres}

La forme $H$-invariante et non d\' eg\' en\' er\' ee $\displaystyle \omega(X,X')=\frac{1}{2}tr(XX')$ induit un isomorphisme de $S(\q_\C)^{H_\C}$ dans $\cqh=\C[Q,S]$. Dans cette partie, nous d\' eterminons les fonctions $H$-invariantes et analytiques sur $\qr$ solutions du syst\`eme 
$$(E_\chi)(\qr)\left\{\begin{array}{c}\partial(Q)\Phi=\chi(Q)\Phi\\ \partial(S)\Phi=\chi(S)\Phi\end{array}\right.,$$
o\`u $\chi$ est un caract\`ere de $\C[Q,S]=\cqh.$ Comme $\cqh$ est isomorphe \`a $\C[\app]^{W_H(\app)}$, un tel
caract\`ere est d\' etermin\' e par  $(\la_1,\la_2)\in\C^2$ avec
$\chi(Q)=\la_1+\la_2$ et $\chi(S)=\la_1\la_2$ et il
 est dit r\' egulier si $\la_1\la_2(\la_1-\la_2)\neq 0.$\me

Dans toute la suite, on fixe un caract\`ere r\' egulier $\chi$
associ\' e \`a $(\la_1,\la_2)\in\C^2$.\me

Soit $\a$ un sous-espace de Cartan. Nous rappelons que la partie
radiale sur $\a$ d'un op\' erateur diff\' erentiel $D$ d\' efini sur
$\qr$, que l'on notera $\Delta_{\mathfrak{a}}( D)$, est l'op\'
erateur diff\' erentiel d\' efini sur $\a^{reg}$ v\' erifiant pour
toute fonction $f$ de $\mathcal{C}^\infty(\q)^H,$ la relation
$$\Delta_{\mathfrak{a}}( D)(f|_{\are})=(Df)|_{\are}.$$

Ainsi, la fonction $\Phi$ est solution de $(E_\chi)(\qr)$ si et seulement si elle satisfait sur $\are$ le syst\`eme
$$(E_\chi)(\are)\left\{\begin{array}{c}\rdq \Phi=(\la_1+\la_2)\Phi\\\rds \Phi=\la_1\la_2\Phi\end{array}\right..$$

Dans cette partie, nous allons  d\' eterminer, pour tout \sec  $\a$
de $\q$, les parties radiales de $\dq$ et $\ds$. Nous d\'
eterminerons ensuite les solutions sur $\are$ du syst\`eme
$(E_\chi)(\are)$.

\subsection{Parties radiales}
Soit $\a$ un sous-espace de Cartan. Pour $\xi\in\a$ et $k$ une
fonction $\w$-invariante de $\ra$ dans $\C$, on note $T_\xi(k)$
l'op\' erateur de Dunkl d\' efini par
$$T_\xi(k)=\partial_\xi+\sum_{\al\in\rap}k_\al\al(\xi)\frac{1-r_\al}{\al},$$
avec $r_\al=id_\a-\frac{2\al}{\al(h_\al)}h_\al.$ L'op\' erateur $T_\xi(k)$ agit  sur $\C[\a]$ et sur $\mathcal{C}^\infty(\a)$. \me

Nous rappelons
tout d'abord  quelques propri\' et\' es remarquables  de ces op\' erateurs de Dunkl dont les preuves sont pr\' ecis\' ees dans \cite{du} et \cite{he}.

 Pour $\xi,\xi'\in\a,$ on  a $[T_\xi(k),T_{\xi'}(k)]=0$.  Ainsi l'application $\xi\mapsto T_\xi(k)$ de $\a$
dans l'ensemble des op\' erateurs de $\mathcal{C}^\infty(\a^{reg})$
s'\' etend en un homomorphisme injectif de $S(\a_\C)$ dans
l'ensemble des op\' erateurs  de $\mathcal{C}^\infty(\a^{reg})$. On
note $T_p(k)$ l'image de $p\in S(\a_\C)$ par cette application.

Lorsque  $p\in S(\a_\C)^{\w}$, l'op\' erateur $T_p(k)$ est un op\' erateur diff\' erentiel sur l'espace des fonctions $\w$-invariantes sur $\a$. On note  $D_p(k)$ la restriction de $T_p(k)$ \`a  l'espace des fonctions $\w$-invariantes sur $\a$. \me

Soit   $Res_\a$ l'homomorphisme de restriction de $\C[\q_\C]^{H_{\C}}$ dans $\C[\a_\C]^{\w}$ donn\' e par $
P\mapsto P|_{\a_\C}.$
La restriction de $\om$ \`a $\a_\C\times\a_\C$ est  non d\' eg\' en\' er\' ee, donc  elle induit un isomorphisme de $S(\a_\C)^{\w}$ dans $\C[\a_\C]^{\w}$ et on notera  \' egalement   $Res_\a$ l'homorphisme de restriction  de $S(\q_\C)^{H_\C}$ dans $S(\a_\C)^{\w}$ que l'on d\' eduit.

\begin{theo}\label{rad} (\cite{to2} paragraphe 1.1) Soit $P\in S(\q)^H$ et $k$ la fonction d\' efinie sur $\ra$ par $k_\al=\frac{\mul}{2}.$ Alors l'action de
$\rdp$ co\"incide avec celle de $ D_{Res_\a(P)}(k) $ sur l'espace
des fonctions $ \w $-invariantes sur $\a.$

\end{theo}

Dans toute la suite, nous consid\' erons la fonction  $k$  d\'
efinie par $k_\al=\frac{\mul}{2},$ pour $\al$ parcourant $\ra.$ On
pose
$$I(k)=\prod_{\al\in\rap}\al^{2k_\al}=\prod_{\al\in\rap}\al^{m_\al}.$$
\begin{rem}\label{Id} On a $|I(k-\frac{1}{2})|=|\de| \mbox{ et } |I(k)|=\Pi$
\end{rem}

\begin{theo}(\cite{op} proposition 3.9(1)).  Pour $p\in S(a_\C)^{\w}$, on a
$$D_p(1-k)=I(k-\frac{1}{2})\circ D_p(k)\circ I(\frac{1}{2}-k)=\de\circ D_p(k)\circ
\de^{-1}.$$
\end{theo}
\begin{dem} Le r\' esultat d'Opdam donne la premi\`ere \' egalit\' e,
la deuxi\`eme en d\' ecoule car  les fonctions $\frac{\de}{|\de|}$
et $\frac{I(k-\frac{1}{2})}{|I(k-\frac{1}{2})|}$ sont localement
constantes sur $\are. $ \end{dem}

Pour $\al\in\ra,$ on pose $$\el:=\frac{1}{4\al}\partial
h_\al(\al\partial h_\al)=\frac{1}{4}\partial
h_\al^2+\frac{1}{\al}\partial h_\al.$$ Nous rappelons que l'on note
$\al_1$ et $\al_2$ les deux racines positives de multiplicit\' e un
de $\a$.

\begin{prop}\label{prad} On a
$$\rdq=\de^{-1}\circ(\elu+\eld)\circ\de=I(1/2-k)\circ(\elu+\eld)\circ I(k-1/2)$$
$$\rds=\de^{-1}\circ(\elu\eld)\circ\de=I(1/2-k)\circ(\elu\eld)\circ I(k-1/2).$$
\end{prop}
\begin{dem} La premi\`ere assertion est d\^ue \`a Dunkl (\cite{du} et \cite{he} th\' eor\`eme 1.6) et s'obtient par un simple calcul que nous reprenons ci-dessous pour obtenir la deuxi\`eme assertion.\me

Pour $P\in S(\q)^H$, on note $p=Res_\a (P)$. Gr\^ace aux th\'
eor\`emes pr\' ec\' edents, sur l'espace des fonctions
$\w$-invariantes sur $\a_\C$, on a les \' egalit\' es suivantes:
 $$\rdp=\de^{-1}\circ D_p(1-k)\circ\de=I(\frac{1}{2}-k)\circ
D_p(1-k)\circ I(k-\frac{1}{2}).$$

Calculons $D_p(1-k)$ pour $P=Q$ et $P=S$. En utilisant l'isomorphisme  $\a_\C\simeq\a_\C^*$ induit par $\omega$, on a $Res_\a(Q)=\frac{\al_1^2+\al_2^2}{4}=\frac{h_{\al_1}^2+h_{\al_2}^2}{4}$ et $Res_\a(S)=\frac{\al_1^2\al_2^2}{16}=\frac{h_{\al_1}^2h_{\al_2}^2}{16}.$

Par d\' efinition, on a $T_\xi(1-k)=\partial_\xi+\frac{1}{2}\al_1(\xi)\frac{1-r_{\al_1}}{\al_1}
+\frac{1}{2}\al_2(\xi)\frac{1-r_{\al_2}}{\al_2}.$ Comme les racines $\al_1$ et $\al_2$ sont  orthogonales et
$\al_i(h_{\al_i})=Q(h_{\al_i})=4$, pour $i\in\{1,2\}$, on obtient 
$\dis T_{h_{\al_i}}(1-k)=\partial
h_{\al_i}+2\frac{1-r_{\al_i}}{\al_i},$ et par suite
$$T_{h_{\al_i}^2}(1-k)=\partial
h_{\al_i}^2+2\partial
h_{\al_i}\circ\frac{1-r_{\al_i}}{\al_i}+2\frac{1-r_{\al_i}}{\al_i}\circ\partial
h_{\al_i}+4\left(\frac{1-r_{\al_i}}{\al_i}\right)^2.$$

Les relations  $r_{\al_i}\circ\partial h_{\al_i}=\partial(r_{\al_i}(
h_{\al_i}))\circ r_{\al_i}=-\partial h_{\al_i}\circ r_{\al_i}$ et 
$\dis \left(\frac{1-r_{\al_i}}{\al_i}\right)^2=0,$  donnent finalement
$$T_{h_{\al_i}^2}(1-k)=\partial
h_{\al_i}^2+\frac{2}{\al_i}\partial
h_{\al_i}+\frac{2}{\al_i}\partial h_{\al_i}\circ r_{\al_i}+2\partial
h_{\al_i}\circ\frac{1-r_{\al_i}}{\al_i}.$$

On obtient donc  les \' egalit\' es suivantes:
$$D_{h_{\al_1}^2+h_{\al_2}^2}(1-k)=\partial
h_{\al_1}^2+\frac{4}{\al_1}\partial h_{\al_1}+
\partial
h_{\al_2}^2+\frac{4}{\al_2}\partial h_{\al_2}=4(\elu+\eld).$$ et
$$D_{h_{\al_1}^2h_{\al_2}^2}(1-k)=T_{h_{\al_1}^2}(1-k)T_{h_{\al_2}^2}(1-k)|_{\C[\a_\C]^{\w}}$$
$$=4T_{h_{\al_1}^2}(1-k)\eld|_{\C[\a_\C]^{\w}}=16\elu\eld=16\eld\elu,$$ puisque  les racines $\al_1$ et
$\al_2$ sont orthogonales. D'o\`u le r\' esultat.\end{dem}

\subsection{Solutions du syst\`eme $E_{\chi}(\are)$}

Soit $\Phi$ une fonction $H$-invariante et analytique sur $\q^{reg}$
solution du syst\`eme $(E_\chi)(\qr)$. Par $H$-invariance,
la fonction $\Phi$ est uniquement  d\' etermin\' ee par ses
restrictions \`a $\are$ pour $\a\in<car(\q)>$ et par le th\'
eor\`eme \ref{rad} et la proposition \ref{prad}, pour tout
$\a\in<car(\q)>$, la fonction $H$-invariante  $(|\de| \Phi)|\are$
est solution du syst\`eme suivant:

$$(S_{\are,\chi})\left\{\begin{array}{l}(\elu+\eld)\tilde{\Psi}=(\la_1+\la_2)\tilde{\Psi}\\
\elu\eld\tilde{\Psi}=\la_1\la_2\tilde{\Psi}\end{array}\right.\qq\mbox{
sur }\are .$$

Pour $j=1,2$, pr\' ecisons   l'expression des op\' erateurs $L_{\al_j}$ en terme des coordonn\' ees sur
$\a\in<car(\q)>$. Soit $\a=\a_{\ep_1,\ep_2}$ avec $\ep_j=\pm$.

La racine $\al_j$ est r\' eelle si et seulement si $\ep_j=+$. Dans ce cas, on a $\al_j(X_{u_1 u_2}^{\ep_1 \ep_2})=2 u_j, \quad \partial h_{\al_j}=2 \frac{\partial}{\partial u_j}$ et $\quad L_{\al_j}=\frac{\partial^2}{\partial u_j^2}+\frac{1}{u_j}\frac{\partial}{\partial u_j}$.

La racine $\al_j$ est imaginaire si et seulement si $\ep_j=-$. Dans ce cas, on a $\al_j(X_{u_1 u_2}^{\ep_1 \ep_2})=2\im  u_j, \quad \partial h_{\al_j}=-2\im  \frac{\partial}{\partial u_j}$ et $\quad L_{\al_j}=-\big(\frac{\partial^2}{\partial u_j^2}+\frac{1}{u_j}\frac{\partial}{\partial u_j}\big)$.

Pour $\a=\a_2$ et  $\al=2(\tau\pm\im\te),$ on a $\partial
h_\al=\frac{\partial}{\partial \tau} \mp\im\frac{\partial}{\partial
\te}$ et $ L_\al=\frac{1}{4}\left( \frac{\partial}{\partial \tau}
\mp\im\frac{\partial}{\partial
\te}\right)^2+\frac{1}{2(\tau\pm\im\te)} \left(
\frac{\partial}{\partial \tau} \mp\im\frac{\partial}{\partial
\te}\right).$\me

On introduit l'op\' erateur diff\' erentiel de Bessel d'une variable
complexe $L_c$ et son analogue r\' eel $L$
$$ L_c=4\left(z\frac{\partial^2}{\partial z^2}+\frac{\partial}{\partial z}\right)\quad\textrm{ et }  L=4\left(t\frac{d^2}{dt^2}+\frac{d}{dt}\right).$$

L'op\' erateur $L$ correspond \`a la partie radiale de l'op\'
erateur de Laplace pour la paire sym\' etrique $(\mathfrak{so}(2,1),
\mathfrak{so}(1,1))$ et s'obtient \`a partir de
$\left(\frac{\partial^2}{\partial
u^2}+\frac{1}{u}\frac{\partial}{\partial u}\right)$ par le changement de variable $t=u^2$.

L'\' equation $L_cy=\la y$ est de type Fuchs en $0$ et ses solutions
s'expriment \`a l'aide des fonctions de Bessel. Plus pr\' ecis\'
ement, on rappelle le r\' esultat suivant:

\begin{prop}\label{fonsol} (\cite{in} paragraphe  16.32 et \cite{fa} appendice A (4)) Pour $\la\in \C$, on pose
$$\Phi_\la(z)=\sum_{n\geq0}\frac{(\la z)^n}{4^n(n!)^2}\qq \textrm{ et }\qq w_\la(z)=\sum_{n\geq0}\frac{a(n)(\la z)^n}{4^n(n!)^2}, \textrm{ o\`u } a(x)= -2 \frac{\Gamma'(x+1)}{\Gamma(x+1)}.$$
\begin{enumerate}
\item Un syst\`eme fondamental de solutions de $L_c y=\la y$ sur $\C-\R_-$ est donn\' ee par la fonction analytique complexe $\Phi_\la$ et la fonction $W_\la(z)=w_\la(z)+\log (z)\Phi_\la(z)$ o\`u $\log$ d\' esigne la
d\' etermination principale du logarithme sur $\C-\R_-.$
\item  Un syst\`eme fondamental de solutions de $L y=\la y$ est donn\' ee par la fonction analytique r\' eelle  $\Phi_\la(t)$ et la fonction $W^r_\la(t)=w_\la(t)+\log |t|\Phi_\la(t)$
\end{enumerate}
\end{prop}

Pour  $\la\in\C^*$, on pose
$$\mathcal{S}ol(L,\la)=\{y\in\mathcal{C}^2(\R^*,\C); Ly=\la y\}$$
et
$$ \mathcal{S}ol(L_c,\la)=\{y :\C-\R_- \to\C,
\textrm{analytique} ;  L_cy=\la y\}.$$

Nous allons exprimer les solutions du syst\`eme $(S_{\are,\chi})$ pour $\a\in <car(\q)>$ en terme des fonctions $\Phi_\la$ et $W_\la$ donn\' ees dans la proposition pr\' ec\' edente.\me

Consid\' erons tout d'abord le cas $\a\in\{\app,\apm,\amm\}$. Par
$H$-invariance des fonctions consid\' er\' ees, il est naturel
d'introduire, pour
$i\in\{1,2\}$, les op\' erateurs diff\' erentiels
$$L_i:=4\left(t_i\frac{ \partial^2}{\partial t_i^2}+\frac{\partial}{\partial t_i}\right).$$
 La description pr\' ec\' edente des op\' erateurs
$L_{\al_i}$ permet d'obtenir  le lemme suivant:

\begin{lem} \label{solm} Soit $\a\in\{\app,\apm,\amm\}$. Si  $\tilde{\Psi}$ est une solution analytique $W_H(\a)$-invariante   de $(S_{\are,\chi})$   alors il existe une fonction $\Psi$ analytique sur $\{(t_1,t_2)\in\R^2/\;
(t_1-t_2)t_1t_2\ne0\}$, telle que:

\begin{enumerate}
\item $\Psi(t_1,t_2)=\Psi(t_2,t_1)$,
\item $\left\{\begin{array}{c}(L_1+L_2)\Psi(t_1,t_2)=(\la_1+\la_2)\Psi(t_1,t_2)\\
L_1L_2\Psi(t_1,t_2)=\la_1\la_2\Psi(t_1,t_2)\end{array}\right. ,$
\item Pour tout $X\in\are$, alors $\displaystyle\tilde{\Psi}(X)=\Psi(u(X), v(X))$ o\`u $\{u(X),v(X)\}$ caract\' erise la $H$-orbite de $X$ (lemme \ref{orbites}).
\end{enumerate}
\end{lem}

\begin{prop}\label{sold}  Les solutions $\Psi$ de classe
$\mathcal{C}^2$ sur $\{(t_1,t_2)\in\R^2/\; (t_1-t_2)t_1t_2\ne0\}$ du
syst\`eme
 $$\left\{\begin{array}{c}[(L_1+L_2)\Psi](t_1,t_2)=(\la_1+\la_2)\Psi(t_1,t_2)\\
L_1L_2\Psi(t_1,t_2)=\la_1\la_2\Psi(t_1,t_2)\end{array}\right..$$
sont sur chaque composante connexe de $\{(t_1,t_2)\in\R^2/\;
(t_1-t_2)t_1t_2\ne0\}$ des combinaisons lin\' eaires des fonctions
$$(t_1,t_2)\mapsto A(t_1)B(t_2)\qq\textrm{et}\qq
(t_1,t_2)\mapsto A(t_2)B(t_1),$$ o\`u $A$ et $B$ parcourent
respectivement $\mathcal{S}ol(L,\la_1)$ et $\mathcal{S}ol(L,\la_2).$

\end{prop}
\begin{dem} Soit $\Psi$ une solution du syst\`eme consid\' er\' e.
Pour $i\in\{1,2\}$, on a donc
$(L_{i}^2-(\la_1+\la_2)L_{i}+\la_1\la_2id)\Psi=0,$
et par suite,   $$\Psi\in \ker(L_{i}^2-(\la_1+\la_2)L_{i}+\la_1\la_2id)=\ker(L_{i}-\la_1
id)\oplus\ker(L_{i}-\la_2 id).$$  Par  la
proposition \ref{fonsol} ({\it 2.}) appliqu\' ee \`a $L_1$, il  existe  une famille
$(g_j)_{j\in\{ 1,\ldots 4\}}$ de fonctions d\' efinies sur $\R$
telles que
$$\Psi(t_1,t_2)=g_1(t_2)\Phi_{\la_1}(t_1)+
g_2(t_2)W^r_{\la_1}(t_1)+ g_3(t_2)\Phi_{\la_2}(t_1)+
g_4(t_2)W^r_{\la_2}(t_1).$$ En inversant la matrice correspondant au
wronskien de $t_1\mapsto\Phi_{\la_i}(t_1)$ et $t_1\mapsto
W^r_{\la_i}(t_1)$, pour $i\in\{1,2\}$, on en d\' eduit que les
fonctions $g_1,g_2,g_3$ et $g_4$ sont de classe $\mathcal{C}^2$ sur
$\{(t_1,t_2)\in\R^2/\;
(t_1-t_2)t_1t_2\ne0\}$.\me

En reprenant le syst\`eme initial et tenant compte la libert\' e
des fonctions $t_1\mapsto\Phi_{\la_i}(t_1)$ et $t_1\mapsto
W^r_{\la_i}(t_1)$, pour $i\in\{1,2\}$, on en d\' eduit que $g_1$ et
$g_2$ sont dans $\mathcal{S}ol(L,\la_2)$ et que $g_3$ et $g_4$ sont
dans
$\mathcal{S}ol(L,\la_1)$.\end{dem}

Etudions maintenant le cas particulier de $\a_2$.\me

Un \' el\' ement de $\a_2$ s'\' ecrit
$X_{\tau,\te}=\left(\begin{array}{c|c}
0 & \begin{array}{cc}\tau & -\te\\\te & \tau\end{array}\\
\hline  \begin{array}{cc}\tau & -\te\\\te & \tau\end{array} &
0\end{array}\right)$ avec $(\tau,\te)\in\R^2$.

On note $\a_2^+$ la chambre de Weyl de $\a_2^{reg}$ form\' ee des $X_{\tau,\te}$ tels que $\tau>0$ et $\te>0$. Ainsi une fonction $W_H(\a_2)$-invariante solution de $S_{\a_2^{reg},\chi}$ est uniquement d\' etermin\' ee par sa restriction \`a $\a_2^+$ et v\' erifie le m\^eme syst\`eme sur $\a_2^+$ que l'on notera $S_{\a_2^+,\chi}$.\me

Pour $i=1,2$, on note $\dis\quad  L_{c,i}:=4\left(z_i\frac{
\partial^2}{\partial z_i^2}+\frac{\partial}{\partial z_i}\right).$

\begin{prop} \label{sol2}Si  $\tilde{\Psi}$ est une solution analytique    de $(S_{\a_2^+,\chi})$   alors il existe un ouvert connexe $U$ de $\C^2$ contenant $\{(z_1,z_2)\in\C^2; z_1=\bar{z_2} \mbox{ et } Im(z_1)>0\}$ et une  fonction $\Psi_c: U\rightarrow \C$ analytique sur $U$ telle que:
\begin{enumerate}
\item $\left\{\begin{array}{c}[(L_{c,1}+L_{c,2})\Psi_c](z_1,z_2)=(\la_1+\la_2)\Psi_c(z_1,z_2)\\
L_{c,1}L_{c,2}\Psi_c(z_1,z_2)=\la_1\la_2\Psi_c(z_1,z_2)\end{array}\right.
\mbox{ sur } U,$
\item Pour $\tau>0$ et $\te>0$ alors on a $\tilde{\Psi}(X_{\tau,\te})=\Psi_c\big( (\tau+\im\te)^2,(\tau-\im\te)^2\big).$
\end{enumerate}
\end{prop}
\begin{dem} On note $ D_+=\{(x,y)\in\R^2; y>0\}$. Par le lemme \ref{orbites},
pour tout $(x,y)\in D_+$, il existe une unique $H$-orbite
$\mathcal{O}$ de $\q^{reg}$ telle que $\mathcal{O}\cap\a_2\neq
\emptyset$ et pour tout $X\in \mathcal{O}$ alors $\{u(X), v(X)\}=\{
x+\im y, x-\im y\}$. 

Ainsi l'application
$\ga$ d\' efinie sur $ \a_2^+$ par $\ga(X_{\tau,\te})=( \tau^2-\te^2, 2\te\tau)$ est un
diff\' eomorphisme analytique  de  $ \a_2^+$ dans $D_+$ d'inverse analytique donn\' e par
$\ga^{-1}(x,y)= X_{\tau,\te}$  avec $\tau=Re(\sqrt{x+\im y})$ 
et  $ \te= Im( \sqrt{x+\im y})$ o\`u la
fonction racine $\sqrt{z}=e^{\frac{1}{2}\log (z)}$ est d\' efinie sur $\C-\R_-$.

On d\' efinit la fonction $\Psi_0$ sur $D_+$ par
$\Psi_0(x,y)=\tilde{\Psi}\circ\ga^{-1}(x,y).$
Cette fonction  est analytique r\' eelle sur $D_+$ et, en
notant $L_0=(x+\im y)\left(\frac{\partial}{\partial x}-\im
\frac{\partial}{\partial y}\right)^2+2\left(\frac{\partial}{\partial
x}-\im \frac{\partial}{\partial y}\right),$ elle satisfait le
syst\`eme suivant sur $D_+$:
$$(S_0(D_+))\left\{\begin{array}{c}(L_0+\overline{L_0})\Psi_0(x, y)=(\la_1+\la_2)\Psi_0(x,y)\\
L_0\overline{L_0}\Psi_0(x,y)=\la_1\la_2\Psi_0(x,y)\end{array}\right..$$

Comme la fonction $\Psi_0$ est  d\' eveloppable en s\' erie enti\`ere au
voisinage de tout point de $D_+$, pour tout $(x_0,y_0)\in
D_+$,  il existe une famille $(a_{m,n})_{(m,n)\in\N^2}$ de nombres
complexes   et $R>0$ tels
 que:

 \no $a)$  la famille $(a_{m,n}r^{m+n})_{(m,n)\in\N^2}$ est sommable  pour tout $r<R$

 \no $b)$ pour $|x-x_0|<R$ et $|y-y_0|<R$ alors $\displaystyle \Psi_0(x,y)=\sum_{m,n\geq0}a_{m,n}(x-x_0)^m(y-y_0)^n,$
 cette expression \' etant ind\' ependante de l'ordre de sommation des termes (\cite{ho2} th\' eor\`eme 2.2.6).

Par suite, il existe une famille $(b_{m,n})_{(m,n)\in\N^2}$ de
nombres complexes telle que, pour tout $r<R$, la famille
$(b_{m,n}r^{m+n})_{(m,n)\in\N^2}$ soit sommable et, en posant
$z=x+\im y$ et $z_0=x_0+\im y_0$, pour $|z-z_0|<R$, on a
$$ \Psi_0(x,y)=\sum_{m,n\geq0}b_{m,n}(z-z_0)^m(\overline{z-z_0})^n.$$

Ceci permet de d\' efinir la fonction $\Psi_c$ sur le disque de
centre $(z_0,\overline{z_0})$ et de rayon $R$ dans $\C^2$ en posant

$$\Psi_c(z_1,z_2)=\sum_{m,n\geq0}b_{m,n}(z_1-z_0)^m(z_2-\overline{z_0})^n.$$

Cette construction est licite au voisinage de tout point
$(z_0,\overline{z_0})$ tels que $(Re(z_0),Im(z_0))\in D_+$. Ainsi,
par  prolongement analytique (\cite{ho2}  th\' eor\`eme 2.2.6 et
remarque qui suit le th\' eor\`eme 2.2.7), il existe un ouvert
connexe $U$ de $\C^2$ contenant $\{(z_1,z_2)\in\C^2; z_1=\bar{z_2}
\mbox{ et } Im(z_1)>0\}$ et une  fonction $\Psi_c: U\rightarrow \C$
analytique sur $U$ telle que, pour tout $x\in\R $ et  $y>0$, l'on
ait
$$\Psi_c(x+\im y,x-\im y)=\Psi_0(x,y)=\tilde{\Psi}\circ\ga^{-1}(x,y),$$

\no ce qui s'\' ecrit encore $\quad \tilde{\Psi}(X_{\tau,\te})=\Psi_c\big( (\tau+\im\te)^2,(\tau-\im\te)^2\big),$
pour $\tau>0$ et $\te>0$.\me

La fonction $\Psi_0$ satisfait le syst\`eme
$(S_0(D_+))$ et donc, par construction  de $\Psi_c$, pour $(x,y)\in
D_+$, on a les relations 
$$\big[(L_{c,1}+L_{c,2})\Psi_c\big](x+\im y,x-\im y)=(L_0+\overline{L_0})\Psi_0(x,y)=(\la_1+\la_2)\Psi_c(x+\im y,x-\im y)$$
$$\textrm{et } 
\big[(L_{c,1}L_{c,2})\Psi_c\big](x+\im y,x-\im y)=(L_0\overline{L_0})\Psi_0(x,y)=(\la_1\la_2)\Psi_c(x+\im y,x-\im y).$$

Maintenant, si $R$ est une fonction analytique sur $U$ alors, pour
$m,n\in\N$ et $(z,\overline{z})\in U$, on a
$$\Big(\frac{\partial^m}{\partial z_1^m}\frac{\partial^n}{\partial z_2^n}R\Big)(z,\overline{z})=\frac{\partial^m}{\partial z^m}\frac{\partial^n}{\partial \overline{z}^n}\big(z\to R(z,\overline{z})\big)$$

Ainsi, par le th\'eor\`eme 2.2.6 de \cite{ho2}, si $R(z,\overline{z})=0$ pour tout $(z,\overline{z})\in U$, alors $R=0$ sur $U$. On obtient donc  que $\Psi_c$ est solution du syst\`eme $$\left\{\begin{array}{c}[(L_{c,1}+L_{c,2})\Psi_c](z_1,z_2)=(\la_1+\la_2)\Psi_c(z_1,z_2)\\
L_{c,1}L_{c,2}\Psi_c(z_1,z_2)=\la_1\la_2\Psi_c(z_1,z_2)\end{array}\right.$$
sur $U$. Ceci ach\`eve la preuve du lemme. \end{dem}

\begin{prop}\label{sold2} Les solutions $\Psi_c$ analytiques sur $U$ du syst\`eme
$$\left\{\begin{array}{c}(L_{c,1}+L_{c,2})\Psi_c(z_1,z_2)=(\la_1+\la_2)\Psi_c(z_1,z_2)\\
(L_{c,1}L_{c,2})\Psi_c(z_1,z_2)=(\la_1\la_2)\Psi_c(z_1,z_2)\end{array}\right.$$
sont des combinaisons lin\' eaires des fonctions $$(z_1,z_2)\mapsto
A(z_1)B(z_2)\qq\textrm{et}\qq (z_1,z_2)\mapsto A(z_2)B(z_1),$$
 o\`u $A$ et $B$ parcourent
respectivement $\mathcal{S}ol(L_c,\la_1)$ et
$\mathcal{S}ol(L_c,\la_2).$
\end{prop}

La preuve de cette proposition  est identique \`a celle de la
proposition  \ref{sold}.

\subsubsection{ Actions sur les int\' egrales orbitales}
Nous concluons cette partie en donnant l'action des op\'erateurs
diff\'erentiels sur les int\'egrales orbitales.

\begin{lem} \label{intrad2} Soit $f\in\D(\q)$ et $\a\in car(\q)$. On note toujours $\al_1$ et $\al_2$ les racines positives de multiplicit\'e $1$ de $\a$. On a alors:

$$(\mathcal{M}_H(\partial(Q)f)|\a^{reg}=(L_{\al_1}+L_{\al_2})(\mathcal{M}_H(f)|\a^{reg})$$
et $$(\mathcal{M}_H(\partial(S)f)|\a^{reg}=(L_{\al_1} L_{\al_2})(\mathcal{M}_H(f)|\a^{reg}).$$

\end{lem}

\begin{dem} Ceci d\'ecoule imm\'ediatement de la normalisation de $\mathcal{M}_H(f)$ et de la proposition \ref{prad} \end{dem}

On reprend les notations de la d\'efinition \ref{notHinv},
c'est-\`a-dire si $X\in \mq^{reg}$ alors
$\mathcal{M}_H(f)(X)=\Mfm(u(X),v(X))$ et si $X=X_{\tau,\te}\in\a_2$
alors $\mathcal{M}_H(f)(X_{\tau,\te})=\Mf2(\tau,\te)$. On note
$\al=2(\tau+\im\te)$ et $\overline{\al}$ les racines positives de
multiplicit\'e $1$ de $\a_2$. On obtient alors:

\begin{cor} \label{intrad} Soit $f\in\D(\q)$ alors on a

$$\mathcal{M}\big(\partial(Q)f\big)_\m=(L_1+L_2)\Mfm \mbox{ et } \mathcal{M}\big(\partial(S)f\big)_\m=L_1L_2 \Mfm$$
et
$$\mathcal{M}\big(\partial(Q)f\big)_2=(L_\al+L_{\overline{\al}})\Mf2 \mbox{ et } \mathcal{M}\big(\partial(S)f\big)_2=L_\al L_{\overline{\al}} \Mf2.$$

\end{cor}

 \section{ Distributions propres invariantes de $L^1_{loc}(\mathcal{U})^H$ et perspectives sur $\Lu$ }

On note $\mathcal{N}$ l'ensemble des \'el\'ements nilpotents de $\q$
et on pose $\mathcal{U}=\q-\mathcal{N}.$\\

 Le but de cette partie est de d\' ecrire l'ensemble des distributions invariantes sur $\mathcal{U}$, propres
 pour un caract\` ere r\' egulier $\chi$ de $\C[\q]^H$ fix\'e,
  donn\' ees par une fonction $\Phi$  localement
int\' egrable $H$-invariante sur $\mathcal{U}$.\\

 On rappelle que pour $T$ une distribution  $H$-invariante   propre sur ${\mathcal U}$
solution du syst\` eme
$$(E_\chi)\left\{\begin{array}{c}\partial(Q)T=\chi(Q)T\\\partial(S)T=\chi(S)T\end{array}\right..$$

\no alors la restriction de $T$ \`a l'ensemble ${\mathcal U}^{reg}=\q^{reg}$  est une fonction analytique  (\cite{se} th\'eor\`eme 5.3(i) de ).\me

Nous allons d\'ecrire les conditions n\' ecessaires et suffisantes
sur $\Phi\in L^1_{loc}(\mathcal{U})^H$ satisfaisant  le syst\`eme
$(E_\chi)(\qr)$ pour que la distribution $T_\Phi$ d\' efinie
sur ${\mathcal D}({\mathcal U})$ par $\dis \langle T_\Phi, f\rangle
=\int_\q \Phi(X) f(X) dX$
soit propre invariante sur $\mathcal{U}$. \me

 Par hypoth\`ese sur   $\Phi$ ,  pour
 $P\in \C[\q]^H$ et $f\in {\mathcal D}({\mathcal U})$, on a $\langle
T_{\partial(P)\Phi},f\rangle=\chi(P) \langle T_\Phi,f\rangle$ ainsi
$T_\Phi$ est propre sur ${\mathcal U}$ pour le caract\` ere $\chi$
si et seulement si pour tout $P\in  \C[\q]^H$ et $f\in {\mathcal
D}({\mathcal U})$, on a
  $$\int_\q[( \partial(P)\Phi) f-\Phi (\partial(P)f)](X) dX=0$$

La formule d'int\' egration de Weyl permet un travail sur chaque
sous-espace de Cartan de $<car(\q)>$.   Compte-tenu du comportement
des int\' egrales orbitales au voisinage des points semi-r\'
eguliers obtenu dans la partie   4  et de la description de   $
(|\delta |\Phi)|a_{reg}$ ( voir  le lemme \ref{solm}   et les
propositions \ref{sold}, \ref{sol2} et \ref{sold2}), une int\'
egration par parties sur chaque sous-espace de Cartan  conduira aux
conditions voulues sur $\Phi$. Ces conditions sont analogues aux
conditions
 de recollement obtenues dans le cas des alg\` ebres r\' eductives par
Harish-Chandra et T. Hira\" i (\cite{hc} et \cite{ hi}) et de celles
de J. Faraut obtenues pour un hyperbolo\" ide \` a une nappe  \cite{
fa3}).

Elles traduisent le comportement de $\Phi$ au voisinage des points
semi-r\' eguliers. Gr\^ace \` a l'application $\varpi$ d\' efinie
dans le paragraphe \ref{cartans} qui renverse l'ordre d'Hira\" i  et la
$H$-invariance des fonctions consid\' er\' ees, on ram\`ene l'\'
etude de $\Phi$  au voisinage des points de $\app\cap\apm$ et
$\app\cap\a_2$. Tout comme pour l'\' etude des int\' egrales
orbitales, ces deux cas
n\' ecessitent une m\' ethode sp\' ecifique.

\subsection{Voisinages invariants des \' el\' ements semi-r\' eguliers}

On rappelle que $\varpi$ d\' esigne l'isomorphisme $H$-\'
equivariant de $\q$ donn\' e  par $\varpi\left(\begin{array}{c|c} 0
& Y\\ \hline Z & 0\end{array}\right)=\left(\begin{array}{c|c} 0 &
Y\\ \hline -Z & 0\end{array}\right)$. Cet isomorphisme renverse
l'ordre d'Hira\" i sur $car(\q).$ Par l'\' etude du paragraphe
\ref{semireg}, tout \' el\' ement semi-r\' egulier appartient \`a
 $H\cdot\p\mq$, ou \`a $H\cdot\ztp$ ou \`a
$H\cdot\varpi(\ztp)$. Par l'\' etude de la partie 3, ces ensembles
sont ouverts puisqu'image d'un ouvert par une application
submersive.\me

On rappelle que les polyn\^omes $Q$, $S$ et $S_0$ sont donn\' es par
$Q(X)=\frac{1}{2} tr(X^2)$, $S(X)=det(X)$ et $S_0= Q^2-4S$.

\begin{lem}\label{ouvert} Nous avons :\begin{enumerate}
\item $$H\cdot\p\mq=\{X\in\q,\; S_0(X)>0\}$$
\item $$H\cdot\ztp=H\cdot\p\z_3\cap\q=\{X\in\q,\;S(X)>0,\; Q(X)>-2\sqrt{|S(X)|}\}$$
\item $$H\cdot\vi(\ztp))=\{X\in\q,\;S(X)>0,\; Q(X)<2\sqrt{|S(X)|}\}.$$
\item $$\mathcal{U}=H\cdot\p\mq\cup H\cdot\ztp\cup H\cdot\vi(\ztp))$$
$$=\{X\in\q,\;
Q(X)\ne0 \;\textrm{ou}\; S_0(X)\ne0\}.$$

\end{enumerate}
\end{lem}

\begin{dem} Soit $X\in\q$ et soit  $X=X_s+X_n$ sa d\' ecomposition de
Jordan. En particulier $[X_s,X_n]=0$ et donc
$X^2=X_s^2+X_n(X_n+2X_s)$ o\`u $X_n(X_n+2X_s)$ est nilpotent et
commute \` a $X_s^2$. On obtient donc $Q(X)=Q(X_s)$. Il est imm\'
ediat que $S(X)=S(X_s)$ et par suite
$S_0(X)=S_0(X_s)$. 

Par ailleurs, par la remarque \ref{rempoly}, on a  $S_0(\mq)\subset
\R_+$ et $S_0(\a_2)\subset \R_-$, le polyn\^ome $S$ prend des
valeurs positives sur $\app, \amm$ et $\a_2$ et des valeurs n\'
egatives sur $\apm$, le polyn\^ome $Q$ prend des valeurs positives
sur $\app$ et n\' egatives sur $\amm$.\me

\no{\it 1.}   L'inclusion $H\cdot\p\mq\subset\{X\in\q,\; S_0(X)>0\}$ est
imm\' ediate.

Soit $X$ dans $\q$ tel que $S_0(X)>0$ et soit
$X=X_s+X_n$ sa d\' ecomposition de Jordan. Puisque
$S_0(X_s)=S_0(X)>0$, il existe $h\in H$ et $\a\in <car(\mq)>$ tel
que $h\cdot X_s\in\a$. Si $X_n= 0$ alors on obtient $X\in
H\cdot\mq$. Si $X_n\neq 0$ alors $X_s$ est semi-r\' egulier et donc
$h\cdot X_s\in\R H_1$ ou $\R\varpi(H_1)$. Comme $h\cdot X_n$ commute
\` a $h\cdot X_s$, on obtient $h\cdot X_n\in\mq$ par d\' efinition
de $\mq$. Comme $\p\mq=\{X\in\mq; S_0(X)\neq 0\}$, on obtient donc
$X\in H\cdot \p\mq$.\me

\no{\it 2.}   Par le lemme \ref{fibt}, on a $H\cdot \p\zt=H\cdot\ztp$ et $\ztp=\{ \left
(\begin{array}{c|c} 0&Y\\ \hline Y&0\end{array}\right );
\big(tr(Y)\big)^2 det(Y)>0\}$.

Montrons d'abord que
$\quad H\cdot\ztp\subset\{X\in\q,\;S(X)>0,\; Q(X)>-2\sqrt{|S(X)|}\}.$ 
Comme ces deux ensembles sont $H$-invariants, il suffit de montrer
que $\ztp$ est inclus dans $\{X\in\q,\;S(X)>0,\;
Q(X)>-2\sqrt{|S(X)|}\}.$

Si  $X=\left (\begin{array}{c|c} 0&Y\\ \hline Y&0\end{array}\right
)\in \ztp$ alors $S(X)=det(Y)^2>0$. On note $\la_1$ et $\la_2$ les
valeurs propres de $Y$. Elles sont soient r\' eelles soient
complexes conjugu\' ees et par hypoth\`ese,  on a
$(\la_1+\la_2)^2\la_1\la_2>0$. Ainsi, on obtient la relation
$Q(X)=\la_1^2+\la_2^2>-2|\la_1\la_2|=-2\sqrt{S(X)}.$\me

Montrons l'inclusion inverse. Soit $X$ dans $\q$ tel que $S(X)>0$ et
$Q(X)>-2\sqrt{|S(X)|}$. Soit $X=X_s+X_n$ sa d\' ecomposition de
Jordan. On a donc $S(X_s)>0$ et $Q(X_s)>-2\sqrt{|S(X_s)|}$.  Ainsi,
il existe $h\in H$ et $\a\in\{ \app,\a_2,\amm\}$ tels que $h\cdot
X_s\in \a$. Si $Q(X_s)\geq 0$ alors $\a\in\{\app,\a_2\}$. Si
$Q(X_s)<0$ alors $S_0(X_s)<0$ et donc $\a=\a_2$. Ainsi, on a
toujours $h\cdot X_s\in\zt$.

Si $X_n=0$ on obtient alors $h\cdot X\in\zt$. Si $X_n\neq 0$ alors
$X_s$ est semi-r\' egulier et donc $h\cdot X_s\in\R H_3$. Comme
$h\cdot X_n$ commute \` a $h\cdot X_s$, il appartient \` a $\zt$ par
d\'
efinition de $\zt$. On conclut donc $h\cdot X\in\p\zt$.\me

\no{\it 3.}   On a
$H\cdot\varpi(\ztp)=\varpi(H\cdot\ztp).$  Comme $S(\varpi(X))=S(X)$ et
$Q(\varpi(X))=-Q(X)$, on obtient $H\cdot\vi(\ztp))=\{X\in\q,\;S(X)>0,\; Q(X)<2\sqrt{|S(X)|}\}$.\me

\no$4.$   Soit $X\in\mathcal{U}.$ Alors $X_s$ est non nul et  les valeurs
$Q(X)=Q(X_s)$ et  $S_0(X)=S_0(X_s)$ ne sont pas simultan\' ement
nulles.
 Supposons que $S_0(X)$ soit non nul. Si
$S_0(X)>0$, alors $X$ est dans  $H\cdot\p\mq$.

Si $S_0(X)<0$, alors on a $Q(X)^2-4S(X)<0$ ce qui donne  $S(X)>0$ et
$(Q(X)-2\sqrt{S(X)})(Q(X)+2\sqrt{S(X)})<0$. Comme
$Q(X)-2\sqrt{S(X)}<Q(X)+2\sqrt{S(X)}$,  on obtient
$Q(X)-2\sqrt{S(X)}<0<Q(X)+2\sqrt{S(X)}$, ce qui implique que $X$ est
dans $H\cdot \ztp \cap H \cdot\varpi(\ztp)$.

Si $S_0(X)=0$, alors $Q(X)$ est non nul et on a $S(X)=Q(X)^2/4>0$ et
$Q(X)=\pm 2\sqrt{S(X)}$. N\'ec\'essairement $X$ est dans
$H\cdot\varpi(\ztp)$ ou $H\cdot\ztp.$

La r\'eciproque est imm\'ediate. \end{dem}

 \subsection{Conditions de recollement}\label{CondRecol}

Nous rappelons  les notations suivantes introduites  dans le paragraphe \ref{OS} et la d\'efinition \ref{notHinv}. \me

Si  $X=\left(\begin{array}{cc} 0& Y\\Z& 0\end{array}\right)$ alors $u(X)$ et $v(X)$ d\'esignent les valeurs propres de $YZ$ et elles caract\' erisent la $H$-orbite de $X$  lorsque $X$ est semi-simple. 

 Si $F$ est une fonction $H$-invariante sur $\q^{reg}$, les  fonctions $F_\m$  d\' efinie sur $(\R^*)^2-diag$,  et $F_2$ d\' efinie sur $(\R^*)^2$ v\' erifient, pour $X\in \mq^{reg}$ alors  $\dis F(X)=F_\mathfrak{m}(u(X),v(X))$  et pour $X_{\tau,\te}\in\a_2$ alors $F(X_{\tau,\te})=F_2(\tau,\te)$.
 Pour $f\in{\mathcal D}(\q)$, on note $\big( \mathcal{M}_H f\big)_\mathfrak{m}=\Mfm$ et $ \big( \mathcal{M}_Hf\big)_2=\Mf2$.\\

Dans tout ce paragraphe, on se placera sous les hypoth\` eses suivantes.\me

 Soit  $\Phi\in L^1_{loc}(\mathcal{U})^H$ une   solution du syst\`eme $(E_\chi)(\q^{reg})$.  Pour conserver les propri\' et\' es de sym\' etrie et anti-sym\'
etrie des fonctions, on exprime le syst\`eme $(E_\chi)$  \`a l'aide
de  la base $\{ Q,S_0\}$ de $\C[\q]^H$.

On note $\tilde{\Psi}=(|\delta| \Phi)|\qr$. Par  l'\' etude du
paragraphe pr\' ec\' edent, pour chaque $\a\in car(\q)$, la fonction
$\tilde{\Psi}|\are$ est une fonction analytique $W_H(\a)$-invariante
solution du syst\`eme

$$(S_{\are,\chi})\left\{\begin{array}{l}(\elu+\eld)\tilde{\Psi}=(\la_1+\la_2)\tilde{\Psi}\\
(\elu-\eld)^2\tilde{\Psi}=(\la_1-\la_2)^2\tilde{\Psi}\end{array}\right.\qq\mbox{
sur }\are ,$$

\no o\`u $\al_1$ et $\al_2$ d\'esignent les racines positives de
multiplicit\'e $1$ de $\a$.

 Par le lemme \ref{solm} et la proposition \ref{sol2} , il existe une fonction
$\Psi_\m$ analytique sur $\{(t_1,t_2)\in\R^2/\;
(t_1-t_2)t_1t_2\ne0\}$ et sym\' etrique par rapport aux deux
variables telle que
$$\left\{\begin{array}{c}[(L_1+L_2)\Psi_\m](t_1,t_2)=(\la_1+\la_2)\Psi_\m(t_1,t_2)\\
\left[(L_1-L_2)^2\Psi_\m\right](t_1,t_2)=(\la_1-\la_2)^2\Psi_\m(t_1,t_2)\end{array}\right.,$$

\noindent  et une fonction $\Psi_c$ analytique sur un ouvert $U$
connexe de $\C^2$ contenant $\{(z_1,z_2)\in\C^2; z_1=\bar{z_2}
\mbox{ et } Im(z_1)>0\}$ telle que
$$\left\{\begin{array}{c}[(L_{c,1}+L_{c,2})\Psi_c](z_1,z_2)=(\la_1+\la_2)\Psi_c(z_1,z_2)\\
\left[(L_{c,1}-L_{c,2})^2\Psi_c\right](z_1,z_2)=(\la_1-\la_2)^2\Psi_c(z_1,z_2)\end{array}\right.
\mbox{ sur } U,$$

avec  pour tout $X\in\qr$:

 $$\tilde{\Psi}(X)=\left\{\begin{array}{ll} \Psi_\m(u(X),v(X))& \mbox{ si }
X\in H\cdot\mq^{reg}\\
\Psi_c(u(X),v(X))&\mbox{ si } X\in H\cdot \a_2^+
\end{array}\right.,$$

\noindent  Plus pr\' ecis\' ement, l'\' el\' ement
$X=X_{\tau,\te}$ appartient \` a $\a_2^+$ si et seulement si   $
\tau>0$ et $ \te>0$, et dans ce cas on a
$\tilde{\Psi}(X_{\tau,\te})=\Psi_c\big( (\tau+\im
\te)^2,(\tau-\im\te)^2\big)$.

  Pour $(\tau,\te)\in(\R^*)^2$, on note $\Psi_2(\tau,\te)= \tilde{\Psi} (X_{\tau, \te}).$ Par invariance de
$\tilde{\Psi}$ sous l'action du groupe de Weyl $W_H(\a_2)$, la
fonction $\Psi_2$ est  paire
 par rapport \`a chaque variable.
  \begin{Def} \label{ipp}   Soient $f$ et $g$ deux fonctions d\' efinies sur un ouvert de $\R^2$. Soit $D$ un op\' erateur
  diff\' erentiel sur $\R^2$ tel  que $Df$ et $Dg$ existent. On notera, lorsque
les int\' egrales consid\' er\' ees convergent
 $$I_\m(D,f,g)=\int_{t_1>t_2} \big( f\; (Dg)-(Df)\; g)(t_1, t_2) dt_1 dt_2$$
$$ I_2(D,f,g)=\int_{\R^2} (\tau^2+\te^2)\big( f\; (Dg)-(Df)\;
g)(\tau,\te) d\tau d\te$$
\end{Def}

Par la formule d'int\' egration de Weyl (lemme \ref{intweyl}) et
l'action des op\' erateurs diff\' erentiels sur les int\' egrales
orbitales (corollaire  \ref{intrad}), on a
$$\begin{array}{c}\langle \partial(Q) T_{\Phi},f\rangle-\langle T_{\partial(Q)
\Phi},f\rangle=I_\m(L_1+L_2,\Psi_\m,\Mfm)+2 I_2(L_\al+L_{\overline{\al}},\Psi_2,\Mf2),\\
\textrm{  et }\langle \partial(S_0) T_{\Phi},f\rangle-\langle T_{\partial(S_0)
\Phi},f\rangle=I_\m\big((L_1-L_2)^2,\Psi_\m,\Mfm)+2
I_2(\big(L_\al-L_{\overline{\al}}\big)^2,\Psi_2,\Mf2).\end{array}$$

Ainsi, $T_\Phi$ est une distribution propre invariante sur
${\mathcal U}$ pour le caract\`ere $\chi$ si et seulement si, pour
tout $f\in\D({\mathcal U})$, les fonctions $\Psi_\m$ et $\Psi_2$
satisfont les relations  suivantes:

$$({\mathcal Rec})\left\{\begin{array}{l}  I_\m(L_1+L_2,\Psi_\m,\Mfm)+2 I_2(L_\al+L_{\overline{\al}},\Psi_2,\Mf2)=0\\
I_\m\big((L_1-L_2)^2,\Psi_\m,\Mfm)+2
I_2(\big(L_\al-L_{\overline{\al}}\big)^2,\Psi_2,\Mf2)=0\end{array}\right.$$

Dans les paragraphes suivants, nous allons \' etudier ces relations
successivement pour $f\in\D(H\cdot\p\mq)$, puis $f\in\D(H\cdot\ztp)$
et enfin $f\in\D(H\cdot\varpi(\ztp))$.

\subsubsection{ Recollement sur $H.\p\mq$}

Par le lemme \ref{ouvert}, si  $f\in \D(H\cdot\p\mq)$ alors le
support de $f$ est contenu dans $\{X\in\q;S_0(X)>0\}$. Comme
$H\cdot\a_2\subset \{X\in\q;S_0(X)\leq 0\}$,  la fonction  $\Mf2$ est
identiquement nulle. 

D'autre part, par le th\' eor\` eme \ref{iom}, l'application
$f\mapsto \Mfm$ est surjective de $\D(H\cdot\p\mq)$ dans
$\mathcal{H}_{\log}^2$. Ainsi,  les relations de recollement
$({\mathcal Rec})$ pour $f\in\D(H\cdot\p\mq)$ sont \'equivalentes
\`a $$I_\m(L_1+L_2,\Psi_\m,F)=0 \textrm{ et }
I_\m\big((L_1-L_2)^2,\Psi_\m,F)=0, \textrm{   pour tout }
F\in\mathcal{H}_{\log}^2.$$

Pour
simplifier les notations, on pose dans tout ce paragraphe $\Psi=\Psi_\m$. \me

Une fonction de $\mathcal{H}_{\log}^2$ s'\' ecrit sous la forme
$a(t_1,t_2)+\log|t_1|b(t_1,t_2)+\log|t_2|b(t_2,t_1)+\log|t_1|\log|t_2|c(t_1,t_2),$
avec  $a$ et $c$  sym\' etriques par rapport aux deux variables
et
$a,b$ et $c$ dans $\D(\R^2-diag)$.

Par le lemme \ref{solm} et la proposition \ref{sold},
sur chaque composante connexe de $(\R^*)^2-diag$, la fonction $\Psi$
s'\' ecrit comme somme de fonctions de  la forme
$(A(t_1)+\log|t_1|B(t_1))(C(t_2)+\log|t_2|D(t_2)),$
o\`u $ A,B,C$ et $D$ sont des fonctions analytiques sur $\R$.\me

Pour exprimer les conditions de recollement, nous introduisons les
notations
suivantes.\me

Soit $\E$ l'espace des fonctions $u$ de la forme
$u(t)=v(t)+\log|t|w(t),$ o\`u $v$ et $w$ sont de classe
$\mathcal{C}^2$ sur $\R^*$ et  admettent, ainsi que leurs d\' eriv\'
ees, des limites  \` a droite et \` a gauche   en $0$.  Pour une
fonction $f$ continue sur $\R^*$, on notera lorsqu'elles existent
$f(0^+)$ et $f(0^-)$ ses limites \` a droite et \` a gauche en $0$.\me

Pour $u\in\E$, on pose
$u^{[1]}(t)=tu'(t)$ et $u^{[0]}(t)=u(t)-\log|t|u^{[1]}(t),$
de telle sorte que
$$\lim_{t\to0^{\pm}}u^{[1]}(t)=w(0^\pm)\qq\textrm{et}\space
\lim_{t\to0^{\pm}}u^{[0]}(t)=v(0^\pm).$$ 
Ainsi, pour $h$ et $u$ dans $\E$ et $\displaystyle
L=4(t\frac{\partial^2}{\partial t^2}+\frac{\partial}{\partial t})$,
 on a
$$  \lim_{t\to0^{\pm}}(u^{[1]}h^{[0]}-u^{[0]}h^{[1]})(t)=(u^{[1]}h^{[0]}-u^{[0]}h^{[1]})(0^{\pm})=
\lim_{t\to0^{\pm}}(t(u'h-uh')(t))$$
$$ \textrm{ et } \quad 
h(t)(Lu)(t)-(Lh)(t)u(t)=4\frac{d}{dt}[h^{[0]}u^{[1]}-h^{[1]}u^{[0]}](t)=4\frac{d}{dt}[t(u'h-uh')(t)].$$

Pour   $f$ et $g$ deux fonctions de classe $\mathcal{C}^1$ sur un
ouvert $U$  de $\R^2$ et $j\in\{1,2\} $, on d\' efinit les op\' erateurs $K_j$ par
$$K_{j}(f,g)(t_1,t_2)=t_j\left(f\; \frac{\partial g}{\partial
t_j}-\frac{\partial f}{\partial t_j}\;g\right)(t_1,t_2).$$ 
En notant  $\displaystyle
L_j=4\big( t_j\frac{\partial^2}{\partial
t_j^2}+\frac{\partial}{\partial t_j}\big)$, on a donc 
$\quad f(L_jg) - (L_jf)g=4\frac{\partial} {\partial t_j}K_j(f,g).$\me

Par ailleurs,  pour $f$ une fonction d\' efinie sur
$\R^{n_1}\times\R^{n_2}$ o\`u  $(n_1, n_2)\in(\N^*)^2$, et
$(x,y)\in\R^{n_1}\times\R^{n_2}$, on note $f_x$ la fonction d\'
efinie sur $\R^{n_2}$ par $f_x(y)=f(x,y)$ et $f^y$ la fonction d\'
efinie sur $\R^{n_1}$ par $f^y(x)=f(x,y).$

\begin{rem}\label{continue} Soit $F\in\mathcal{H}_{\log}^2$.
 Par la description de $\mathcal{H}_{\log}^2$ et les hypoth\`eses sur $\Psi$,  les fonctions
$F_{t_1}$ et $\Psi_{t_1}$ pour $t_1\in\R^*$ et les fonctions
$F^{t_2}$ et $\Psi^{t_2}$ pour $t_2\in\R^*$ appartiennent \` a $\E$.

En particulier, pour tout $t\in\R^*$ et $j\in\{0,1\}$ , les limites
$(\Psi_t)^{[j]}(0^\pm)$ et $(\Psi^t)^{[j]}(0^\pm)$ existent.

Par ailleurs, pour tout $t\in\R^*$ et $j\in\{0,1\}$, les fonctions
$x\mapsto (F_t)^{[j]}(x)$ et $x\mapsto (F^t)^{[j]}(x)$ sont
continues en $0$.

\end{rem}

\begin{rem}\label{K} Pour $F\in\mathcal{H}_{\log}^2$. on a

$$K_1(\Psi,F)(t_1,t_2)=((\Psi^{t_2})^{[0]}(F^{t_2})^{[1]}-(\Psi^{t_2})^{[1]}(F^{t_2})^{[0]})(t_1)$$
et
$$K_2(\Psi,F)(t_1,t_2)=((\Psi_{t_1})^{[0]}(F_{t_1})^{[1]}-(\Psi_{t_1})^{[1]}(F_{t_1})^{[0]})(t_2).$$

\end{rem}

\begin{rem}\label{symetrie}
\begin{enumerate}
\item Les fonctions $\Psi$ et $F$ \' etant sym\' etriques, pour tout $(t_1,t_2)\in(\R^*)^2-diag$ et pour $j\in\{0,1\}$, elles v\' erifient  les relations
$$(F_{t_1})^{[j]}(t_2)=(F^{t_2})^{[j]}(t_1)\;\;\;\textrm{et}\;\;\;(\Psi_{t_1})^{[j]}(t_2)=(\Psi^{t_2})^{[j]}(t_1)$$
et donc
 on a $$K_1(\Psi,F)(t_1,t_2)=K_2(\Psi,F)(t_2,t_1).$$
\item  Comme $F$ est une fonction  \`a support born\' e, nulle au voisinage de tout \' el\' ement diagonal  de   $\R^2$, il en est de m\^eme des fonctions $K_j(\Psi,F)$ pour $j\in\{1,2\}.$

\end{enumerate}
\end{rem}

\begin{lem}\label{ippqm} Pour tout
$F\in\mathcal{H}_{\log}^2$, les fonctions $t\mapsto
K_2(\Psi,F)(t,0^{\pm})$ sont int\' egrables sur $\R$ et on a
$$I_\m(L_1+L_2,\Psi,F)=4\int_\R[K_2(\Psi,F)(t,0^{-})-K_2(\Psi,F)(t,0^{+})] dt.$$
\end{lem}

\begin{dem} On a \begin{eqnarray*}
I_\m(L_1+L_2,\Psi,F)&=&\int_{t_1>t_2}\left(\Psi(L_1+L_2)F-
F(L_1+L_2)\Psi\right)(t_1,t_2) dt_1
dt_2\\
&=&\int_{t_1>t_2}\left([\Psi L_1F- FL_1\Psi]+[\Psi L_2F-
FL_2\Psi]\right)(t_1,t_2) dt_1 dt_2.
\end{eqnarray*}

Soit $k\in\{1,2\}$. Montrons que $ \Psi L_kF-
FL_k\Psi=4\frac{\partial} {\partial t_k}K_k(\Psi,F)$ est
int\'egrable sur $\R^2$. Pour une fonction $f$ de classe $\mathcal{C}^2$ sur un ouvert de $\R^*$, nous
avons l'\' egalit\' e
$L(\log|t|f(t))=8f'(t)+4\log|t|(f'(t)+tf''(t)).$
Ainsi, pour chaque composante connexe ${\mathcal C}$ de
$(\R^*)^2-diag$, il existe des fonctions $G_{i,j}\in {\mathcal
D}(\R^2-diag)$ o\`u $0\leq i,j\leq2$  telles que, pour $(t_1,t_2)\in
{\mathcal C}$, on ait
 $$\Big(\Psi \big(L_kF\big)- F\big(L_k\Psi\big)\Big)(t_1,t_2)=\sum_{0\leq
i,j\leq2}G_{i,j}(t_1,t_2)(\log|t_1|)^{i}(\log|t_2|)^j.$$ Par suite,
la fonction  $\dis
\Psi L_kF- FL_k\Psi=4\frac{\partial} {\partial t_k}K_k(\Psi,F)$ est int\'egrable sur $\R^2$.  Les propri\'et\'es sur le support de $F$ ( remarque \ref{symetrie} (2.)) et  le th\' eor\`eme de Fubini assurent l'int\' egrabilit\' e des fonctions $K_k(\Psi,F)(t,0^\pm)$.
On obtient donc
$$
I_\m(L_1+L_2,\Psi,F)$$
$$=\int_{t_1>t_2}[\Psi L_1F- FL_1\Psi](t_1,t_2)
dt_1 dt_2+\int_{t_1>t_2}[\Psi L_2F- FL_2\Psi](t_1,t_2) dt_1 dt_2$$
$$
=4\int_{t_1>t_2}\frac{\partial} {\partial t_1}K_1(\Psi,F)(t_1,t_2)
dt_1 dt_2+4\int_{t_1>t_2}\frac{\partial} {\partial
t_2}K_2(\Psi,F)(t_1,t_2) dt_1 dt_2.
$$

Par ailleurs, le  th\' eor\`eme de Fubini et  les propri\'et\'es sur
le support de $F$ donne 
$\dis \int_{t_1>t_2>0}\frac{\partial} {\partial t_1}K_1(\Psi,F)(t_1,t_2)
dt_1 dt_2=0,$ et donc,  on obtient
\begin{eqnarray*}& & \hspace{2,5cm} \int_{t_1>t_2}\frac{\partial}
 {\partial t_1}K_1(\Psi,F)(t_1,t_2) dt_1
dt_2 \\
&=&\int_{t_1>0>t_2}\frac{\partial} {\partial
t_1}K_1(\Psi,F)(t_1,t_2) dt_1 dt_2+\int_{0>t_1>t_2}\frac{\partial}
{\partial t_1}K_1(\Psi,F)(t_1,t_2) dt_1 dt_2\\
&=&-\int_{-\infty}^0K_1(\Psi,F)(0^+,t)dt+\int_{-\infty}^0K_1(\Psi,F)(0^-,t)dt,
\end{eqnarray*}

Par la sym\'etrie des fonctions
$\Psi$ et $F$ (remarque \ref{symetrie}), on a $K_1(\Psi,F)(t_1,t_2)=K_2(\Psi,F)(t_2,t_1)$
ce qui implique
$$\int_{t_1>t_2}\frac{\partial} {\partial t_1}K_1(\Psi,F)(t_1,t_2)
dt_1
dt_2=-\int_{-\infty}^0K_2(\Psi,F)(t,0^+)dt+\int_{-\infty}^0K_2(\Psi,F)(t,0^-)dt.$$											

De m\^eme, on obtient la relation

$$\int_{t_1>t_2}\frac{\partial} {\partial t_2}K_2(\Psi,F)(t_1,t_2)
dt_1
dt_2=-\int^{+\infty}_0K_2(\Psi,F)(t,0^+)dt+\int^{+\infty}_0K_2(\Psi,F)(t,0^-)dt.$$

Ces relations donnent alors 
le r\' esultat voulu.\end{dem}

\begin{prop}\label{recolm} Les trois assertions suivantes sont \' equivalentes:
\begin{enumerate}

\item Pour tout $F\in  \mathcal{H}_{\log}^2$, on a $I_\m(L_1+L_2,\Psi,F)=0$,

\item Pour $j\in\{0,1\}$ et $t_1\in\R^*$ alors
$(\Psi_{t_1})^{[j]}(0^{+})=(\Psi_{t_1})^{[j]}(0^{-}).$

\item Pour tout $f\in\D(H\cdot\p\mq)$ et pour tout $t\in\R^*$, on a

$$K_2(\Psi,\Mfm)(t,0^+)=K_2(\Psi,\Mfm)(t,0^-).$$

\end{enumerate}
\end{prop}

\begin{rem}\label{symK} Gr\^ace \`a la remarque \ref{symetrie}, on
constate que
\begin{enumerate} \item le point {\it 2.} de la proposition pr\'ec\'edente est
\'equivalent \`a l'assertion suivante:

({\it 2'.}) pour $j\in\{0,1\}$ et $t_2\in\R^*$ alors
$(\Psi^{t_2})^{[j]}(0^{+})=(\Psi^{t_2})^{[j]}(0^{-}).$

\item le point {\it3.} de la proposition pr\'ec\'edente est
\'equivalent \`a l'assertion suivante:

({\it 3'.}) pour tout $f\in\D(H\cdot\p\mq)$ et pour tout $t\in\R^*$,
on a

$$K_1(\Psi,\Mfm)(0^+,t)=K_1(\Psi,\Mfm)(0^-,t).$$
\end{enumerate}
\end{rem}\me

\no  $Preuve\; de\; la\; proposition\; \ref{recolm} :$  Supposons que $I_\m(L_1+L_2, \Psi,F)=0$ pour tout $F\in \mathcal{H}_{\log}^2$.
Par le   lemme \ref{ippqm}, ceci s'\' ecrit 
$$ \int_\R[K_2(\Psi,F)(t,0^-)-K_2(\Psi,F)(t,0^+)]dt
=0,$$ ce qui s'\' ecrit \' egalement, gr\^ace \`a la  remarque
\ref{K}

$$ \int_\R\Big[(F_{t_1})^{[1]}(0)\big((
\Psi_{t_1})^{[0]}(0^-)-(\Psi_{t_1})^{[0]}(0^+)\big)
 -(F_{t_1})^{[0]}(0)\big(
(\Psi_{t_1})^{[1]}(0^-)-(\Psi_{t_1})^{[1]}(0^+)\big)\Big]
dt_1=0.$$

 Nous allons appliquer cette \' egalit\' e \` a des fonctions test pour obtenir les
relations voulues sur $\Psi$.
Soit   $\chi\in\D(\R^*)$ quelconque. Il existe $\varepsilon\in]0,1[$ tel que
$supp(\chi)\subset[-\varepsilon^{-1}, -\varepsilon]\cup[\varepsilon, \varepsilon^{-1}].$ 
Soit $\varphi\in\D(\R)$ telle  que
$supp(\varphi)\subset[\frac{\varepsilon}{2}, 2\varepsilon^{-1}]$ et
$\varphi\equiv1$ sur
$[\varepsilon;\varepsilon^{-1}].$
On d\' efinit la fonction $A$ de $\mathcal{H}_{\log}^2$ par
$$A(t_1,t_2)=\chi(|t_1-t_2|)(\varphi(t_1)+\varphi(t_2))+\chi(-|t_1-t_2|)(\varphi(-t_1)+\varphi(-t_2)).$$
Ainsi, pour tout $t_1\in\R$,  on a 
$(A_{t_1})^{[1]}(0^{\pm})=0 $ et 
$(A_{t_1})^{[0]}(0^{\pm})=A(t_1,0)$ et de  plus, 
si $ t_1>0$ alors $A(t_1,0)=\chi(t_1)\varphi(t_1)=\chi(t_1) $ et  si $ t_1<0$ alors $A(t_1,0)=\chi(-|t_1|)\varphi(-t_1)=\chi(t_1).$

On obtient donc  $$I_\m(L_1+L_2, \Psi,A)=4
\int_\R\chi(t_1)((\Psi_{t_1})^{[1]}(0^{+})-(\Psi_{t_1})^{[1]}(0^{-}))dt_1=0.$$
Cette  \' egalit\' e est  vraie pour toute fonction
$\chi\in\D(\R^*)$  ce qui donne  la
premi\`ere conclusion
$$(\Psi_{t_1})^{[1]}(0^{+})=(\Psi_{t_1})^{[1]}(0^{-}),\forall t_1\in\R^*.$$

Maintenant, pour $(t_1,t_2)$ dans $(\R^*)^2,$ on pose
$C(t_1,t_2)=(\log|t_1|+\log|t_2|)A(t_1,t_2)\in \mathcal{H}_{\log}^2$ .
Pour $ t_1\in\R^*$,  on a $(C_{t_1})^{[1]}(0^{\pm})=A(t_1,0)=\chi(t_1)$ et
$(C_{t_1})^{[0]}(0^{\pm})=\log|t_1|\chi(t_1).$ De plus, comme $(\Psi_{t_1})^{[1]}(0^{+})=(\Psi_{t_1})^{[1]}(0^{-})$,
on en d\' eduit 
$$I_\m(L_1+L_2, \Psi,C)=4\int_\R\left(\chi(t_1)((\Psi_{t_1})^{[0]}(0^{+})-(\Psi_{t_1})^{[0]}(0^{-}))\right)dt_1.$$
Nous concluons comme pr\' ec\' edemment que
$$(\Psi_{t_1})^{[0]}(0^{+})-(\Psi_{t_1})^{[0]}(0^{-})=0,\forall t_1\in\R^*.$$
Nous obtenons ainsi l'implication $(1)\Rightarrow (2)$.\me

On suppose maintenant  que pour $j\in\{0,1\}$ et $t_1\in\R^*$ alors
$
(\Psi_{t_1})^{[j]}(0^+)=(\Psi_{t_1})^{[j]}(0^-).$
Par la remarque \ref{K}, pour $F\in\mathcal{H}^2_{\log},$ on a alors
imm\' ediatement
$K_2(\Psi,F)(t_1,0^+)=K_2(\Psi,F)(t_1,0^-)$ pour tout $ t_1\in\R^*$
ce qui donne l'assertion $(3)$.\me

L'implication $(3)\Rightarrow (1)$ est imm\' ediate par l'expression
de $I_\m(L_1+L_2,\Psi,F)$ donn\'ee au lemme \ref{ippqm} et la
surjectivit\'e l'application $f\mapsto \Mfm$  de $\D(H\cdot\p\mq)$
dans $\mathcal{H}_{\log}^2$.  \hskip0,5cm \carre

\begin{Def} \begin{enumerate}
\item Pour $f$ et $g$ deux fonctions d'une variable $x$ r\' eelle ou complexe, on note
$$S^+(f,g)(x_1,x_2)=f(x_1)g(x_2)+f(x_2)g(x_1)$$
et $$[f,g](x_1,x_2)=f(x_1)g(x_2)-f(x_2)g(x_1).$$
\item Si $(f_j)_{j\in J}$ est une famille finie d'un espace vectoriel, on note
$Vect\langle f_j; j\in J\rangle$ le sous-espace vectoriel engendr\'
e par la famille
 $(f_j)_{j\in J}$.
\end{enumerate}
\end{Def}

\begin{cor}\label{expsolm}
La fonction  $\Psi$ satisfait l'une des conditions \' equivalentes
de la proposition \ref{recolm} si et seulement si  
$$\Psi\in
Vect\langle S^+(A,B)(t_1,t_2), signe(t_1-t_2)[A,B](t_1, t_2);(A,B)\in
\mathcal{S}ol(L,\la_1)\times\mathcal{S}ol(L,\la_2)\rangle.$$

\end{cor}

\begin{dem} Par le lemme \ref{solm},  la fonction $\Psi$ s\' ecrit,
sous les formes suivantes:

\no Pour
$i\in\{1,2\}$, il existe  $(f_i^+,g_i^+)\in
\mathcal{S}ol(L,\la_2)\times\mathcal{S}ol(L,\la_1)$ telles que, pour $t_1>t_2>0,$  
$$\Psi(t_1,t_2)=f_1^+(t_1)\Phi_{\la_1}(t_2)+f_2^+(t_1)W^r_{\la_1}(t_2)+
g_1^+(t_1)\Phi_{\la_2}(t_2)+g_2^+(t_1)W^r_{\la_2}(t_2),$$
Pour
$i\in\{1,2\}$, Il existe $(f_i^-,g_i^-)\in
\mathcal{S}ol(L,\la_2)\times\mathcal{S}ol(L,\la_1)$ telles que, pour $t_1>0>t_2,$ 
$$\Psi(t_1,t_2)=f_1^-(t_1)\Phi_{\la_1}(t_2)+f_2^-(t_1)W^r_{\la_1}(t_2)+
g_1^-(t_1)\Phi_{\la_2}(t_2)+g_2^-(t_1)W^r_{\la_2}(t_2),$$ 
La condition $(\Psi_{t_1})^{[1]}(0^+)=(\Psi_{t_1})^{[1]}(0^-)$
implique que
$f_2^+(t_1)+g_2^+(t_1)=f_2^-(t_1)+g_2^-(t_1)$ pour tout $ t_1\in\R^*_+.$
Comme $\mathcal{S}ol(L,\la_2)\cap\mathcal{S}ol(L,\la_1)=\{0\},$
alors, pour tout $t_1\in\R_+^*$, on obtient
$f_2^+(t_1)=f_2^-(t_1)$  et 
$g_2^+(t_1)=g_2^-(t_1).$ On pose alors
$f_2(t_1)=f_2^+(t_1)=f_2^-(t_1)$ et
$g_2(t_1)=g_2^+(t_1)=g_2^-(t_1)$.

De m\^eme la condition
$(\Psi_{t_1})^{[0]}(0^+)=(\Psi_{t_1})^{[0]}(0^-)$ implique que, pour tout $t_1\in\R^*_+$ alors 
$f_1(t_1):=f_1^+(t_1)=f_1^-(t_1)$ et 
$g_1(t_1):=g_1^+(t_1)=g_1^-(t_1).$
 Ainsi   sur
$\{(t_1,t_2),\;t_1>t_2,\; t_1>0\}$, nous avons
$$\Psi(t_1,t_2)=f_1(t_1)\Phi_{\la_1}(t_2)+f_2(t_1)W^r_{\la_1}(t_2)+
g_1(t_1)\Phi_{\la_2}(t_2)+g_2(t_1)W^r_{\la_2}(t_2).$$ 
En raisonnant
de la m\^eme mani\`ere  sur l'ensemble
$\{(t_1,t_2),\;t_1>t_2,\; t_2<0\}$ en utilisant les relations
$(\Psi^{t_2})^{[j]}(0^+)=(\Psi^{t_2})^{[j]}(0^-),\;j\in\{0,1\}$, nous obtenons, pour $\;t_1>t_2$ et $t_2<0$
 la m\^eme expression pour $\Psi$, c'est-\`a-dire 
$\Psi(t_1,t_2)=f_1(t_1)\Phi_{\la_1}(t_2)+f_2(t_1)W^r_{\la_1}(t_2)+
g_1(t_1)\Phi_{\la_2}(t_2)+g_2(t_1)W^r_{\la_2}(t_2).$

Comme $\Psi$ est sym\' etrique par rapport aux deux variables, on
peut donc l'\' ecrire comme combinaison lin\' eaire des fonctions
$\mathbf{1}_{t_1\leq t_2} A(t_1)B(t_2) +\mathbf{1}_{t_1>t_2}
A(t_2)B(t_1)$ avec  $(A,B)$ ou $(B,A)$ dans $
\mathcal{S}ol(L,\la_1)\times\mathcal{S}ol(L,\la_2)$. On en d\' eduit
ais\' ement le
corollaire. \end{dem}
\begin{prop}\label{exImmoins} Soit  $\Psi\in Vect\langle S^+(A,B)(t_1,t_2),
signe(t_1-t_2)[A,B](t_1, t_2); (A,B)\in
\mathcal{S}ol(L,\la_1)\times\mathcal{S}ol(L,\la_2)\rangle$.

Alors, pour tout $F\in\mathcal{H}^2_{\log},$ on a
 la relation

$$I_\m((L_1-L_2)^2,\Psi,F)=0.$$

\end{prop}

\begin{dem} Nous allons montrer que les int\' egrales
$I_\m(L_1-L_2,(L_1-L_2)\Psi,F)$ et $I_\m(L_1-L_2,\Psi,(L_1-L_2)F)$
existent et sont nulles. Ceci donnera le r\' esultat voulu puisque,
dans ce cas, on aura
$$I_\m((L_1-L_2)^2,\Psi,F)=I_\m(L_1-L_2,(L_1-L_2)\Psi,F)+I_\m(L_1-L_2,\Psi,(L_1-L_2)F)=0.$$

On proc\`ede comme dans la preuve du lemme \ref{ippqm}.
 Par hypoth\`ese sur $\Psi$ et d\' efinition  de
$\mathcal{H}^2_{\log},$ pour $g=\Psi$ ou $F$, il existe des fonctions
$(g_{i,j})_{0\leq i,j\leq 1}$  de classe ${\mathcal C}^\infty$ sur
$\R^2$ telles que, pour $(t_1,t_2)\in(\R^*)^2-diag$, l'on ait
$$g(t_1,t_2)= \sum_{0\leq
i,j\leq1}g_{i,j}(t_1,t_2)(\log|t_1|)^{i}(\log|t_2|)^j .$$

Comme, pour  toute  fonction $f$   de classe $\mathcal{C}^2$ sur un
ouvert de $\R^*$, pour tout $t\in\R^*$, on a $\dis
L(\log|t|f(t))=8f'(t)+4\log|t|(f'(t)+tf''(t)),$ les fonctions $
L_k\Psi$ et $L_k F$, pour $k=1,2$ ont une expression du m\^eme type.
Ainsi, si $g$ d\' esigne l'une des fonctions $\Psi$, $F$,
$(L_1-L_2)\Psi$ ou $(L_1-L_2)F$,  les fonctions
$(t_1,t_2)\mapsto(g_{t_1})^{[j]}(t_2)$  sont continues en tout point
$(t,0)$ pour $t\in
\R^*$. 

Comme  $F$ est nulle en dehors d'un compact de $\R^2-diag$, le m\^eme raisonnement que dans la  
preuve du lemme \ref{ippqm} donne que, pour
$(g,h)=\big( (L_1-L_2)\Psi,F\big)$ ou $\big( \Psi,(L_1-L_2)F\big)$,
la fonction   $\dis g (L_kh)-(L_k g) h=4\frac{\partial} {\partial
t_k}K_k(g,h)$ est   int\'egrable sur $\R^2$  et on a la relation

$$I_\m(L_1-L_2, g,h)=4\int_{-\infty}^0[K_1(g,h)(0^-,t)-K_1(g,h)(0^+,t)]dt$$
$$\qq-4\int_0^{+\infty} [K_2(g,h)(t,0^{-})-K_2(g,h)(t,0^{+})] dt.$$

Comme $\Psi$ et $F$ sont sym\' etriques, l'une des deux fonctions
$g$ ou $h$ est sym\' etrique et l'autre est antisym\' etrique .
Ainsi, pour $(t_1,t_2)\in(\R^*)^2-diag$, on a
$$K_1(g,h)(t_1,t_2)=-K_2(g,h)(t_2,t_1).$$
et on obtient
$$I_\m(L_1-L_2,g,h)=4\int_\R[K_2(g,h)(t,0^{+})-K_2(g,h)(t,0^{-})] dt.$$

La remarque \ref{K} et les propri\' et\' es de $g$ et $h$  d\'
ecrites pr\' ec\' edemment donnent
$K_2(g,h)(t,0^{+})=K_2(g,h)(t,0^{-})$. Ainsi, on obtient
$$I_\m(L_1-L_2,\Psi,(L_1-L_2)F)=I_\m(L_1-L_2,(L_1-L_2)\Psi)=0.$$\end{dem}

 \subsubsection{ Recollement sur $H\cdot\ztp$}

Nous \'etudions dans ce paragraphe les conditions n\'ecessaires et suffisantes pour que les fonctions  $\Psi_\m$ et $\Psi_2$ v\'erifient les relations $({\mathcal Rec})$ pour tout $f\in \mathcal{D}(H\cdot\ztp)$.\me

Les \'el\'ements r\'eguliers de $H\cdot\ztp$ sont conjugu\'es par $H$ \`a un \'el\'ement de $\app$ ou de $\a_2$. La $H$-orbite d'un \'el\'ement semi-simple de $\app$ est caract\'eris\'ee par deux r\'eels positifs $t_1$ et $t_2$, param\'etrisation utilis\'ee dans le paragraphe pr\'ec\'edent et celle d'un \'el\'ement semi-simple de $\a_2$ par deux nombres complexes conjugu\'es $(\tau+\im\te)^2$ et $ (\tau-\im\te)^2.$

Ici, nous utiliserons la param\'etrisation suivante des \'el\'ements
de $\app$ plus coh\'erente avec celle de $\a_2$. Pour
$(\tau,\te)\in\R^2$, on note  $X_{\tau,\te}^r=\left(
    \begin{array}{c|c}
    0&\begin{array}{cc} \tau+\te&0\\ 0&\tau-\te\end{array}\\ \hline
    \begin{array}{cc} \tau+\te&0\\ 0&\tau-\te\end{array}&0
    \end{array}
    \right),$
de telle sorte que $\app=\{X_{\tau,\te}^r,\;(\tau,\te)\in\R^2\}.$
Avec ces notations, la $H$-orbite de $X_{\tau,\te}^r$ est donc
l'orbite caract\'eris\'ee par les deux r\'eels positifs
$(\tau+\te)^2$ et $(\tau-\te)^2$ et on a
$|\de|(X_{\tau,\te}^r)|=4|\tau\te|.$\\

Pour  $f\in  \mathcal{D}(H\cdot\ztp)$ et $\Psi_\m$ fix\'ee comme
pr\'ec\'edemment, on d\'efinit les fonctions $(\Mfm)_r$ et $
(\Psi_\m)_r$ sur $\{(\tau,\te)\in (\R^*)^2; \tau^2\neq\te^2\}$ par

    $$\begin{array}{l}(\Mfm)_r(\tau,\te)=\mathcal{M}_H(f)(X_{\tau,\te}^r)=\Mfm((\tau+\te)^2,(\tau-\te)^2)\\
    (\Psi_\m)_r(\tau,\te)=\tilde{\Psi}(X_{\tau,\te}^r)=\Psi_\m((\tau+\te)^2,(\tau-\te)^2)\end{array}.$$

Comme $\tilde{\Psi}$ et ${\mathcal M}_H(f)$ sont invariantes par le
groupe de Weyl $W_H(\app)$, les fonctions $(\Psi_\m)_r$ et
$(\Mfm)_r$ sont sym\' etriques et  paires en chaque variable.

Pour   $\dis \Phi=\frac{\tilde{\Psi}}{|\de|}$, la formule
d'int\'egration de Weyl ( lemme \ref{intweyl}) donne alors

$$\int_{H\cdot \app}  \Phi (X) f(X)
 =\int_{\R^2}  \Phi (X_{\tau, \te}^r)
    \mathcal{M}_H(f)(X_{\tau,\te}^r) |4\tau \te| |\tau^2-\te^2| d\tau
    d\te$$
    $$=\int_{\R^2}  (\Psi_\m)_r (\tau,\te)
    (\Mfm)_r(\tau,\te)|\tau^2-\te^2| d\tau
    d\te=\int_{t_1>t_2>0} \Psi_\m (t_1, t_2)\Mfm (t_1,t_2) d t_1 dt_2 .$$

      \begin{Def}
    Soient $f$ et $g$ deux fonctions des variables $(\tau,\te)\in
    \{(x,y)\in(\R^*)^2,\;x^2\ne y^2\}$. Soit $D$ un op\' erateur diff\'
    erentiel sur $\R^2$ tel que $Df$ et $Dg$ existent. On notera,
    lorsque l'int\' egrale consid\' er\' ee converge

     $$I_\m^r(D,f,g)=\int_{\R^2} |\tau^2-\te^2|\big( f\; (Dg)-(Df)\;
    g)(\tau,\te) d\tau d\te.$$
        \end{Def}

Sur chaque sous-espace de Cartan $\a$, les composantes radiales des
op\'erateurs $\partial(Q)$ et $\partial(S_0)$ sont donn\'ees en
terme des op\'erateurs $L_{\ga}=\frac{1}{4\ga}\partial
    h_{\ga}\big( \ga
    \partial h_{\ga}\big)$  (proposition \ref{prad}) o\`u $\ga$ est une racine positive de multiplicit\'e $1$ de $\a$.

    Sur $\app$, les racines positives de multiplicit\'e $1$ sont donn\'ees par $\al_1(X_{\tau,\te}^r)=2(\tau+\te)$ et $\al_2(X_{\tau,\te}^r)=2(\tau-\te)$.
Par ailleurs, les racines positives de multiplicit\'e $1$ de $\a_2$
sont donn\'ees par $\al(X_{\tau,\te})=2(\tau+\im\te)$
 et sa conjugu\'ee ${\overline \al}$.  Ainsi, les op\'erateurs $L_\ga$, pour $\ga(X)=2(\tau+\ep\te)$ avec $\ep=\pm 1$ si $X=X_{\tau,\te}^r\in\app$ et $\ep=\pm\im$ si $X=X_{\tau,\te}\in\a_2$, s'expriment de la mani\`ere suivante
 $$L_{\ga}=\frac{1}{4(\tau+\ep\te)}(\frac{\partial}{\partial\tau}+\overline{\ep}\frac{\partial}{\partial\te})\big((\tau+\ep\te)(\frac{\partial}{\partial\tau}+\overline{\ep}\frac{\partial}{\partial\te})\big).$$

      Gr\^ace aux  lemme  \ref{intrad2} et  corollaire \ref{intrad} on obtient donc:
     \begin{lem}\label{chgvar} Pour $f\in\D(H\cdot\ztp)$, on a

    $$\begin{array}{ll} & I_\m(L_1+L_2,\Psi_\m,\Mfm)=
    I_\m^r(L_{\al_1}+L_{\al_2},(\Psi_m)_r,(\Mfm)_r)\\
    \textrm{et } & I_\m\big((L_1-L_2)^2,\Psi_\m,\Mfm)=
    I_\m^r(\big(L_{\al_1}-L_{\al_2})^2,(\Psi_m)_r,(\Mfm)_r).\end{array}$$

    \end{lem}

       Ainsi les relations de recollement $({\mathcal Rec})$ s'\'ecrivent pour $f\in{\mathcal D}(H\cdot\ztp)$,

    $$\left\{\begin{array}{l} I_\m^r(L_{\al_1}+L_{\al_2},(\Psi_m)_r,(\Mfm)_r)+2 I_2( L_\al+L_{\overline{\al}},\Psi_2,\Mf2)=0\\
     I_\m^r\big((L_{\al_1}-L_{\al_2})^2,(\Psi_m)_r,(\Mfm)_r\big)+2 I_2\big(( L_\al-L_{\overline{\al}})^2,\Psi_2,\Mf2)=0\end{array}\right.$$\me

Pour \'etudier ces relations, on introduit les notations suivantes.
Pour $g$ et $h$ deux fonctions de classe ${\mathcal C}^2$  sur un
ouvert $U$ de $\R^2$, on d\'efinit les op\'erateurs $R_\tau$ et
$R_\te$ par
$$R_{\tau}(g,h)(\tau,\te)=\left(g\; \frac{\partial h}{\partial
\tau}-\frac{\partial g}{\partial \tau}\;h\right)(\tau,\te)\mbox{ et
}R_{\te}(g,h)(\tau,\te)=\left(g\; \frac{\partial h}{\partial
\te}-\frac{\partial g}{\partial \te}\;h\right)(\tau,\te).$$

Soit   $\ga(X)=2(\tau+\ep\te)$ avec $\ep=\pm 1$ si
$X=X_{\tau,\te}^r\in\app$ et $\ep=\pm\im$ si
$X=X_{\tau,\te}\in\a_2$. On a alors $\dis\partial h_\ga
=(\frac{\partial}{\partial \tau}+\overline{ \ep}
\frac{\partial}{\partial \te})$ et  donc $\dis g (\partial h_\ga
h)-(\partial h_\ga  g)h=R_\tau(g,h) +\overline{\ep} R_\te(g,h)$.

Ainsi, pour $(\tau,\te)\in U$  avec $\tau^2\neq\ep^2\te^2$, on
obtient
  $$|\tau^2-\ep^2\te^2|\big[g( L_\ga h)- (L_\ga g)h\big](\tau,\te)=\frac{1}{4}(\frac{\partial}{\partial \tau}+\overline{\ep} \frac{\partial}{\partial \te})\Big( |\tau^2-\ep^2\te^2|\big[R_\tau(g,h) +\overline{\ep}  R_\te(g,h)\big]\Big)(\tau,\te),$$
et donc, on a les relations suivantes
$$\begin{array}{lc} &  |\tau^2-\te^2|\left(g(L_{\al_{1}}+L_{\al_{2}})(h)-h(L_{\al_{1}}+L_{\al_{2}})(g)\right)(\tau,\te)\\
& =\frac{1}{2}\big[\dfrac{\partial}{\partial\tau}\left(|\tau^2-\te^2|R_\tau(g,h)\right)(\tau,\te)
+\dfrac{\partial}{\partial\te}\left(|\tau^2-\te^2|R_\te(g,h)\right)(\tau,\te)\big],\\
& \\
\textrm {et }& |\tau^2+\te^2|\left(g(L_{\al}+L_{\overline{\al}})(h)-h(L_{\al}+L_{\overline{\al}})(g)\right)(\tau,\te)\\
& =\frac{1}{2}\big[\dfrac{\partial}{\partial\tau}\left((\tau^2+\te^2)R_\tau(g,h)\right)(\tau,\te)
-\dfrac{\partial}{\partial\te}\left((\tau^2+\te^2)R_\te(g,h)\right)(\tau,\te)\big].\end{array}$$

Nous rappelons maintenant les propri\' et\' es des fonctions  $(\Psi_\m)_r$, $\Psi_2$, $(\Mfm)_r$ et $\Mf2$ que nous utiliserons.

  Par  le th\'eor\`eme
    \ref{ioa2},  pour tout $f\in\D(H\cdot\ztp)$, il existe une unique fonction $G_f=a(\tau, w)+|w|^{1/2}b(\tau,w)\in\mathcal{H}_Y^{pair}$ o\`u $a$ et $b$ sont paires
    par rapport \`a la premi\`ere variable et appartiennent \`a $\D(\R^*\times\R \cap\{(\tau,w); \tau^2> w\})$, telle que, pour  $\tau\te\ne0$, on ait\me

$\left\{\begin{array}{l}(\Mfm)_r(\tau,\te)=\mathcal{M}_H(f)(X_{\tau,\te}^r)=G_f(\tau,\te^2)=a(\tau,\te^2)
\quad \textrm{ pour } \tau^2-\te^2> 0\\
(\Mfm)_r(\tau,\te)=\mathcal{M}_H(f)(X_{\te,\tau}^r)=G_f(\te,\tau^2)=a(\te,\tau^2)\quad 
\textrm{ pour } \tau^2-\te^2< 0\\
\Mf2(\tau,\te)=\mathcal{M}_H(f)(X_{\tau,\te})=G_f(\tau,-\te^2)=a(\tau,-\te^2)+|\te|b(\tau,-\te^2) .
   \end{array}\right.$\me

  De plus, l'application  $ f\mapsto G_f$ est surjective de $\D(H\cdot\ztp)$ dans $\mathcal{H}_Y^{pair}$.

    \begin{rem}\label{supportMf}
   La fonction    $(\Mfm)_r$  est sym\' etrique et paire en chaque variable. Elle se prolonge en une  fonction de   classe $\mathcal{C}^\infty$ \`a support compact contenu dans $\{(\tau,\te)\in\R^2; \tau^2\neq \te^2\}$. La fonction $\Mf2$ se prolonge  en une fonction continue \`a support compact contenu dans $\R^*\times\R$ et pour tout $D\in\C[\dfrac{\partial}{\partial \tau}, \dfrac{\partial}{\partial \te}]$ et $\tau\neq 0$, les limites $(D\Mf2)(\tau,0^\pm)$ existent.
   \end{rem}
 La fonction  $(\Psi_m)_r(\tau,\te)$ est sym\' etrique et paire en chaque variable.  Gr\^ace au lemme \ref{solm}, pour $\tau>\te>0$, elle est somme de fonctions de la forme
    $[A((\tau+\te)^2)+\log(\tau+\te)^2B((\tau+\te)^2)][C((\tau-\te)^2)+\log(\tau-\te)^2D((\tau-\te)^2)],$
    o\`u $A,\;B,\;C$ et $D$ sont analytiques sur $\R$.

  La fonction  $\Psi_2$ est paire en chaque variable et la  proposition \ref{sold}  permet d'\'ecrire $\Psi_2$  sur $(\R_+^*)^2$ comme somme de fonctions de  la forme
    $[A_c((\tau+\im\te)^2)+\log(\tau+\im\te)^2B_c((\tau+\im\te)^2)][C_c((\tau-\im\te)^2)+\log(\tau-\im\te)^2D_c((\tau-\im\te)^2)],$
    o\`u $A_c,\;B_c,\;C_c$ et $D_c$ sont analytiques sur $\C$.
    
 \begin{rem}\label{regularitepsi} En particulier, pour tout $D\in\C[\dfrac{\partial}{\partial \tau}, \dfrac{\partial}{\partial \te}]$ et $\tau\neq 0$, les limites   $D(\Psi_\m)_r(\tau,0^+)$, $D(\Psi_\m)_r(0^+,\tau)$ et $(D \Psi_2)(\tau,0^\pm)$ existent.

   \end{rem}

   \begin{lem}\label{ippmr} Soit $f\in\D( H\cdot\ztp)$.
 On a les relations suivantes
 $$ I_\m^r(L_{\al_1}+L_{\al_2},(\Psi_\m)_r,(\Mfm)_r)=-4\int_0^{+\infty} \tau^2
R_\te((\Psi_\m)_r,(\Mfm)_r)(\tau,0^+)
    d\tau$$
$$ I_2(L_\al+L_{\overline{\al}},\Psi_2,\Mf2)=2\int_0^{+\infty} \tau^2
R_\te(\Psi_2,\Mf2)(\tau,0^+)
    d\tau.$$

  \end{lem}

\begin{dem}   Par les remarques \ref{supportMf} et \ref{regularitepsi}, pour $R=R_\tau$ ou $R_\te$ et $D\in\C[\dfrac{\partial}{\partial \tau}, \dfrac{\partial}{\partial \te}]$, la fonction  $D\big[|\tau^2-\te^2|R\big((\Psi_\m)_r, (\Mfm)_r\big)\big]$ se prolonge en une fonction de classe ${\mathcal C}^\infty$  sur $\R_+^2$ \`a support compact contenu dans $\D(\R^2_+-diag)$. Par parit\' e en $\tau$ et $\te$ de $(\Psi_\m)_r$ et $(\Mfm)_r$, cette fonction est int\' egrable sur $\R^2$.  Ainsi, on a

$$I_\m^r(L_{\al_1}+L_{\al_2},(\Psi_\m)_r,(\Mfm)_r)=2\int_{\begin{array}{c}\tau>0\\ \te>0\end{array}}
\left(\dfrac{\partial}{\partial\tau}\left(|\tau^2-\te^2|R_\tau((\Psi_\m)_r,(\Mfm)_r)\right)(\tau,\te)\right.
$$$$\left.+\dfrac{\partial}{\partial\te}\left(|\tau^2-\te^2|R_\te((\Psi_\m)_r,(\Mfm)_r)\right)(\tau,\te)\right)
d\tau d\te$$

$$=-2\int_0^{+\infty} \te^2 R_\tau((\Psi_\m)_r,(\Mfm)_r)(0^+,\te) d\te-2\int_0^{+\infty} \tau^2 R_\te((\Psi_\m)_r,(\Mfm)_r)(\tau,0^+) d\tau.$$

Comme les fonctions $(\Psi_\m)_r$ et $(\Mfm)_r$ sont sym\' etriques
en $(\tau,\te)$, on obtient
$$I_\m^r(L_{\al_1}+L_{\al_2},(\Psi_\m)_r,(\Mfm)_r)=-4\int_0^{+\infty} \tau^2 R_\te((\Psi_\m)_r,(\Mfm)_r)(\tau,0^+) d\tau,$$
ce qui donne la premi\`ere \' egalit\' e.

De m\^eme, par les remarques \ref{supportMf} et \ref{regularitepsi},
les fonctions $\dfrac{\partial}{\partial\tau}
[(\tau^2+\te^2)R_\tau(\Psi_2,\Mf2)]$ et
$\dfrac{\partial}{\partial\te}[ (\tau^2+\te^2)R_\te(\Psi_2,\Mf2)]$
sont paires en chaque variable et int\' egrables sur $\R^2$ et on a
$$I_2(L_\al+L_{\overline{\al}},\Psi_2,\Mf2)=2\int_{\begin{array}{c}\tau>0\\ \te>0\end{array}}
\dfrac{\partial}{\partial\tau}\left((\tau^2+\te^2)R_\tau(\Psi_2,\Mf2)\right)(\tau,\te)d\tau
d\te$$
$$-2\int_{\begin{array}{c}\tau>0\\ \te>0\end{array}}\dfrac{\partial}{\partial\te}\left((\tau^2+\te^2)R_\te(\Psi_2,\Mf2)\right)(\tau,\te)
d\tau d\te.$$

 On obtient la deuxi\`eme relation  car  $\Mf2$ est nulle au voisinage des points $(0,\te)\in\{0\}\times \R^*$. \end{dem}

   \begin{prop}\label{recolztp} Les trois assertions suivantes sont \' equivalentes:
     \begin{enumerate}
     \item Pour tout $f\in\D(H\cdot\ztp)$, on a

    $$I_\m(L_1+L_2,\Psi_\m,\Mfm)+2I_2(L_\al+L_{\overline{\al}},\Psi_2,\Mf2)=0.$$

    \item Pour tout $\tau>0$, on a
        $$\left\{\begin{array}{c}\frac{\partial}{\partial\te}(\Psi_\m)_r(\tau,0^+)-
    \frac{\partial}{\partial\te}\Psi_2(\tau,0^+)=0\\
    \Psi_2(\tau,0^+)=0\;\;\qq\qq\end{array}\right..$$
    \item Pour tout $f\in\D(H\cdot\ztp)$ et pour tout $\tau>0$, on a
     $$R_\te((\Psi_\m)_r,(\Mfm)_r)(\tau,0^+)-R_\te(\Psi_2,\Mf2)(\tau,0^+)=0.$$
    \end{enumerate}
    \end{prop}

    \begin{dem}  On suppose l'assertion {\it 1.} v\' erifi\' ee. Par les
    lemmes \ref{chgvar} et \ref{ippmr}, pour tout $f\in\D(H\cdot\ztp)$,
    on a donc
    $$-\int_0^{+\infty}\tau^2 R_\te((\Psi_\m)_r,(\Mfm)_r)(\tau, 0^+) d\tau+
    \int_0^{+\infty} \tau^2 R_\te(\Psi_2,\Mf2)(\tau,0^+) d\tau=0.$$
        Par surjectivit\' e de l'application $f\mapsto G_f$ de
    $\D(H\cdot\ztp)$ dans $\mathcal{H}_Y^{pair}$ (th\' eor\`eme
     \ref{ioa2}), l'assertion  {\it 1.} est donc \' equivalente \` a

    $$\int_0^{+\infty}\tau^2\big[ a(\tau,0)\big( \frac{\partial}{\partial\te}(\Psi_\m)_r
    -\frac{\partial
    \Psi_2}{\partial\theta}\big)(\tau,0^+)+b(\tau,0)\Psi_2(\tau,
    0^+)\big] d\tau=0$$ ceci pour tout $a$ et $b$ dans
    $\D((\R^*\times\R)\cap\{(\tau,w); \tau^2>w\})$ paires par rapport
    \`a la premi\`ere
    variable.

     Soit $\chi\in\D(\R^*_+)$
    quelconque. Il existe $\varepsilon\in]0,1[$ tel que $supp(\chi)$
    soit inclus dans $[\varepsilon,\varepsilon^{-1}]$. Soit
    $\phi\in\D(\R)$ telle que $\phi\equiv1$ sur
    $[-\varepsilon^2/2;\varepsilon^2/2]$ et $\phi\equiv0$ sur
    $\R-[-\varepsilon^2;\varepsilon^2]$.

    En prenant la fonction $a(t,w)$ nulle et $b(t,w)$ d\' efinie par
    $b(t,w)=\frac{\chi(t)+\chi(-t)}{t^2}(1-\phi(t^2-w))\phi(w),$ on
    obtient $\displaystyle b(\tau,0)=\frac{1}{\tau^2}\chi(\tau)$ pour
    $\tau>0$ et donc
    $$\int_0^{+\infty}\chi(\tau)\Psi_2(\tau,0^+)d\tau=0.$$
    Cette derni\`ere \'egalit\' e \' etant vraie pour toute fonction
    $\chi$ de $\D(\R^*_+)$, alors
    $$\Psi_2(\tau,0^+)=0,\qq\forall\tau\in\R^*_+.$$

    Maintenant, en prenant la fonction
    $a(t,w)=\frac{\chi(t)+\chi(-t)}{t^2}(1-\phi(t^2-w))\phi(w)$ et la
    fonction $b(t,w)$ nulle, nous obtenons alors

    $$\int_0^{+\infty}\chi(\tau)\left(\frac{\partial
    (\Psi_\m)_r}{\partial\theta}(\tau,0^+)-\frac{\partial
    \Psi_2}{\partial\theta}(\tau,0^+)\right)d\tau=0.$$

    Cette derni\`ere \' egalit\' e est vraie pour toute fonction $\chi$
    de $\D(\R^*_+)$, ainsi on a
    $$\frac{\partial
    (\Psi_\m)_r}{\partial\theta}(\tau,0^+)-\frac{\partial
    \Psi_2}{\partial\theta}(\tau,0^+)=0.$$
    Ceci donne l'implication ${\it 1.}\Rightarrow {\it 2.}$\me

    Maintenant on suppose  que
    $\dis \Psi_2(\tau,0^+)=0 \mbox{ et } \frac{\partial
    (\Psi_\m)_r}{\partial\theta}(\tau,0^+)-\frac{\partial
    \Psi_2}{\partial\theta}(\tau,0^+)=0.$
  Soit $f\in\D(H\cdot\ztp)$. Nous avons les propri\' et\' es
    $$ (\Mfm)_r(\tau,0^+)= \Mf2(\tau,0^+),\quad  \frac{\partial
    }{\partial\theta}(\Mfm)_r(\tau,0^+)=0,\textrm{ et } \frac{\partial
    }{\partial\theta}\Mf2(\tau,0^+)\textrm{  est finie.}$$

    Les hypoth\` eses sur $\Psi_2$ et $\Psi_\m$ donne alors
    $$R_\te(\Psi_2,\Mf2)(\tau,0^+)-R_\te((\Psi_\m)_r,(\Mfm)_r)(\tau,0^+)$$
    $$=-\frac{\partial
    \Psi_2}{\partial\theta}(\tau,0^+)\Mf2(\tau,0^+)+\frac{\partial
    (\Psi_\m)_r}{\partial\theta}(\tau,0^+)(\Mfm)_r(\tau,0^+)=0,$$

    d'o\` u l'implication ${\it 2.}\Rightarrow {\it 3.}$\me
    
    L'implication ${\it 3.}\Rightarrow {\it 1.}$ est imm\' ediate par
    les lemmes \ref{chgvar} et \ref{ippmr}.
   \end{dem}

    Nous d\' ecrivons maintenant l'espace des fonctions $\Psi_\m$ et
    $\Psi_2$
    satisfaisant les assertions de la proposition \ref{recolztp}.

     Pour $k\in\{1,2\}$, on note $$A_{1,\la_k}(z)=\Phi_{\la_k}(z)\qq\textrm{et}\qq
    A_{2,\la_k}(z)=W_{\la_k}(z)=w_{\la_k}(z) + \log(z) \Phi_{\la_k}(z),$$
     le syst\` eme fondamental de solutions de $L_cy=\la_k y$ sur $\C-\R_-$ (proposition
    \ref{fonsol}).
    
 La restriction de ces fonctions \` a $\R_+^*$ est un syst\` eme fondamental de
    solution de $Ly=\la_k y$ sur $\R^*_+$ c'est-\` a-dire, pour $t>0$,
    on a $W^r_{\la_k}(t)=W_{\la_k}(t)$. On omettra l'exposant $"r"$ dans
    la suite de ce paragraphe.

     \begin{lem}\label{indepdte} Sur $\R^*_+$, la famille de fonctions $\{
    (A_{i,\la_1}A_{j,\la_2})(t)\}_{1\leq i,j\leq 2}$ est libre.
     \end{lem}

    \begin{dem} Soit $\al_{i,j}\in \C$ tels que, pour tout $t>0$, l'on ait
    $\sum_{1\leq i,k\leq 2}\al_{i,j}A_{i,\la_1}(t)A_{j,\la_2}(t)=0$.
    Cette relation s'\' ecrit aussi
     $$\big[\al_{1,1}\Phi_{\la_1}\Phi_{\la_2}+\al_{1,2}\Phi_{\la_1}
    w_{\la_2}+\al_{2,1}w_{\la_1}\Phi_{\la_2}+\al_{2,2}
    w_{\la_1}w_{\la_2}\big](t)$$
     $$+ \log(t)\big[(\al_{1,2}+\al_{2,1})
    \Phi_{\la_1}\Phi_{\la_2} +\al_{2,2}(\Phi_{\la_1}
    w_{\la_2}+w_{\la_1}\Phi_{\la_2})\big](t)+\log(t)^2\al_{2,2}\Phi_{\la_1}\Phi_{\la_2}(t)=0.$$
     Les fonctions $\Phi_{\la_k}$ et $w_{\la_k}$ sont continues et non nulles en $0$. Ainsi, on
    obtient ais\' ement $$\al_{2,2}= \al_{1,2}+\al_{2,1}=0\quad \textrm{ et }\quad  \al_{1,1}\Phi_{\la_1}\Phi_{\la_2}+\al_{1,2}\big(\Phi_{\la_1}
    w_{\la_2}-w_{\la_1}\Phi_{\la_2} \big)=0.$$

    Comme  $w_{\la_1}(0)=w_{\la_2}(0)=2\gamma$ o\`u $\gamma$  est la
    constante d'Euler  et $\Phi_{\la_1}(0)=\Phi_{\la_2}(0)=1$, on
    obtient
    $\al_{1,1}=0.$
    
    Par ailleurs, comme $\displaystyle
    w_{\la_k}(t)=2\gamma+(1+2\gamma)\frac{\la_k t}{4}+o_{t\to0^+}(t)$ et
    $\displaystyle \Phi_{\la_k}(t)=1+\frac{\la_k t}{4}+o_{t\to0^+}(t)$,
    on a
    $\Phi_{\la_1}(t)w_{\la_2}(t)-\Phi_{\la_2}(t)w_{\la_1}(t)=\frac{t}{4}(\la_2-\la_1)+o_{t\to0^+}(t)$
    ce qui donne
    $\al_{1,2}=0.$\end{dem}

    \begin{cor}\label{solztp} Les fonctions  $\Psi_\m$ et $\Psi_2$ satisfont l'une des
    conditions \' equivalentes de la proposition \ref{recolztp} si et seulement si       il existe  $\Psi^+\in Vect\langle S^+(A,B); (A,B)\in
    \mathcal{S}ol(L,\la_1)\times \mathcal{S}ol(L,\la_2)\rangle$
 et $ \Psi_c\in Vect \langle [A,B]; (A,B)\in
    \mathcal{S}ol(L_c,\la_1)\times \mathcal{S}ol(L_c,\la_2)\rangle $
    telles que

     $$\begin{array}{l}\mbox{ pour } t_1>t_2>0 \textrm{ alors }
    \Psi_\m(t_1,t_2)=\Psi^+(t_1,t_2)+\im\Psi_c(t_1,t_2) \\
    \mbox{ pour }\tau>0\mbox{ et }\te>0 \textrm{ alors }
     \Psi_2(\tau,\te)=\Psi_c\big((\tau+\im\te)^2,(\tau-\im\te)^2\big) \end{array}$$

    \end{cor}

   \begin{dem} Par la proposition \ref{sold2}, il existe des  nombres
    complexes $\al_{i,j}$ et $\be_{i,j}$ tels que
    $\Psi_c(z_1,z_2)=\sum_{1\leq
    i,j\leq2}\al_{i,j}A_{i,\la_1}(z_1)A_{j,\la_2}(z_2)+\be_{i,j}A_{i,\la_1}(z_2)A_{j,\la_2}(z_1).$
     Ainsi, pour $\tau>0$, la condition  $\Psi_2(\tau,0)=0$ implique que,
    $$\sum_{1\leq
    i,j\leq2}(\al_{i,j}+\be_{i,j})A_{i,\la_1}(\tau^2)A_{j,\la_2}(\tau^2)=0.
    $$
    Par le lemme pr\' ec\' edent, on obtient donc $
    \al_{i,j}+\be_{i,j}=0$ pour $1\leq i,j\leq2$ et par suite
    $$\Psi_c(z_1,z_2)=\sum_{1\leq
    i,j\leq2}\al_{i,j}\big(A_{i,\la_1}(z_1)A_{j,\la_2}(z_2)-A_{i,\la_1}(z_2)A_{j,\la_2}(z_1)\big).$$

    Etudions \`a pr\' esent les cons\' equences de la condition
    $\frac{\partial
    }{\partial\theta}(\Psi_\m)_r(\tau,0^+)-\frac{\partial
    }{\partial\theta}\Psi_2(\tau,0^+)=0$.

   Par la proposition \ref{sold}, il existe des nombres complexes
    $a_{i,j}$ et $b_{i,j}$ tels que pour $\tau>\te\geq0$, on ait:
    
  \no$\dis (\Psi_\m)_r(\tau,\te)=\sum_{1\leq
    i,j\leq2}a_{i,j}A_{i,\la_1}((\tau+\te)^2)A_{j,\la_2}((\tau-\te)^2)+\sum_{1\leq
    i,j\leq2}b_{i,j}A_{i,\la_1}((\tau-\te)^2)A_{j,\la_2}((\tau+\te)^2).$

    On obtient donc

    $\dis \frac{\partial
    }{\partial\theta}(\Psi_\m)_r(\tau,0^+)=2\tau\sum_{1\leq
    i,j\leq2}(a_{i,j}-b_{i,j})
    (A_{i,\la_1}^{'}(\tau^2)A_{j,\la_2}(\tau^2)-A_{i,\la_1}(\tau^2)A_{j,\la_2}^{'}(\tau^2))$
    
    et
    $\dis \frac{\partial
    }{\partial\theta}\Psi_2(\tau,0^+)=4\im\tau\sum_{1\leq
    i,j\leq2}\al_{i,j}
    (A_{i,\la_1}^{'}(\tau^2)A_{j,\la_2}(\tau^2)-A_{i,\la_1}(\tau^2)A_{j,\la_2}^{'}(\tau^2)).$
    
    La condition $\frac{\partial
    }{\partial\theta}(\Psi_\m)_r(\tau,0^+)-\frac{\partial
    }{\partial\theta}\Psi_2(\tau,0^+)=0$ donne alors, pour $t>0$
    $$\sum_{1\leq
    i,j\leq2}(a_{i,j}-b_{i,j}-2\im \al_{i,j})
    (A_{i,\la_1}^{'}A_{j,\la_2}-A_{i,\la_1}A_{j,\la_2}^{'})(t)=0.$$

    En d\' erivant cette expression, on obtient
    $$\sum_{1\leq
    i,j\leq2}(a_{i,j}-b_{i,j}-2\im\al_{i,j})
    (A_{i,\la_1}^{''}A_{j,\la_2}-A_{i,\la_1}A_{j,\la_2}^{''})(t)=0.$$
    Ainsi, pour $L=4(t\dfrac{\partial^2}{\partial t^2}+\dfrac{\partial}{\partial t})$ et  pour tout $t>0$, nous avons
    $$0=\sum_{1\leq
    i,j\leq2}(a_{i,j}-b_{i,j}-2\im\al_{i,j})
    ([LA_{i,\la_1}]A_{j,\la_2}-A_{i,\la_1}[LA_{j,\la_2}])(t)$$
 $$
   = \sum_{1\leq
    i,j\leq2}(a_{i,j}-b_{i,j}-2\im\al_{i,j})(\la_1-\la_2)A_{i,\la_1}A_{j,\la_2}(t).$$

    Le lemme \ref{indepdte} donne alors
    $a_{i,j}-b_{i,j}-2\im\al_{i,j}=0 \mbox{ pour } 1\leq i,j\leq2.$ Ceci permet d'obtenir l'expression suivante de $\Psi_\m$ pour  $t_1>t_2>0$:

    $$\Psi_\m(t_1,t_2)=\sum_{1\leq
    i,j\leq2}a_{i,j}A_{i,\la_1}(t_1)A_{j,\la_2}(t_2)
    +b_{i,j}A_{i,\la_1}(t_2)A_{j,\la_2}(t_1)$$
   
    $$= \sum_{1\leq
    i,j\leq2}\frac{a_{i,j}+b_{i,j}}{2}\big(A_{i,\la_1}(t_1)A_{j,\la_2}(t_2)+
    A_{i,\la_1}(t_2)A_{j,\la_2}(t_1)\big)+\im\Psi_c(t_1,t_2).$$
    ce qui donne les expressions de $\Psi_2$ et $\Psi_\m$ voulues. La r\' eciproque est imm\' ediate par simple calcul. \end{dem}

    \begin{prop}\label{recolS0ztp}
    On suppose qu'il existe  $\Psi^+\in Vect\langle S^+(A,B); (A,B)\in
    \mathcal{S}ol(L,\la_1)\times \mathcal{S}ol(L,\la_2)\rangle$
 et $ \Psi_c\in Vect \langle [A,B]; (A,B)\in
    \mathcal{S}ol(L_c,\la_1)\times \mathcal{S}ol(L_c,\la_2)\rangle $
    telles que

     $$\begin{array}{l}\mbox{ pour } t_1>t_2>0 \textrm{ alors }
    \Psi_\m(t_1,t_2)=\Psi^+(t_1,t_2)+\im\Psi_c(t_1,t_2) \\
    \mbox{ pour }\tau>0\mbox{ et }\te>0 \textrm{ alors }
     \Psi_2(\tau,\te)=\Psi_c\big((\tau+\im\te)^2,(\tau-\im\te)^2\big) \end{array}$$

    Alors, pour tout $f\in\D(H\cdot\ztp)$, on a

    $$I_\m((L_1-L_2)^2,\Psi_\m,\Mfm)+2I_2((L_\al-L_{\overline{\al}})^2,\Psi_2,\Mf2)=0.$$
    \end{prop}
%%%%%%%%%%%%%%%%%%%%%%%%%%%%

\begin{dem} Soit $f\in\D(H\cdot\ztp)$. Il existe deux fonctions $a$
et $b$  de  $\D(\R^*\times\R \cap\{(\tau,w); \tau^2> w\})$ paires en
la premi\`ere variable telles que:

\no  pour $\tau>\te\geq 0$, alors $(\Mfm)_r(\tau,\te)=a(\tau,\te^2)$

\no et pour $\tau>0$ et $\te\geq 0$ alors $\Mf2(\tau,\te)=a(\tau,-\te^2)+\te b(\tau,-\te^2).$\\

 On proc\`ede comme dans la preuve du lemme  \ref{ippmr}.
On rappelle que $I_\m((L_1-L_2)^2,\Psi_\m,\Mfm)=
I_\m^r\big(((L_{\al_1}-L_{\al_2})^2,(\Psi_\m)_r,(\Mfm)_r\big)$.
\medskip

On consid\`ere tout d'abord  $D=L_{\al_1}-L_{\al_2}$ et $(g,h)=((\Psi_\m)_r,D(\Mfm)_r)$ ou $(D(\Psi_\m)_r),(\Mfm)_r)$. Par l'expression des $L_{\al_j}$ pour $j=1,2$, on a d'une part
$$L_{\al_1}-L_{\al_2}=\fr{\partial\tau}\fr{\partial\te}-\frac{\te}{\tau^2-\te^2}\fr{\partial\tau}
    +\frac{\tau}{\tau^2-\te^2}\fr{\partial\te},$$

\no et d'autre part, pour $A$ et $B$ deux fonctions de classe
$\mathcal{C}^2$  sur un ouvert $U$ de $\R^2$ et pour $(\tau,\te)\in
U$  avec $\tau^2\neq\te^2$, on a la relation suivante

 $$|\tau^2-\te^2|\left(A(L_{\al_{1}}-L_{\al_{2}})(B)-A(L_{\al_{1}}-L_{\al_{2}})(B)\right)(\tau,\te)$$
 $$=
\frac{1}{2}\big[\dfrac{\partial}{\partial\te}\left(|\tau^2-\te^2|R_\tau(A,B)\right)(\tau,\te)
+\dfrac{\partial}{\partial\tau}\left(|\tau^2-\te^2|R_\te(A,B)\right)(\tau,\te)\big].$$

Les remarques \ref{supportMf} et \ref{regularitepsi} assurent
l'int\' egrabilit\' e des deux fonctions
$\dfrac{\partial}{\partial\tau}\left(|\tau^2-\te^2|R_\te(g,h)\right)(\tau,\te)$
et
$\dfrac{\partial}{\partial\te}\left(|\tau^2-\te^2|R_\tau(g,h)\right)(\tau,\te)$
sur $\R_+^2$. Comme  $D(\Mfm)_r$ et $D(\Psi_\m)_r$ sont  impaires en
chaque variable, la fonction  $g(Dh)-(Dg)h$ est paire en chaque
variable. On obtient alors

$$I_\m^r(L_{\al_1}-L_{\al_2},g,h)=2\int_{\begin{array}{c}\tau>0\\ \te>0\end{array}}
\left(\dfrac{\partial}{\partial\tau}\left(|\tau^2-\te^2|R_\te(g,h)\right)(\tau,\te)\right.
$$$$\left.+\dfrac{\partial}{\partial\te}\left(|\tau^2-\te^2|R_\tau(g,h)\right)(\tau,\te)\right)
d\tau d\te.$$

puis, en utilisant la sym\' etrie des fonctions $g$ et $h$,
$$I_\m^r(L_{\al_1}-L_{\al_2},g,h)=-2\int_0^{+\infty} \te^2 R_\te(g,h)(0^+,\te) d\te-2\int_0^{+\infty} \tau^2 R_\tau(g,h)(\tau,0^+) d\te$$
$$=-4\int_0^{+\infty} \tau^2 R_\tau(g,h)(\tau,0^+) d\te$$

Pour $\tau>\te\geq 0$ on a $(\Mfm)_r(\tau,\te)=a(\tau,\te^2)$. Ainsi
la fonction  $h=(L_{\al_1}-L_{\al_2})(\Mfm)_r$ v\' erifie
$h(\tau,0^+)=\dfrac{\partial}{\partial \tau}h(\tau, 0^+)=0$ pour
tout $\tau>0$ et donc
$I_\m^r(L_{\al_1}-L_{\al_2},(\Psi_\m)_r,(L_{\al_1}-L_{\al_2})(\Mfm)_r)=0.$
On en d\' eduit l'\' egalit\' e suivante.
$$I_\m^r((L_{\al_1}-L_{\al_2})^2,(\Psi_\m)_r,(\Mfm)_r)=I_\m^r(L_{\al_1}-L_{\al_2},(L_{\al_1}-L_{\al_2})(\Psi_\m)_r,(\Mfm)_r)$$
$$=-4\int_0^{+\infty} \tau^2 \big[(L_{\al_1}-L_{\al_2}) (\Psi_m)_r(\tau,0^+) \dfrac{\partial}{\partial \tau}a(\tau,0)-\dfrac{\partial}{\partial \tau }[(L_{\al_1}-L_{\al_2})(\Psi_m)_r](\tau,0^+) a(\tau,0)\big] d\tau.$$

  Consid\' erons maintenant $D=L_\al-L_{\overline{\al}}$ et $(g,h)=(\Psi_2,D \Mf2)$ ou $(g,h)=(D\Psi_2,\Mf2)$. De m\^eme que pr\' ec\' edemment, on a
$$L_{\al}-L_{\overline{\al}}= -\im \fr{\partial\tau}\fr{\partial\te}-\im\frac{\te}{\tau^2+\te^2}\fr{\partial\tau}
    -\im \frac{\tau}{\tau^2+\te^2}\fr{\partial\te},$$
\no et  pour $A$ et $B$ deux fonctions de classe $\mathcal{C}^2$
sur un ouvert $U$ de $\R^2$ et pour $(\tau,\te)\in U$  avec
$\tau^2\neq\te^2$, on a la relation suivante
$$(\tau^2+\te^2)\left(A(L_{\al}-L_{\overline{\al}})(B)-A(L_{\al}-L_{\overline{\al}})(B)\right)(\tau,\te)$$
$$=
\frac{-\im}{2}\big[\dfrac{\partial}{\partial\te}\left((\tau^2+\te^2)R_\tau(A,B)\right)(\tau,\te)
+\dfrac{\partial}{\partial\tau}\left((\tau^2+\te^2)R_\te(A,B)\right)(\tau,\te)\big].$$
Les remarques \ref{supportMf} et \ref{regularitepsi} et les propri\'
et\' es de parit\' e des fonctions consid\' er\' ees assurent
l'int\' egrabilit\' e des deux fonctions
$\dfrac{\partial}{\partial\tau}\left((\tau^2+\te^2)R_\te(g,h)\right)(\tau,\te)$
et
$\dfrac{\partial}{\partial\te}\left((\tau^2+\te^2)R_\tau(g,h)\right)(\tau,\te)$
sur $\R^2$ ce qui permet d'obtenir
$$I_2(D,g,h)= 2\im\int_0^{+\infty} \tau^2 R_\tau(g,h) (\tau,0^+) d\tau+2\im \int_0^{\infty} \te^2 R_\te(g,h)(0^+,\te) d\te.$$
Pour $(\tau,\te)\in\R_+^2$, on a $\Mf2(\tau,\te)=a(\tau,-\te^2)+\te
b(\tau,-\te^2)$ avec $a,b\in \D(\R^*\times\R \cap\{(\tau,w); \tau^2>
w\})$.  En particulier, on a $R_\te(g,h)(0^+,\te)=0$ pour $\te>0$ et
par suite
$$I_2(D,g,h)= 2\im\int_0^{+\infty} \tau^2 R_\tau(g,h) (\tau,0^+) d\tau.$$
Par hypoth\`ese, il existe une fonction  $ \Psi_c\in Vect \langle
[A,B]; (A,B)\in
    \mathcal{S}ol(L_c,\la_1)\times \mathcal{S}ol(L_c,\la_2)\rangle $
    telle que,  pour $\tau>0$ et $\te>0$ alors $
     \Psi_2(\tau,\te)=\Psi_c\big((\tau+\im\te)^2,(\tau-\im\te)^2\big) .$ On a donc $ \Psi_2(\tau,0^+)=\dfrac{\partial}{\partial \tau}\Psi_2(\tau, 0^+)=0$ ce qui donne $R_\tau(\Psi_2, (L_{\al}-L_{\overline{\al}})\Mf2)(\tau,0^+)=0,$ et par suite
$$I_2((L_\al-L_{\overline{\al}})^2,\Psi_2,\Mf2)=I_2(L_\al-L_{\overline{\al}},(L_\al-L_{\overline{\al}})\Psi_2,\Mf2)$$
$$=2\im\int_0^\infty \tau^2\big[ (L_\al-L_{\overline{\al}})\Psi_2(\tau,0^+) \dfrac{\partial}{\partial\tau}a(\tau, 0) - \dfrac{\partial}{\partial\tau}[ (L_\al-L_{\overline{\al}})\Psi_2](\tau,0^+)  a(\tau,0)\big] d\tau.$$
Par hypoth\`ese, la fonction $\Psi_\m$ est donn\' ee par
$\Psi_\m(t_1,t_2)=(\Psi^++\im\Psi_c)(t_1,t_2)$  pour $t_1>t_2>0$
avec  $\Psi^+\in Vect\langle S^+(A,B); (A,B)\in
    \mathcal{S}ol(L,\la_1)\times \mathcal{S}ol(L,\la_2)\rangle$.
    Par ailleurs, on a les \' egalit\' es suivantes:
$$\begin{array}{c} (L_{\al_1}-L_{\al_2})(\Psi_\m)_r=\big( (L_1-L_2)\Psi_\m\big)_r=\big( (L_1-L_2)\Psi^+\big)_r+\im\big( (L_{1,c}-L_{2,c})\Psi_c\big)_r\\
 (L_\al-L_{\overline{\al}})\Psi_2(\tau,\te)=\big( (L_{1,c}-L_{2,c})\Psi_c\big)\big( (\tau+\im \te)^2,(\tau-\im\te)^2\big) \textrm { pour } (\tau,\te)\in\R^*_+\end{array}.$$
Ainsi, on obtient
 $$ (L_{\al_1}-L_{\al_2})(\Psi_\m)_r(\tau,0^+) =\im  (L_\al-L_{\overline{\al}})\Psi_2(\tau,0^+)$$
  $$\dfrac{\partial}{\partial\tau}\big[(L_{\al_1}-L_{\al_2})(\Psi_\m)_r\big](\tau,0^+)=\im\dfrac{\partial}{\partial\tau}\big[  (L_\al-L_{\overline{\al}})\Psi_2\big](\tau,0^+)  .$$
 Gr\^ace aux  expressions obtenues pr\' ec\' edemment de $I_\m^r((L_{\al_1}-L_{\al_2})^2,(\Psi_\m)_r,(\Mfm)_r)$ et $I_2((L_\al-L_{\overline{\al}})^2,\Psi_2,\Mf2)$, on en d\' eduit
 $$I_\m^r((L_{\al_1}-L_{\al_2})^2,(\Psi_\m)_r,(\Mfm)_r)+2I_2((L_\al-L_{\overline{\al}})^2,\Psi_2,\Mf2)=0.$$\end{dem}

\subsubsection{ Recollement sur $H\cdot\varpi(\ztp)$}

 En utilisant  l'isomorphisme $\varpi$ et les r\' esultats du paragraphe pr\' ec\' edent, nous \'etudions  \`a pr\' esent les conditions n\'ecessaires et
suffisantes pour que les fonctions  $\Psi_\m$ et $\Psi_2$
v\'erifient les relations $({\mathcal Rec})$ pour tout $f\in
\mathcal{D}(H\cdot\varpi(\ztp))$.\me

Pour $P\in\C[Q,S]$ et $f\in\D(\mathcal{U})$, on a de mani\`ere imm\'
ediate

$$\langle \partial(P) T_{\Phi},f\rangle-\langle T_{\partial(P)
\Phi},f\rangle=\int_\q[\Phi(X) \partial(P)f(X)-
\partial(P)\Phi(X)f(X)]dX$$
$$=\int_\q[\Phi\circ\varpi(X)
( \partial(P)f)\circ\varpi(X)-(
\partial(P)\Phi)\circ\varpi(X)f\circ\varpi(X)]dX.$$

Par ailleurs, pour toute application $g$ de classe
$\mathcal{C}^\infty$ sur un ouvert de $\q$, on a
$$\partial(Q)(g\circ\varpi)=-\partial(Q)(g)\circ\varpi$$
et
$$\partial(S)(g\circ\varpi)=\partial(S)(g)\circ\varpi.$$

Ainsi, comme $\Phi$ est solution du syst\`eme $(E_\chi)$ sur
$\q^{reg}$,  la fonction $\Phi\circ \varpi$ est solution du
syst\`eme  $(E_{\chi^-})$ sur $\q^{reg}$ o\`u $\chi^-$ est le
caract\`ere de $\C[Q,S_0]$ d\' etermin\' e par
$\chi^-(Q)=-\chi(Q)=-(\la_1+\la_2)$ et
$\chi^-(S_0)=\chi(S_0)=(\la_1-\la_2)^2$.

D'autre part, comme  $|\de|$ est $\varpi$-invariant, on a
$|\de|(\Phi\circ\varpi)=(|\de|\Phi)\circ\varpi=\widetilde{\Psi}\circ\varpi$.

On remarque que les fonctions $u$ et $v$ caract\' erisant les
orbites semi-simples (avec $u+v=Q$ et $uv=S$)  satisfont les
relations suivantes:
$$\begin{array}{l} u\circ\varpi (X)=-v(X)  \textrm{ pour }  X\in\mq^{reg}\\
\left\{ \begin{array}{l} u\circ\varpi
(X_{\tau,\te})=u(X_{\te,\tau})\\
v\circ\varpi (X_{\tau,\te})=v(X_{\te,\tau})\end{array}\right.
\textrm{ pour } X_{\tau,\te}\in \a_2^{reg}\end{array}.$$ Ainsi, on
obtient
$$\left\{ \begin{array}{l}  \tilde{\Psi}\circ\varpi
(X)=\Psi_\m(-v(X),-u(X))=\Psi_\m(-u(X),-v(X))  \textrm{ pour }  X\in\mq^{reg}\\
 \tilde{\Psi}\circ\varpi
(X_{\tau,\te})=\tilde{\Psi}(X_{\te,\tau})=\Psi_2(\te,\tau)=\Psi_2(\pm\te,\pm\tau),$\textrm{
pour }
 $X_{\tau,\te}\in\a_2. \end{array}\right.$$

On introduit alors les fonctions $\check{\Psi}_\m$ et
$\check{\Psi}_2$ d\' efinies par

$$\begin{array}{ll}\check{\Psi}_\m(t_1,t_2)=\Psi_\m(-t_1,-t_2)& \textrm{ pour } (t_1,t_2)\in\R^2 \textrm{  tel que }(t_1-t_2)t_1t_2\ne0\\
\check{\Psi}_2(\tau,\te)=\Psi_2(\te,\tau)&  \textrm{ pour }
(\tau,\te)\in(\R^*)^2\end{array}$$ de telle sorte que
$(|\de|\Phi\circ\varpi)_\m=\check{\Psi}_\m$ et
$(|\de|\Phi\circ\varpi)_2=\check{\Psi}_2.$\medskip

Par ailleurs, pour tout $f\in\D(\q)$, on a
$\mathcal{M}_H(f\circ\varpi)=\mathcal{M}_H(f)\circ\varpi$  et
$f\in\D(H\cdot\varpi(\ztp)$ si et seulement si
 $f\circ\varpi\in \D(H\cdot\ztp)$. On obtient alors le r\' esultat suivant:

 \begin{cor} Les  assertions suivantes sont \' equivalentes:
 \begin{enumerate}
\item La fonction $\Phi$, solution du syst\`eme $(E_\chi)$ sur $\q^{reg}$, v\' erifie  $\langle \partial(P) T_{\Phi},f\rangle-\langle T_{\partial(P)
\Phi},f\rangle=0$ pour  tout $f\in \mathcal{D}(H\cdot\varpi(\ztp))$
et  $P\in\C[Q,S_0]$,
\item  Les  fonctions $\Psi_\m$ et $\Psi_2$
v\'erifient les relations $({\mathcal Rec})$ pour tout $f\in
\mathcal{D}(H\cdot\varpi(\ztp)),$

\item Les  fonctions $\check{\Psi}_\m$ et $\check{\Psi}_2$
v\'erifient les relations $({\mathcal Rec})$ pour tout $f\in
\mathcal{D}(H\cdot\ztp).$
\end{enumerate}

\end{cor}

La proposition \ref{recolztp} permet alors d'obtenir le r\' esultat
suivant:

\begin{cor}\label{recolvarpiztp} Les trois assertions suivantes sont \' equivalentes:
  \begin{enumerate}
     \item Pour tout $f\in\D(H\cdot\varpi(\ztp))$, on a

    $$I_\m(L_1+L_2,\Psi_\m,\Mfm)+2I_2(L_\al+L_{\overline{\al}},\Psi_2,\Mf2)=0.$$

    \item Pour tout $\tau>0$, on a
        $$\left\{\begin{array}{c}\frac{\partial}{\partial\te}(\check{\Psi}_\m)_r(\tau,0^+)-
    \frac{\partial}{\partial\te}\check{\Psi}_2(\tau,0^+)=0\\
    \check{\Psi}_2(\tau,0^+)=0\;\;\qq\qq\end{array}\right..$$
    \item Pour tout $f\in\D(H\cdot\varpi(\ztp))$ et pour tout $\tau>0$, on a
     $$R_\te((\check{\Psi}_\m)_r,\big({\mathcal M}(f\circ\varpi)_\m\big)_r)(\tau,0^+)-R_\te(\check{\Psi}_2,{\mathcal M}(f\circ\varpi)_2)(\tau,0^+)=0.$$
    \end{enumerate}
\end{cor}

Comme la fonction $\Phi\circ\varpi$ est solution du syst\`eme
$(E_{\chi^-})$ sur $\q^{reg}$, on obtient alors gr\^ace au
corollaire \ref{solztp} et \`a la proposition \ref{recolS0ztp}:
\begin{cor} \label{solvarpiztp}
\begin{enumerate}
\item Les fonctions  $\Psi_\m$ et $\Psi_2$ satisfont l'une des
    conditions \' equivalentes du corollaire  \ref{recolvarpiztp} pr\' ec\' edent si et seulement si       il existe  $\Psi_{\chi^-}^+\in Vect\langle S^+(A,B); (A,B)\in
    \mathcal{S}ol(L,-\la_1)\times \mathcal{S}ol(L,-\la_2)\rangle$
 et $ \Psi_{c,\chi^-}\in Vect \langle [A,B]; (A,B)\in
    \mathcal{S}ol(L_c,-\la_1)\times \mathcal{S}ol(L_c,-\la_2)\rangle $
    telles que
 $$\begin{array}{l}\mbox{ pour } t_1>t_2>0 \textrm{ alors }
   \Psi_\m(-t_1,-t_2)= \check{\Psi}_\m(t_1,t_2)=\Psi_{\chi^-}^+(t_1,t_2)+\im \Psi_{c,\chi^-}(t_1,t_2) \\
    \mbox{ pour }\tau>0\mbox{ et }\te>0 \textrm{ alors }
   \Psi_2(\te,\tau)= \check{\Psi}_2(\tau,\te)= \Psi_{c,\chi^-}\big((\tau+\im\te)^2,(\tau-\im\te)^2\big) \end{array}.$$
\item Dans ce cas,  pour tout $f\in\D(H\cdot\varpi(\ztp))$, on a

    $$I_\m((L_1-L_2)^2,\Psi_\m,\Mfm)+2I_2((L_\al-L_{\overline{\al}})^2,\Psi_2,\Mf2)=0.$$
    \end{enumerate}
    \end{cor}

\subsubsection{Recollement sur $\mathcal{U}$}

La synth\`ese des r\' esultats des trois  paragraphes pr\' ec\'
edents nous permet de donner maintenant  les conditions
n\'ecessaires et suffisantes pour que les fonctions  $\Psi_\m$ et
$\Psi_2$ v\'erifient les relations $({\mathcal Rec})$ pour tout
$f\in \mathcal{D}(\mathcal{U})$. Nous gardons les notations
$\check{\Psi}_\m$ et $\check{\Psi}_2$ du paragraphe pr\'ec\'edent.\me

Commen\c cons par rechercher les fonctions $\Psi_2$ compatibles avec
les r\'esultats des corollaires \ref{solztp} et \ref{solvarpiztp}.

\begin{prop}\label{psi2}
La fonction $\Psi_2$ satisfait, pour tout $\tau\in\R^*$, les
relations
$$\left\{ \begin{array}{l} \Psi_2(\tau,0^+)=0\\ \check{ \Psi}_2(\tau,0^+)= \Psi_2(0^+,\tau)=0\end{array}\right.$$
si et seulement si il existe $\Psi_c\in Vect \langle
[\Phi_{\la_1},\Phi_{\la_2}], [\Phi_{\la_1},W_{\la_2}]+
[W_{\la_1},\Phi_{\la_2}]\rangle$ telle que, pour tout
$(\tau,\te)\in(\R_+^*)^2$, on ait
$ \Psi_2(\tau,\te)=\Psi_c\left((\tau+\im\te)^2,(\tau-\im\te)^2\right).$
\end{prop}

\begin{dem} Par le corollaire \ref{solztp}, la condition $
\Psi_2(\tau,0^+)=0$ implique qu'il existe $\Psi_c$ dans $ Vect
\langle [A,B]; (A,B)\in
    \mathcal{S}ol(L_c,\la_1)\times \mathcal{S}ol(L_c,\la_2)\rangle $ telle que , pour tout $(\tau,\te)\in(\R_+^*)^2$, on ait
$ \Psi_2(\tau,\te)=\Psi_c\left((\tau+\im\te)^2,(\tau-\im\te)^2\right).$

Ainsi, on veut  d\' eterminer les fonctions $\Psi_c$ de cette forme
telles que, pour tout $\tau>0$, on ait  $$\lim_{\te\rightarrow 0^+}
\Psi_c\big((\te+\im\tau)^2, (\te-\im\tau)^2\big)=0.$$

On remarque que lorsque   $A$ et $B$ sont des fonctions analytiques
sur $\C$, alors, pour $\tau>0$,  on a
$\dis \lim_{\te\rightarrow 0^+} [A,B]\big((\te+\im\tau)^2, (\te-\im\tau)^2\big)=0.$
Pour $\tau>0$, on a $\displaystyle\lim_{\te\rightarrow 0^+}
\log\big((\te+\im\tau))^2\big)=\log(\tau^2)+\im\pi$ et
$\displaystyle\lim_{\te\rightarrow 0^+}
\log\big((\te-\im\tau)^2\big)=\log(\tau^2)-\im\pi$, Ainsi, en posant
$z=(\te+\im\tau)^2$,  pour $\tau>0$, on obtient
$$\displaystyle\lim_{\te\rightarrow 0^+}[A, \log(\cdot) B](z, \overline{z})=-\lim_{\te\rightarrow 0^+}[ \log(\cdot)A, B](z, \overline{z})=-2i\pi A(-\tau^2) B(-\tau^2)$$
$$\textrm{ et } \hspace{3cm} \lim_{\te\rightarrow 0^+}[\log(\cdot)A, \log(\cdot) B](z, \overline{z})=0.$$

Ainsi, la fonction $\Psi_2$ paire par rapport \`a chaque variable,
d\' efinie sur $(\R^*_+)^2$ par
 $\Psi_2(\tau,\te)= [\Phi_{\la_1},\Phi_{\la_2}]\big((\tau+\im\te)^2, (\tau-\im\te)^2\big)$ v\' erifie donc la condition $\check{\Psi}_2(\tau,0^+)=0.$\me

Consid\' erons
$\Psi_c=\al_{12}[\Phi_{\la_1},W_{\la_2}]+\al_{21}[W_{\la_1},\Phi_{\la_2}]+\al_{22}[W_{\la_1},W_{\la_2}].$
On a donc
$$\Psi_c=\al_{12} ([\Phi_{\la_1},w_{\la_2}]+[\Phi_{\la_1}, \log(\cdot)\Phi_{\la_2}])+\al_{21} ([w_{\la_1},\Phi_{\la_2}]+[ \log(\cdot)\Phi_{\la_1},\Phi_{\la_2}])$$
$$+\al_{22}\big([w_{\la_1},w_{\la_2}]+[ \log(\cdot)\Phi_{\la_1},w_{\la_2}]+[w_{\la_1},\log(\cdot)\Phi_{\la_2}]+[ \log(\cdot)\Phi_{\la_1},\log(\cdot)\Phi_{\la_2}]\big)$$
et par suite
$$\check{\Psi}_2(\tau,0^+)=-2\im\pi(\al_{12}-\al_{21})\Phi_{\la_1}(-\tau^2)\Phi_{\la_2}(-\tau^2)$$
$$-
2\im\pi\al_{22}\big(w_{\la_1}(-\tau^2)\Phi_{\la_2}(-\tau^2)-w_{\la_2}(-\tau^2)\Phi_{\la_1}(-\tau^2)\big).$$

Comme  $\displaystyle
w_{\la_k}(t)=2\gamma+(1+2\gamma)\la_k t/4+o_{t\to0^-}(t)$ et
$\displaystyle \Phi_{\la_k}(t)=1+\la_k t/4+o_{t\to0^-}(t)$
la condition $\check{\Psi}_2(\tau,0^+)=0$  donne  $\al_{12}-\al_{21}=\al_{22}=0,$
et donc $\Psi_c=\al_{1,2}( [\Phi_{\la_1},W_{\la_2}]+ [W_{\la_1},\Phi_{\la_2}])$.

Finalement, les conditions $\Psi_2(\tau,0^+)=\check{\Psi}_2(\tau,
0^+)= 0$ implique $\Psi_c\in  Vect \langle
[\Phi_{\la_1},\Phi_{\la_2}], [\Phi_{\la_1},W_{\la_2}]+
[W_{\la_1},\Phi_{\la_2}]\rangle$. La r\' eciproque est imm\'
ediate.\end{dem}

\begin{rem} Pour $\la\in\C^*$,  les solutions fondamentales $\Phi_\la$ et $W_\la=w_\la+\log(\cdot)\Phi_\la$ de $\mathcal{S}ol(L_c,\la)$ v\' erifient, pour tout   $z\in\C-\R$,
$$\Phi_\la(-z)=\Phi_{-\la}(z)\quad \textrm{ et } \quad w_\la(-z)=w_{-\la}(z).$$

Par suite, pour $(z_1,z_2)\in(\C-\R)^2$, on a les relations
suivantes
$$[\Phi_{-\la_1},\Phi_{-\la_2}](z_1,z_2)= [\Phi_{\la_1},\Phi_{\la_2}](-z_1,-z_2),$$
$$\big( [\Phi_{-\la_1},W_{-\la_2}]+ [W_{-\la_1},\Phi_{-\la_2}]\big)(z_1,z_2)$$
$$= \big( [\Phi_{\la_1},w_{\la_2}]+ [w_{\la_1},\Phi_{\la_2}]\big)(-z_1,-z_2)+\log(z_1z_2) [\Phi_{\la_1},\Phi_{\la_2}](-z_1,-z_2).$$

Ainsi, pour $\Psi_c\in Vect \langle [\Phi_{\la_1},\Phi_{\la_2}],
[\Phi_{\la_1},W_{\la_2}]+ [W_{\la_1},\Phi_{\la_2}]\rangle$, la
fonction $\Psi_2$ paire en chaque variable d\' efinie sur
$(\R^*_+)^2$ par  $
\Psi_2(\tau,\te)=\Psi_c\left((\tau+\im\te)^2,(\tau-\im\te)^2\right)$
v\' erifie $$\check{\Psi}_2(\tau,\te)=\Psi_{c,\chi^-}\big(
(\tau+\im\te)^2, (\tau-\im\te)^2\big)$$ o\`u $\Psi_{c,\chi^-}\in
Vect \langle [A,B]; (A,B)\in
    \mathcal{S}ol(L_c,-\la_1)\times \mathcal{S}ol(L_c,-\la_2)\rangle$ satisfait   $\Psi_{c,\chi^-}(z_1,z_2)=-\Psi_c(-z_1,-z_2)$ sur $(\C-\R)^2.$
    Le r\' esultat du lemme pr\' ec\' edent est donc compatible avec celui du corollaire \ref{solvarpiztp}.
\end{rem}

\begin{Def} On d\' efinit la fonction ${\mathcal S}_{\la_1,\la_2}$   sur $(\C^*)^2$ par
$${\mathcal S}_{\la_1,\la_2}(z_1,z_2)=\big([\Phi_{\la_1},w_{\la_2}]+[w_{\la_1},\Phi_{\la_2}]\big)(z_1,z_2)+\log |z_1z_2|[\Phi_{\la_1},\Phi_{\la_2}](z_1,z_2)$$
 de telle sorte que
 \begin{enumerate} \item

${\mathcal S}_{\la_1,\la_2}
(z,\overline{z})=\big([\Phi_{\la_1},W_{\la_2}]+
[W_{\la_1},\Phi_{\la_2}]\big)(z,\overline{z})$, pour tout $z\in
\C-\R$.

\item ${\mathcal S}_{\la_1,\la_2}(t_1,t_2)=\big( [\Phi_{\la_1},W^r_{\la_2}]+ [W^r_{\la_1},\Phi_{\la_2}]\big)(t_1,t_2)$, pour tout $(t_1,t_2)\in(\R^*)^2$.
\end{enumerate}
\end{Def}
\begin{theo}\label{recU} Les fonctions $\Psi_\m$ et $\Psi_2$ v\'erifient les conditions $(\mathcal{R}ec)$ pour tout $f\in\D({\mathcal U})$ si et seulement si  il existe
\begin{enumerate}\item $\Psi^+\in Vect\langle S^+(A,B); (A,B)\in
    \mathcal{S}ol(L,\la_1)\times \mathcal{S}ol(L,\la_2)\rangle$
  \item  $ \Psi^-\in Vect\langle [\Phi_{\la_1},\Phi_{\la_2}],
{\mathcal S}_{\la_1,\la_2}\rangle$\end{enumerate}
    telles que,
pour tout $(t_1,t_2)\in(\R^*)^2-diag$ et $(\tau,\te)\in(\R_+^*)^2$,
on ait les relations
$$\begin{array}{c} \Psi_\m(t_1,t_2)=
\Psi^+(t_1,t_2)+\im signe(t_1-t_2)\Psi^-(t_1,t_2)\\
\textrm{et } \;\;\;  \Psi_2(\tau,\te)=\Psi^-\left((\tau+\im\te)^2,(\tau-\im\te)^2\right).\end{array}$$
\end{theo}

\begin{dem} On suppose tout d'abord que $\Psi_\m$ et $\Psi_2$
v\'erifient les conditions $(\mathcal{R}ec)$ pour tout
$f\in\D({\mathcal U})$. Par la proposition  \ref{solztp} et le
corollaire \ref{recolvarpiztp}, on a $\Psi_2(\tau,0^+)=
\Psi_2(0^+,\tau)=0$ pour tout $\tau\in\R^*$. La proposition
\ref{psi2} donne alors l'existence de $ \Psi^-\in Vect\langle
[\Phi_{\la_1},\Phi_{\la_2}], {\mathcal S}_{\la_1,\la_2}\rangle$
telle que pour tout $(\tau,\te)\in(\R^*)^2$, on ait $
\Psi_2(\tau,\te)=\Psi^-\left((\tau+\im\te)^2,(\tau-\im\te)^2\right).$
Par les corollaires \ref{recolm} et  \ref{solztp} et les propri\'
et\' es de sym\' etries de $\Psi_\m$, il existe $\Psi^+\in
Vect\langle S^+(A,B); (A,B)\in
    \mathcal{S}ol(L,\la_1)\times \mathcal{S}ol(L,\la_2)\rangle$ telle que, pour tout $(t_1,t_2)\in(\R^*)^2-diag$, on ait  $\Psi_\m(t_1,t_2)=
\Psi^+(t_1,t_2)+\im signe(t_1-t_2)\Psi^-(t_1,t_2)$. Maintenant les propri\' et\' es des  fonctions $[\Phi_{\la_1},\Phi_{\la_2}]$ et ${\mathcal S}_{\la_1,\la_2}$ montrent que les fonctions $\Psi_\m$ et $\Psi_2$ satisfont alors les propri\' et\' es du corollaire \ref{solvarpiztp}.\me

R\' eciproquement, on suppose que les fonctions $\Psi_\m$ et
$\Psi_2$ sont donn\' ees par $\Psi_\m(t_1,t_2)= \Psi^+(t_1,t_2)+\im
signe(t_1-t_2)\Psi^-(t_1,t_2)$, pour tout
$(t_1,t_2)\in(\R^*)^2-diag$ et  par
$\Psi_2(\tau,\te)=\Psi^-\left((\tau+\im\te)^2,(\tau-\im\te)^2\right)$,
pour tout $(\tau,\te)\in(\R_+^*)^2$ o\`u $\Psi^+\in Vect\langle
S^+(A,B); (A,B)\in
    \mathcal{S}ol(L,\la_1)\times \mathcal{S}ol(L,\la_2)\rangle$
    et $ \Psi^-\in Vect\langle [\Phi_{\la_1},\Phi_{\la_2}],
{\mathcal S}_{\la_1,\la_2}\rangle$. Ainsi, les r\' esultats des paragraphes pr\' ec\' edents assurent que les fonctions $\Psi_\m$ et $\Psi_2$ satisfont les relations $({\mathcal R}ec)$ pour toute fonction $g$ dont le support est contenu soit dans $H\cdot\p\mq$ (corollaire \ref{expsolm} et proposition \ref{exImmoins}), soit dans $H\cdot\ztp$ (corollaire \ref{solztp} et proposition \ref{recolS0ztp}), ou encore dans $H\cdot\vi(\ztp)$ (corollaire \ref{solvarpiztp}).\me

Soit $f\in\D({\mathcal U})$. Comme ${\mathcal U}=H\cdot\p\mq\cup
H\cdot\ztp\cup H\cdot\vi(\ztp)$, le th\'eor\`eme de partition de
l'unit\'e assure l'existence de  fonctions $\varphi_\m,\;\varphi_3$
et $\varphi^-_3$ de classe $\mathcal{C}^\infty$ sur $\q$ telles que:
(a)  $ Support(\varphi_\m)\subset H\cdot\p\mq, \quad
Support(\varphi_3)\subset H\cdot\ztp,\quad \textrm{ et } \quad
Support(\varphi^-_3)\subset H\cdot\vi(\ztp),$

(b) $\textbf{1}_{supp(f)}\leq \varphi_\m+\varphi_3+\varphi^-_3\leq 1.$

\noindent Ainsi, on a  $f=\varphi_\m f+\varphi_3f+\varphi^-_3f$ et par suite  $\Psi_\m$ et $\Psi_2$ satisfont les relations $({\mathcal R}ec)$ pour la fonction $f$ par ce qui pr\' ec\`ede. \end{dem}

\begin{rem} L'expression de la fonction $\Psi_2$ d\' ependait jusqu'\`a pr\' esent  du choix de la d\' etermination du logarithme choisie pour d\' ecrire le syst\`eme fondamental de solutions de ${\mathcal S}ol_{L_c, \la}$  \`a travers la fonction $\Psi_c$. L'expression obtenue pour $\Psi_2$ dans le th\' eor\`eme ci-dessus en terme de la fonction ${\mathcal S}_{\la_1,\la_2}$ montre qu'elle est ind\' ependante de ce choix.
\end{rem}

Soient $\Psi_\m$ et $\Psi_2$ satisfaisant les hypoth\`eses du th\'
eor\`eme pr\' ec\' edent. Les fonctions $\tilde{\Psi}$ et
$\Phi=\tilde{\Psi}/|\delta| $ s'expriment
alors de la mani\`ere suivante en terme des fonctions $u$ et $v$,
qui  caract\' erisent les orbites semi-simples. On rappelle que ces
fonctions satisfont les relations suivantes:
 $ u(X)+v(X)=Q(X)=\frac{1}{2} tr(X^2)$, $ u(X)v(X)=S(X)=det(X)$ et $u(X)-v(X)=\delta(X).$ Comme pour $X\in\mq$, on a $u(X)-v(X)=\de(X)=|\de(X)|\geq 0$ et  pour $X\in\a_2$, on a $|\de|(X)=-\im\de(X)$  et  $u(X)=\overline{v(X)}$, nous obtenons:

$$\tilde{\Psi}(X)=\left\{\begin{array}{ll} \Psi^+(u(X),v(X))+\im \Psi^-(u(X),v(X))& \mbox{ si }
X\in\mq^{reg}\\
\Psi^-(u(X),v(X))&\mbox{ si } X\in\a_2^+  \end{array}\right.,$$

Par ailleurs, pour tout $X\in\q^{reg}$, on a
$$\mathcal{S}_{\la_1,\la_2}(u(X),v(X))=
\big([\Phi_{\la_1},w_{\la_2}]+
[w_{\la_1},\Phi_{\la_2}]\big)(u(X),v(X))$$
$$+(\log|S(X)|)
[\Phi_{\la_1},\Phi_{\la_2}](u(X),v(X)).$$

On introduit alors les fonctions suivantes

$$F_{ana}:=\frac{[\Phi_{\la_1},\Phi_{\la_2}](u,v)}{u-v},$$

$$F_{sing}:=\frac{[\Phi_{\la_1},w_{\la_2}](u,v)+[w_{\la_1},\Phi_{\la_2}](u,v)}{u-v}+
\log|uv|\frac{[\Phi_{\la_1},\Phi_{\la_2}](u,v)}{u-v}$$ 

$$\textrm{ et }\hspace{2cm}F^+_{A,B}:= Y(S_0)\frac{ S^+(A,B)(u,v)}{u-v},\;\;\;\;$$ avec
$(A,B)\in\{(\Phi_{\la_1},\Phi_{\la_2}),\;(\Phi_{\la_1},W^r_{\la_2}),\;(W^r_{\la_1},\Phi_{\la_2}),\;(W^r_{\la_1},W^r_{\la_2})\}$
( le polyn\^ome $S_0=Q^2-4S$ est \`a valeurs r\'
eelles sur $\q$ et $S_0(X)\geq 0$ si et seulement si $ X\in
H\cdot\m\cap\q.$)

\begin{theo}\label{solgen}
\begin{enumerate} \item Les fonctions $F_{ana},\;F_{sing}$ et $F^+_{A,B}$,
avec
$(A,B)\in\{(\Phi_{\la_1},\Phi_{\la_2}),\;(\Phi_{\la_1},W^r_{\la_2}),\;(W^r_{\la_1},\Phi_{\la_2}),\;(W^r_{\la_1},W^r_{\la_2})\}$
sont localement int\' egrables sur ${\mathcal U}$.

\item Ces fonctions forment une  base des distributions propres invariantes  sur
$\mathcal{U}$ donn\' ees par une fonction de $L^1_{loc}(\mathcal{U})^H$.\\
\end{enumerate}
\end{theo}

\begin{dem} Soit $F$ l'une des fonctions $F_{ana},\;F_{sing}$ ou
$F^+_{A,B}$.  Il suffit de montrer que,  pour toute fonction $f$
positive de $\mathcal{D}(\mathcal{U})$, l'int\'egrale $\int_\q
f\left|F\right|$ est finie. En notant,  comme pr\' ec\' edemment
$\tilde{\Psi}=|\delta| F$, ceci revient \`a montrer gr\^ace \`a la
formule d'int\'egration de Weyl que, pour toute fonction $f$
positive de $\mathcal{D}(\mathcal{U})$, la somme
$\dis \Sigma(f)=\int_{t_1>t_2}|\Psi_\m|(t_1,t_2)\Mfm(t_1,t_2) dt_1
dt_2+8\int_{\tau>0\;\; \te>0}
|\Psi_2|(\tau,\te) \Mf2(\tau,\te) (\tau^2+\te^2)d\te d\tau $ est
finie.\me

Par partition de l'unit\' e, il existe des fonctions $f_\m,\;f_3$ et
$g_3$ de classe $\mathcal{C}^\infty$ sur $\q$ de support inclus
respectivement dans $H\cdot\p\mq$, $H\cdot\ztp$ et $H\cdot\vi(\ztp)$
telles que $f=f_\m+f_3+g_3$. La convergence des int\' egrales
$\Sigma(g)$ pour $g=f_\m$, $f_3$ ou $g_3$ d\' ecoule des remarques
\ref{continue} , \ref{supportMf} et \ref{regularitepsi} (voir
preuves des lemmes \ref{ippqm} et \ref{ippmr}) ce qui donne la
premi\`ere assertion.

La deuxi\`eme assertion est alors une cons\' equence imm\' ediate du th\' eor\`eme \ref{recU}. \end{dem}
\subsection{Perspectives sur les distributions propres invariantes
de $\Lu$}

 Nous allons tout  d'abord un lemme  pr\'eciser
la r\'egularit\'e des fonctions trouv\'ees dans le th\'eor\`eme
\ref{solgen}.

\begin{lem} Les fonctions  $F_{ana}$ et $F_{sing}$ s'expriment sous la forme suivante:

$$F_{ana}=f_{ana}(Q,S_0) \quad \textrm{ et } \quad F_{sing}=f_s(Q,S_0)+\log|S|g_s(Q,S_0)$$ o\`u les fonctions $f_{ana}$, $f_s$ et $g_s$ sont holomorphes sur $\C$.
\end{lem}

\begin{dem}  Soient $f$ et $g$ deux fonctions holomorphes sur  $\C$
dont la  s\'erie enti\`ere en z\'ero est de  rayon de convergence
infini. Alors, pour $(x,h)\in\C^2$, on a
$$f(x+h)g(x-h)-f(x-h)g(x+h)=h\int_{-1}^1(f'(x+th)g(x-th)-f(x+th)g'(x-th))dt.$$
Ainsi, la fonction $(x,h)\mapsto\frac{f(x+h)g(x-h)-f(x-h)g(x+h)}{h}$
est holomorphe  sur $\C^2$. D'autre part, cette fonction est paire
par rapport \`a $h$, donc elle admet un d\' eveloppement  en s\'erie
enti\`ere sur $\C^2$ de la forme $\sum_{0\leq
k,l}a_{k,l}x^kh^{2l}$. \me

Maintenant, pour $X\in\q$, on consid\`ere  $x=\frac{Q}{2}(X)$ et
$h=\frac{\de}{2}(X)$ de telle sorte que $u(X)=x+h$ et $v(X)=x-h$. On
obtient alors
$$\frac{[f,g](u,v)}{\de}=\sum_{0\leq
k,l}\frac{a_{k,l}}{2^{k+2l}}Q^k\de^{2l}=\sum_{0\leq
k,l}\frac{a_{k,l}}{2^{k+2l}}Q^k S_0^{l}. $$

Les d\' efinitions de $F_{ana}$ et $F_{sing}$ donnent alors le r\'
esultat voulu .\end{dem}

\begin{rem} La fonction $\delta=u-v$ correspond \`a un choix de racine carr\' e du polyn\^ome $S_0=Q^2-4S$. Le lemme pr\' ec\' edent exprime les fonctions $F_{ana}$ et $F_{sing}$ uniquement en terme des polyn\^omes $Q$ et $S$. Ces fonctions sont donc ind\' ependantes du choix de cette racine carr\' e.
\end{rem}

\begin{cor} L'espace des distributions propres invariantes sur $\q$ d\' efinies par une fonction analytique sur $\q$ est engendr\' e par $F_{ana}$.
\end{cor}

  \begin{dem} Comme les op\' erateurs $\partial(Q)$ et $\partial(S)$
sont polynomiaux, une simple int\' egration par parties sur $\q$
assure que $F_{ana}$ d\' efinit bien une distribution propre
invariante sur $\q$. Le r\' esultat voulu est alors imm\' ediat.
\end{dem}

\begin{prop} La fonction
$F_{sing}$ est dans $\Lu$.\end{prop}

\begin{dem} Il suffit de montrer que $\log|S|\in\Lu$ ou encore  que pour toute fonction positive $\chi$
dans $\D(\q)$, l'int\'egrale $\int_\q\big|\log|S|\big|(X) \chi(X) dX$ est
finie.\me

On a la d\'ecomposition $\q=\q_1\oplus\q_2,$ o\`u
$\dis \q_1=\left\{\left (\begin{array}{c|c} 0&A\\ \hline
0&0\end{array}\right );\;A\in M_{2}(\R)\right\}
\;\;\textrm{et}\;\;\q_2=\left\{\left (\begin{array}{c|c} 0&0\\
\hline B&0\end{array}\right );\;B\in M_{2}(\R)\right\}.$
Soient $\mathcal{Q}_1$ (respectivement $\mathcal{Q}_2$) la forme
quadratique d\'efinie sur $\q_1$ (respectivement  $\q_2$) par
$$\mathcal{Q}_1\left(\left (\begin{array}{c|c} 0&A\\ \hline
0&0\end{array}\right)\right)=\det(A)\;\;\left(\textrm{respectivement}\;\;
\mathcal{Q}_2\left(\left (\begin{array}{c|c} 0&0\\ \hline
B&0\end{array}\right )\right)=\det(B)\right).$$
Ainsi on a
$$
\int_\q\big|\log|S(X)|\big|\chi(X)dX=\int_{\q_1\times\q_2}\big|\log|\mathcal{Q}_1(X_1)|+\log|\mathcal{Q}_2(X_2)|\big|\chi(X_1+X_2)dX_1dX_2$$
$$\leq\int_{\q_1\times\q_2}\left(\big|\log|\mathcal{Q}_1(X_1)|\big|+\big|\log|\mathcal{Q}_2(X_2)|\big|\right)\chi(X_1+X_2)dX_1dX_2$$$$=
\int_{\R^2}\left(\big|\log|t_1|\big|+\big|\log|t_2|\big|\right)M_{\mathcal{Q}_1,\mathcal{Q}_2}\chi(t_1,t_2)dt_1dt_2,$$
gr\^ace \`a la relation \ref{rel2} et au lemme \ref{double}.
Comme les formes quadratiques $\mathcal{Q}_1$ et $\mathcal{Q}_2$
sont de signature $(2,2)$ , alors le th\'eor\`eme \ref{comp2} permet
d'affirmer que $M_{\mathcal{Q}_1,\mathcal{Q}_2}\chi(t_1,t_2)$
s'\'ecrit sous la forme
$\chi(t_1,t_2)=a(t_1,t_2)+|t_1|b(t_1,t_2)+|t_2|c(t_1,t_2)+|t_1 t_2|d(t_1,t_2),$
o\`u $a,\;b,\;c$ et $d$ sont dans $\D(\R^2)$. Ainsi la derni\`ere
int\'egrale est convergente. D'o\`u le r\'esultat.\end{dem}

La m\' ethode de descente utilis\' ee pour l'\' etude de l'int\'
egrale orbitale ${\mathcal M}_H(f)$ d'une fonction  $f$ de $\D(\q)$
ne permet pas de  d\' ecrire le comportement de ${\mathcal M}_H(f)$
au voisinage de $0$.

Ainsi,  il ne nous est pas possible de  pr\' eciser
l'int\'egrabilit\'e sur $\q$ ou non des fonctions de la forme
$F_{A,B}^+$ avec
$(A,B)\in\{(\Phi_{\la_1},\Phi_{\la_2}),\;(\Phi_{\la_1},W^r_{\la_2}),\;(W^r_{\la_1},\Phi_{\la_2}),\;(W^r_{\la_1},W^r_{\la_2})\}$.

Pour la m\^eme raison, les int\' egrations par parties effectu\' ees
dans le paragraphe  \ref{CondRecol} ne sont pas forc\' ement licites pour
$f\in\D(\q)$. Par suite,   bien que la fonction $F_{sing}$ soit
localement int\' egrable sur $\q$, nous ne pouvons pas affirmer que
cette fonction d\' efinit une distribution propre invariante sur
$\q$ tout entier.

\end{document}